\documentclass{article}

\usepackage[utf8]{inputenc}

\usepackage{algorithm}
\usepackage{algorithmic}
\usepackage{amsmath}
\usepackage{amssymb}  
\usepackage{amsthm}
\usepackage{arydshln}
\usepackage{authblk}

\usepackage{bbold}
\usepackage{bm}
\usepackage{booktabs}
\usepackage{bookmark}
\usepackage{calc}              
\usepackage{cleveref}
\usepackage{enumitem}
\usepackage{extarrows}
\usepackage{float}
\usepackage{geometry}
\usepackage{graphicx}          
\usepackage{hyperref}
\usepackage{indentfirst}
\usepackage{latexsym}
\usepackage{lineno}
\usepackage{mathtools}
\usepackage{multirow}          
\usepackage{subcaption}
\usepackage{tikz}              
\usepackage{xcolor}            
\usetikzlibrary{shapes.geometric, arrows, chains}

\title{A hierarchical approach for multicontinuum homogenization in high contrast media}

\providecommand{\keywords}[1]{\textbf{\textit{Keywords:}} #1}

\author[a,b]{Wei Xie}
\author[b]{Viet Ha Hoang}
\author[a]{Yin Yang\thanks{Corresponding author}}
\author[a]{Yunqing Huang}
\affil[a]{School of Mathematics and Computational Science, Xiangtan University, National Center for Applied Mathematics in Hunan, Xiangtan, Hunan, China 411105}
\affil[b]{Division of Mathematical Sciences, School of Physical and Mathematical Sciences, Nanyang Technological University, Singapore 637371}

\date{}

\begin{document}
\maketitle

\renewcommand{\thefootnote}{}
\footnotetext{E-mail addresses: xiew@smail.xtu.edu.cn (Wei Xie), 
vhhoang@ntu.edu.sg (Viet Ha Hoang),
yangyinxtu@xtu.edu.cn (Yin Yang),
huangyq@xtu.edu.cn (Yunqing Huang)}
\renewcommand{\thefootnote}{\arabic{footnote}}

\begin{abstract}
A recently developed upscaling technique, the multicontinuum homogenization method, has gained significant attention for its effectiveness in modeling complex multiscale systems. 
This method defines multiple continua based on distinct physical properties and solves a series of constrained cell problems to capture localized information for each continuum. 
However, solving all these cell problems on very fine grids at every macroscopic point is computationally expensive, which is a common limitation of most homogenization approaches for non-periodic problems. 
To address this challenge, we propose a hierarchical multicontinuum homogenization framework. The core idea is to define hierarchical macroscopic points and solve the constrained problems on grids of varying resolutions. 
The local solutions are decomposed into the linear interpolation of contributions inherited from preceding levels and an additional correction term.
This combination is substituted into the original constrained problems, and the correction term is resolved using finite element (FE) grids of varying sizes, depending on the level of the macropoint.
By normalizing the computational cost of fully resolving the local problem to $\mathcal{O}(1)$,
we establish that our approach incurs a cost of $\mathcal{O}(L \eta^{(1-L)d})$, highlighting substantial computational savings across hierarchical layers $L$, coarsening factor $\eta$, and spatial dimension $d$. Numerical experiments validate the effectiveness of the proposed method in media with slowly varying properties, underscoring its potential for efficient multiscale modeling.
\end{abstract}


\section{Introduction}
Multiscale phenomena are common in engineering and industrial applications, 
involving processes that interact across different scales. 
For example, in fluid flow through porous media, features like fractures, vugs, and micropores are crucial for determining properties such as permeability and fluid transport, which are critical to subsurface flow modeling and material design. 
In composite materials, the arrangement of fibers and matrix influences heat conduction, affecting the thermal behavior of engineering materials. 
Similarly, in electromagnetic wave propagation, the structure of materials like layered composites or metamaterials affects the dispersion and absorption of waves. 
However, traditional numerical methods, such as the finite element method and the finite volume method, often have difficulty resolving the fine-scale details of these multiscale materials due to high computational costs and complexity. 
Developing an efficient multiscale method is therefore necessary.

Among the various approaches developed for multiscale problems, an important class focuses on constructing local basis functions on fine grids and solving the resulting algebraic systems on coarse grids. This approach balances computational efficiency with the preservation of fine-scale details. Key methods in this category include the localized orthogonal decomposition method \cite{maalqvist2014localization,henning2014localized}, which generates localized multiscale basis functions through systematic decomposition of solution spaces, effectively capturing fine-scale information while reducing computational cost. The multiscale finite element method (MsFEM) \cite{hou1997multiscale,hou1999convergence,ye2020convergence,jiang2017reduced,huang2001partition} employs precomputed basis functions to account for fine-scale heterogeneities, enabling accurate coarse-scale simulations. 
The generalized multiscale finite element method (GMsFEM) \cite{efendiev2013generalized,calo2016randomized,alikhanov2025multiscale,xie2025time,xie2025multiscale,chung2014adaptive,ye2024robust,fu2024efficient} extends the multiscale finite element method (MsFEM) by introducing symmetric and adaptive basis function construction. 
The constraint energy minimizing generalized multiscale finite element method (CEM-GMsFEM) \cite{chung2018constraint,ye2023constraint,wang2024multiscale,xie2024cem} further constructs local basis functions by solving energy minimization problems within oversampled domains, and exhibits convergence with respect to the coarse mesh size.

In contrast, some problems necessitate a focus solely on macroscopic variables, circumventing the need to resolve fine-scale details explicitly. Upscaling techniques establish a connection between microscopic and macroscopic scales by assuming that microscale solutions can be represented using different ansatzes of the macroscale solutions.
Examples include homogenization \cite{hornung2012homogenization, cioranescu1999introduction} and the heterogeneous multiscale method \cite{abdulle2012heterogeneous, henning2009heterogeneous}. 
In this paper, we focus on the multicontinuum homogenization method \cite{chung2024multicontinuum, efendiev2023multicontinuum, xie2025multicontinuum, efendiev2025multicontinuum, bai2026multicontinuum}. 
This method defines multiple continua to represent distinct physical regions. 
The macroscopic coefficient is obtained by solving a series of constraint problems in an oversampled representative volume element (or target coarse block) region. The constraints in the cell problems are used to represent different behaviors for different continua, which has been widely investigated in nonlocal multicontinuum upscaling (NLMC) methods \cite{chung2018non, chung2021nonlocal}, as well as in the CEM-GMsFEM \cite{chung2018constraint}. 
In fact, the accuracy of the upscaling coefficient depends heavily on the finite element (FE) grid used to solve the cell problems. 
In practice, since local geometries or material properties are typically non-periodic, resolving all cell problems on very fine grids is computationally expensive.
This motivates the development of a hierarchical strategy for computing the upscaling coefficient.

In this work, we propose a hierarchical multicontinuum homogenization framework. The method consists of three main components. First, the macropoints (centers of coarse blocks) are organized into a hierarchy of grids. Next, distinct finite element spaces are assigned to different levels. Finally, the local solutions are constructed inductively, with each solution expressed as a combination of a linear interpolation from macropoints at preceding levels and a correction term. Substituting this representation into the original constrained problem yields a reduced problem for the correction term, which becomes the main computational task. Since the correction problems are solved in nested FE spaces, the overall cost is significantly lower than solving all problems on the same fine grid. We also provide a detailed cost comparison between the original and hierarchical formulations. Numerical experiments in two and three dimensions confirm the efficiency and accuracy of the proposed framework: when physical properties vary moderately, the hierarchical approach achieves high accuracy, even in cases with dramatic variations, the results remain acceptable.

Over the past decades, hierarchical ideas \cite{park2020hierarchical, brown2013efficient, brown2017hierarchical} and multilevel methods \cite{briggs2000multigrid, xu2017algebraic, cheng2022scalable, efendiev2015multilevelcicp, li2018multilevel, ye2025highly, efendiev2015multilevelmms} have been extensively developed and widely applied in large-scale scientific computing, including reservoir simulation and porous media flow. Most of these approaches, however, either explicitly resolve fine-grid solutions or rely on scale-separation assumptions to construct efficient macroscopic models or to establish convergence analyses. In contrast, the main contribution of this work is an efficient upscaling technique that bypasses direct fine-grid computation while producing a computable and accurate macroscopic model, without assuming scale separation or periodicity of the permeability field. 
Our approach involves only two levels. The first level is used to construct the computable macroscopic model, while the second level is employed to compute an efficient macroscopic coefficient. To reduce the computational cost at the second level, a hierarchical algorithm is applied. 
To the best of our knowledge, this is also the first work to present three-dimensional numerical examples for the multicontinuum homogenization method.

This paper is organized as follows.
In Section \ref{sec:preliminaries}, we present the model problem and review the multicontinuum homogenization method.
In Section \ref{sec:hmh}, we detail the hierarchical computation of the upscaling coefficient and discuss the computational savings.
The numerical results are presented in Section \ref{sec:numerres}.
Finally, the paper concludes with a summary in Section \ref{sec:conclusions}.

\section{Preliminaries} \label{sec:preliminaries}
In this section, we present the model problem and review the multicontinuum homogenization framework introduced in \cite{efendiev2023multicontinuum,chung2024multicontinuum}.
Let $\Omega \subset \mathbb{R}^d$ ($d=2,3$) be a convex region. 
We consider the following equation:
\begin{equation}
\mathcal{A}^{\epsilon} u = f 
\quad \textnormal{in}~\Omega,
\label{eq:Au_pde}
\end{equation}
where $\mathcal{A}^{\epsilon}$ is a linear differential operator, and $u$ vanishes on the boundary.
For simplicity, we consider the scalar elliptic equation: 
$\mathcal{A}^{\epsilon}(u) = - \nabla \cdot (\kappa \nabla u)$. 
The weak formulation of Eq. \eqref{eq:Au_pde} is
\begin{equation}
\int_{\Omega} \kappa \nabla u \nabla v = 
\int_{\Omega} f v, \quad
\forall v \in V,
\label{eq:weak_galerkin}
\end{equation}
where $V$ is $H_0^1(\Omega)$. 
We use $h$ to denote the mesh size required to resolve the heterogeneities induced by the coefficients or geometries and to achieve high accuracy. The fine-scale finite element partition is assumed to be general and may consist of either structured or unstructured meshes, provided it forms a shape-regular family. For clarity of presentation, we illustrate the method using rectangular elements.
To avoid the unaffordable algebra system from that mesh size, we use the multicontinuum homogenization method to build an upscaling model and solve it at the macroscopic scale. 

\begin{figure}[h!]
\centering
\begin{tikzpicture}
\draw[help lines,step=0.5cm] (0,0) grid (6,5);
\foreach \i in {0,...,11} {
\foreach \j in {0,...,9} {
\node at (\i*0.5+0.25,\j*0.5+0.2) {*};
}}
\draw[ultra thick,->] (4.5,1.25) -- (7.8,2.5);
\draw[->] (9.8,2.3) -- (11.8,1.5);
\draw[->] (10.5,2.9) -- (11.8,2.5);
\draw[->] (11,3.5) -- (11.8,3.5);
\draw[ultra thick] (4,1) rectangle (4.5,1.5);
\draw (8,1) rectangle (11,4);
\draw (8.5,1.5) rectangle (10.5,3.5);
\draw (9.2,2.2) rectangle (9.8,2.8);
\node[anchor=south] at (6.5,0) {$\Omega$};
\node at (12.2,3.5) {$K_p$};
\node at (12.2,2.5) {$R_p^+$};
\node at (12.2,1.5) {$R_p$};
\end{tikzpicture}
\caption{Illustration of RVE $R_p$, oversampling domain $R_p^+$ and coarse block $K_p$.}
\label{fig:macropic_point_rve}
\end{figure}
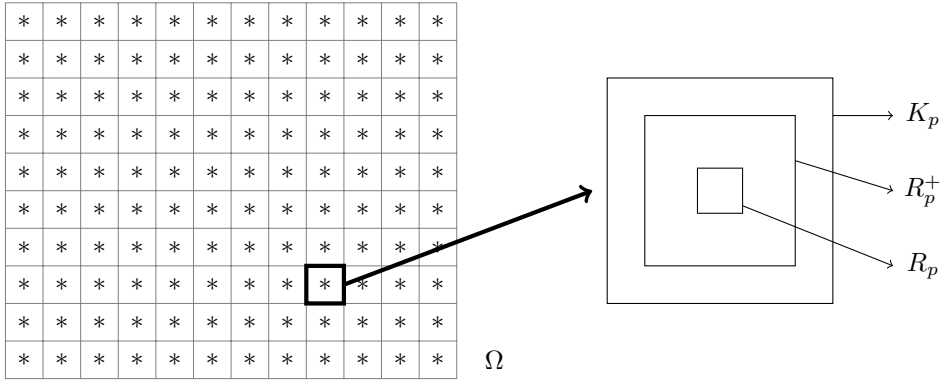

We present the computational grid used for multicontinuum homogenization. The computational domain $\Omega$ is partitioned into a series of coarse blocks, and we denote the center of a coarse block $K_p$ by $x_p^*$. The coarse partition is in principle arbitrary and may be adapted to the heterogeneities of the operator $\mathcal{A}^\epsilon$. In the examples below, we employ rectangular coarse blocks solely for ease of exposition.
In the following, the superscript $*$ is used to indicate a macropoint. To capture the local heterogeneity at each macropoint $x_p^*$, a representative volume element (RVE) $R_p$ is defined.
To mitigate the influence of artificial boundary conditions, the RVE is extended via oversampling to $R_p^+$. This procedure is particularly necessary in high-contrast cases, and the number of oversampling layers depends on the coarse mesh size and the contrast of the medium, as analyzed in \cite{maalqvist2014localization,chung2018constraint}.

The continuum is defined through different characteristic functions $\psi_i^p$ in the corresponding representative volume element $R_p$.
Specifically, $\psi_i^p = 1$ in continuum $i$ and $0$ elsewhere. 
We use two continua for illustration, but this approach can be extended to multiple continua. We then assume that the unknown variable $u$ can be expanded as a series of products of macroscopic variables $U_i$ and local microscopic functions $\phi$'s:
\[
u = \phi_i U_i + \phi_i^m \partial_m U_i 
+ \phi_i^{mn} \partial_{mn}^2 U_i + \cdots
\quad \textnormal{in}~R_p,
\]
where $\partial_m = \frac{\partial}{\partial x_m}$ 
and $\partial_{mn}^2 = \frac{\partial^2}{\partial x_m \partial x_n}$. 
We adopt the Einstein summation convention for simplicity, for example:
\[
\phi_i U_i = \sum_{i=1}^2 \phi_i U_i, ~
\phi_i^m \partial_m U_i = \sum_{i=1}^2 \sum_{m=1}^d \phi_i^m \partial_m U_i,
\]
unless otherwise stated, in which case the summation sign will be omitted.
The macroscopic function $U_i$ is smooth and defined on the computational domain $\Omega$. Its value at the center $x_p^*$ is given by
$U_i (x_p^*) = \int_{R_p} \psi_i u / \int_{R_p} \psi_i$.

The formulation of the local solutions $\phi$ is inspired by \cite{chung2018constraint,chung2021nonlocal}, where they are cast as constrained cell problems. Specifically, the local solutions are obtained by minimizing the energy associated with the PDE operator $\mathcal{A}^{\epsilon}$ in a local domain $R_p^+$, subject to constraints that preserve certain integral properties of the macroscopic functions. In this work, we restrict our attention to local problems corresponding to constant and linear macroscopic behaviors. In particular, we approximate $u$ as $u \approx \phi_i U_i + \phi_i^m \partial_m U_i$.

The first-type saddle point problem is formulated by imposing constraints designed to capture the macroscopic average behavior of multiple continua.

\begin{equation}
\begin{aligned}
\int_{R_p^+} \kappa \nabla \phi_i \cdot \nabla v - 
\sum_{j,q} \frac{\beta_{ij}^q}{\int_{R_q} \psi_j} 
\int_{R_q} \psi_j v = 0, \\
\int_{R_q} \phi_i \psi_j = \delta_{ij} \int_{R_q} \psi_j,
\quad
\forall R_q \in R_p^+.
\end{aligned}
\label{eq:phi_avar}
\end{equation}

The second-type saddle point problem imposes constraints to represent the macroscopic linear behavior of each continuum.

\begin{equation}
\begin{aligned}
\int_{R_p^+} \kappa \nabla \phi_i \cdot \nabla v - 
\sum_{j,q} \frac{\beta_{ij}^q}{\int_{R_q} \psi_j} 
\int_{R_q} \psi_j v = 0, \\
\int_{R_q} \phi_i^m \psi_j = \delta_{ij} \int_{R_q} (x_m-c_{mj})\psi_j,
\quad
\forall R_q \in R_p^+,
\end{aligned}
\label{eq:phi_grad}
\end{equation}
where $c_{mj}$ statisfy $\int_{K_p} (x_m-c_{mj}) \psi_j^{p} =0.$
In Eqs. \eqref{eq:phi_avar} and \eqref{eq:phi_grad}, 
$\beta$ denotes the Lagrange multiplier associated with the specified constraints.
We assume the macroscopic function $U_i$ are smooth in the RVE, 
the variation is small compare to the local functions $\phi$.
We have 

\begin{equation}
\begin{aligned}
\int_{R_p} \kappa \nabla (\phi_i U_i) \nabla v \approx 
U_i(x_p^*) \int_{R_p} \kappa \nabla \phi_i \nabla v, \\
\int_{R_p} \kappa \nabla (\phi_i^m \nabla_m U_i) \nabla v \approx 
\nabla_m U_i(x_p^*) \int_{R_p} \kappa \nabla \phi_i^m \nabla v.
\end{aligned}
\label{eq:phi_approx_intgral}
\end{equation}
Based on the approximation \eqref{eq:phi_approx_intgral} and setting $v \approx \phi_j V_j + \phi_j^n \partial_n V_j$ in Eq. \eqref{eq:weak_galerkin}, we have the following result:

$$
\begin{aligned}
\sum_p \frac{|K_p|}{|R_p|} \int_{R_p} f \phi_j = 
U_i V_j \sum_p \frac{|K_p|}{|R_p|} 
\int_{R_p} \kappa \nabla \phi_i \nabla \phi_j
+ V_j \nabla_m U_i \sum_p \frac{|K_p|}{|R_p|} 
\int_{R_p} \kappa \nabla \phi_i^m \nabla \phi_j \\
+ U_i \nabla_n V_j \sum_p \frac{|K_p|}{|R_p|} 
\int_{R_p} \nabla \phi_i \nabla \phi_j^n
+ \nabla_m U_i \nabla_n V_j \sum_p \frac{|K_p|}{|R_p|} 
\int_{R_p} \kappa  \nabla \phi_i^m \nabla \phi_j^n,
\end{aligned}
$$
where $U_i$, $V_j$, and their gradients are taken at specific points, such as the midpoint of the coarse block. 
Moreover, the RVE index $p$ has been omitted here for simplicity, as it is independent across different macropoints. However, it will be explicitly included in the discussion of the hierarchical method.

Then, we obtain the following upscaling multicontinuum system in strong form:

\begin{equation}
B_{ji} U_i + B_{ji}^m \partial_m U_i - \partial_n \overline{B}_{ji}^n U_i - \partial_n (B_{ji}^{mn} \partial_m U_i) = b_j,
\label{eq:macro_eq}
\end{equation}
where the coefficients are piecewise-constant vectors or matrices, 

\begin{equation}
\begin{aligned}
B_{ji} = \frac{|K_p|}{|R_p|} \int_{R_p} \kappa
\nabla \phi_i \nabla \phi_j, ~
B_{ji}^m = \frac{|K_p|}{|R_p|} \int_{R_p} \kappa
\nabla \phi_i^m\nabla \phi_j, \\
\overline{B}_{ji}^n = \frac{|K_p|}{|R_p|} \int_{R_p} \kappa
\nabla \phi_i \nabla \phi_j^n, ~
B_{ji}^{mn} = \frac{|K_p|}{|R_p|} \int_{R_p} \kappa
\nabla \phi_i^m \nabla \phi_j^n, ~
b_j = \frac{|K_p|}{|R_p|} \int_{R_p} f \phi_j,
\end{aligned}
\label{eq:Bs_defi}
\end{equation}

\section{Hierarchical multicontinuum homogenization} \label{sec:hmh}

The main goal of hierarchical multicontinuum homogenization is to use a different FE grid size to solve the cell problems, while ensuring that all of them achieve the finest grid accuracy. 
The local solutions \eqref{eq:phi_avar} and \eqref{eq:phi_grad} can be summarized as follows:
\begin{equation}
\Phi_p = \arg\min\limits_{v \in V(R_p^+)} \mathcal{J}^{\epsilon}(v),
\quad \textnormal{s.t. constraints,}
\label{eq:phi_abstract}
\end{equation}
where $\mathcal{J}^{\epsilon}(v) = \tfrac{1}{2} \int_{R_p^+} \kappa |\nabla v|^2$.
The constraints are detailed in Eqs.~\eqref{eq:phi_avar} and \eqref{eq:phi_grad}.
Hereafter, $\Phi_p$ is used to represent the local solutions at macropoint $x_p^*$, regardless of the specific type of local solutions.
Since solving the cell problems for every macroscopic point using a very fine mesh remains a huge challenge, we design a hierarchical multicontinuum homogenization method to overcome this difficulty. 

For each macropoint $x_p^*$, the local solution $\Phi_p$ is represented as a decomposition into two components:
(i) an interpolation term $\mathcal{I}_p(\Phi)$, which resolves the cell problems on a very fine mesh inherited from preceding levels, and
(ii) a correction term $\Xi_p$, computed on a coarser mesh to account for local variations.
This decomposition reads
\begin{equation}
\Phi_p = \mathcal{I}_p(\Phi) + \Xi_p. 
\label{eq:PhiXiI}
\end{equation}

The computational process of our algorithm involves three steps:
First, we build a dense hierarchical macrogrid $\{ \mathcal{S}_n \}_{n=1}^L$ for the macropoint set $\mathcal{T}$.
Second, we define a nested FE space $\{ V_n \}_{n=1}^L$ to match the hierarchical macrogrid $\{ \mathcal{S}_n \}_{n=1}^L$.
Finally, we substitute Eq.~\eqref{eq:PhiXiI} into the constraint problem \eqref{eq:phi_abstract} and solve for the correction term $\Xi_p$ within the corresponding FE space.
The framework of our method is summarized in Algorithm~\ref{alg:hmh}.

In the following, we present the details of the hierarchical computation of the local solution.

\begin{algorithm}[h!]
\caption{Hierarchical multicontinuum homogenization method.} 
\begin{algorithmic}[1]
\REQUIRE Macropoint set $\mathcal{T}$, coarsing factor $\eta$, fine grid mesh size $h$, levels $L$.
\STATE Divide $\mathcal{T}$ into a dense hierarchy of macroscopic point sets $\{\mathcal{S}_n\}_{n=1}^L$.
\STATE Define a nested FE space $\{ V_n \}_{n=1}^L$ matching the 
hierarchical macroscopic points $\{ \mathcal{S}_n \}_{n=1}^L$.
\FOR {$1 \le n \le L$}
\FOR {$x_p^* \in \mathcal{S}_n$}
\STATE Substitute the approximation \eqref{eq:PhiXiI} into \eqref{eq:phi_abstract} and solve the correction term $\Xi_p$ in $V_n$.
\STATE Compute the upscaling coefficients in $x_p^*$ by \eqref{eq:Bs_defi}.
\ENDFOR
\ENDFOR
\STATE Solving the upscaled multicontinuum system \eqref{eq:macro_eq}.
\ENSURE Upscaled solutions $U_i$.
\end{algorithmic}
\label{alg:hmh}
\end{algorithm}

\subsection{Build hierarchy of macrogrids}
A hierarchical macrogrid can help us judiciously decide the FE mesh size for cell problems at different macrogrid points. 
We first build a nested macroscopic grid (cf. \cite{brown2013efficient, park2020hierarchical, brown2017hierarchical}) for the macropoint set $\mathcal{T}$, denoted as

$$
\mathcal{T}_1 \subset 
\mathcal{T}_2 \subset
\cdots \subset 
\mathcal{T}_L = \mathcal{T}.
$$
The grid is constructed inductively, starting with an initial mesh $\mathcal{T}_1$ of size $\eta^L H$.
Here, $H$ denotes the original coarse mesh size, i.e., the coarse mesh size at the last level $\mathcal{T}_L$, and $\eta$ is the grid coarsening factor.
At each subsequent step, the refinement $\mathcal{T}_n$ is obtained from the preceding grid $\mathcal{T}_{n-1}$, with the mesh size reduced to $\eta^{L-n}H$.
For the given $n$-th level, we ensure that there is a macropoint in the coarse block with mesh size $\eta^{L-n} H$. To avoid repeated computations, we remove the overlapping macropoints to build the hierarchical macrogrids $\{ \mathcal{S}_n \}_{n=1}^L$, that is,

\begin{equation}
\mathcal{S}_n = \mathcal{T}_n \big\backslash
\left(
\bigcup_{k< n} \mathcal{S}_k
\right).
\label{eq:defi_Sn}
\end{equation}
Specifically, for the first level, we define $\mathcal{S}_1 = \mathcal{T}_1$.
From \eqref{eq:defi_Sn}, we can also obtain the following property: $\mathcal{S}_{n+1} = \mathcal{T}_{n+1} \backslash \mathcal{T}_n$.
Constructed in this way, the macrogrids have the dense property. That is, for any macroscopic point $x_p^* \in \mathcal{S}_n$, there exists at least one point $x_q^* \in \bigcup_{k<n}~ \mathcal{S}_k$ such that $\mathrm{dist}(x_p^*, x_q^*) < O(\eta^{L-n} H)$. This paves the way for our next nested FE space definition and also the linear interpolation operator.

\begin{figure}[h!]
\centering
\begin{tikzpicture}
\foreach \i in {0,1} {
\foreach \j in {0,1} {
\node at (\i*2+1.25,\j*2+1.2) {*};
}}
\foreach \i in {0,...,3} {
\foreach \j in {0,...,3} {
\node at (\i+4.75,\j+0.2) {*};
}}
\foreach \i in {0,...,7} {
\foreach \j in {0,...,7} {
\node at (\i*0.5+9.25,\j*0.5+0.2) {*};
}}
\draw[help lines,step=0.5cm] (0,0) grid (4,4);
\draw[help lines,step=0.5cm,shift={(4.5,0)}] (0,0) grid (4,4);
\draw[thick,step=0.5cm,shift={(9,0)}] (0,0) grid (4,4);
\draw[thick,step=2cm] (0,0) grid (4,4);
\draw[thick,shift={(4.5,0)}] (0,0) grid (4,4);
\node at (2,-1) {(a) $\mathcal{T}_1$};
\node at (6.5,-1) {(b) $\mathcal{T}_2$};
\node at (11,-1) {(c) $\mathcal{T}_3$};
\end{tikzpicture}
\caption{3-level nested macrogrids.}
\label{fig:grid_3level_T}
\end{figure}

In Figures \ref{fig:grid_3level_T} and \ref{fig:grid_3level_S},
we illustrate an example of three-level nested macrogrids $\{ \mathcal{T}_n \}_{n=1}^3$ and the corresponding hierarchical structures
$\{ \mathcal{S}_n \}_{n=1}^3$ in a unit square.
Here, the final coarse mesh size is $H=1/8$, and the coarsening factor is $\eta=2$.
As the coarse mesh decreases from $4H$ to $2H$, and then to $H$, 
a macroscopic point remains within the corresponding coarse block in $\mathcal{T}_n$.

\begin{figure}[h!]
\centering
\begin{tikzpicture}
\foreach \i in {0,1} {
\foreach \j in {0,1} {
\node at (\i*2+1.25,\j*2+1.2) {*};
}}
\foreach \i in {0,1} {
\foreach \j in {0,...,3} {
\node at (\i*2+4.75,\j+0.2) {*};
}}
\foreach \i in {0,1} {
\foreach \j in {0,1} {
\node at (\i*2+5.75,\j*2+0.2) {*};
}}
\foreach \i in {0,...,3} {
\foreach \j in {0,...,7} {
\node at (\i+9.75,\j*0.5+0.2) {*};
}}
\foreach \i in {0,...,3} {
\foreach \j in {0,...,3} {
\node at (\i+9.25,\j+0.7) {*};
}}
\draw[thick,step=0.5cm] (0,0) grid (4,4);
\draw[thick,step=0.5cm,shift={(4.5,0)}] (0,0) grid (4,4);
\draw[thick,step=0.5cm,shift={(9,0)}] (0,0) grid (4,4);
\node at (2,-1) {(a) $\mathcal{S}_1$};
\node at (6.5,-1) {(b) $\mathcal{S}_2$};
\node at (11,-1) {(c) $\mathcal{S}_3$};
\end{tikzpicture}
\caption{3-level hierarchical macrogrids.}
\label{fig:grid_3level_S}
\end{figure}

\subsection{Build nested FE spaces}
For a given macrogrid point $x_p^* \in \mathcal{S}_n$,
we denote the corresponding FE space as $V_n$.
More specifically, the FE space $V_n$ is defined on the FE grid with a mesh size of $h \eta^{n-1}$.
The abstract problem \eqref{eq:phi_abstract} for $x_p^* \in \mathcal{S}_n$ will be solved using the mesh size $h \eta^{n-1}$.
It is clear that $V_L \subset V_{L-1} \subset \cdots \subset V_1 \subset V$.
This implies that as the hierarchical level $n$ increases from $1$ to $L$, the FE grid becomes coarser.
By combining the definition of $\{ \mathcal{S}_n \}_{n=1}^L$, it can be observed that only a few macropoints are equipped with the finest FE space, while the majority are equipped with moderately fine FE spaces.
It should be noted that $\eta^L h$ remains sufficient to resolve the heterogeneities at all scales.

\subsection{Calculating the correction term}

We now relate the hierarchical macrogrids $\{ \mathcal{S}_n \}_{n=1}^L$, the nested FE spaces $\{ V_n \}_{n=1}^L$, and our assumption in \eqref{eq:PhiXiI}. 
For a given macropoint $x_p^* \in \mathcal{S}_n$, we define a collection $\{ x_1^*, x_2^*, \cdots, x_r^* \} \subset \bigcup_{k<n} \mathcal{S}_k$. In the following, we refer to this collection as the neighboring macropoints of $x_p^*$.
The macropoints in this collection must be sufficiently close to $x_p^*$. More specifically, for any macropoint $x_q^*$ in this collection, we have $\mathrm{dist}(x_p^*, x_q^*) < O(\eta^{L-n} H)$. In the following, we provide a detailed explanation of the assumption in \eqref{eq:PhiXiI},

\begin{equation}
\phi_{i,p} = \xi_{i,p} + \sum_{t=1}^r c_t \phi_{i,t}, \quad
\phi_{i,p}^m = \xi_{i,p}^m + \sum_{t=1}^r c_t \phi_{i,t}^m,
\label{eq:phiXiI_all}
\end{equation}
where $\phi_{i,t}$ represents the average behavior of the $i$-th continuum at the neighboring macropoint $x_t^*$, and $\phi_{i,t}^m$ represents the linear behavior of the $i$-th continuum along the $m$-th axis at the neighboring macropoint $x_t^*$.
For $x_p^* \in \mathcal{S}_1$, as no pre-constructed level exists, we define $\mathcal{I}_p = 0$.
Note that the linear interpolation coefficients satisfy $\sum_{t=1}^r c_t = 1$.
For the macropoint $x_p^* \in \mathcal{S}_n$, we substitute the local solutions $\phi_{i,p}$ from \eqref{eq:phiXiI_all} into \eqref{eq:phi_avar}, and we obtain:

\begin{equation}
\begin{aligned}
\int_{R_p^+} \kappa \nabla \xi_{i,p} \cdot \nabla v - 
\sum_{j,q} \frac{\beta_{ij}^q}{\int_{R_q} \psi_j} 
\int_{R_q} \psi_j v = - \sum_{t} c_t
\int_{R_p^+} \kappa \nabla \phi_{i,t} \cdot \nabla v , \\
\int_{R_q} \xi_{i,p} \psi_j = 
\int_{R_q} (\delta_{ij} - \sum_t c_t \phi_{i,t}) \psi_j,
\quad
\forall R_q \in R_p^+.
\end{aligned}
\label{eq:phi_avar_corrector}
\end{equation}
After substituting the local solutions $\phi_{i,p}^m$ from \eqref{eq:phiXiI_all} into \eqref{eq:phi_grad}, we obtain another correction term:

\begin{equation}
\begin{aligned}
\int_{R_p^+} \kappa \nabla \xi_{i,p}^m \cdot \nabla v - 
\sum_{j,q} \frac{\beta_{ij}^{mq}}{\int_{R_q} \psi_j} 
\int_{R_q} \psi_j v = - \sum_{t} c_t
\int_{R_p^+} \kappa \nabla \phi_{i,t} \cdot \nabla v , \\
\int_{R_q} \xi_{i,p}^m \psi_j = 
\int_{R_q} (\delta_{ij}(x_m-c_{mj}) - \sum_t c_t \phi_{i,t}^m) \psi_j,
\quad
\forall R_q \in R_p^+.
\end{aligned}
\label{eq:phi_grad_corrector}
\end{equation}

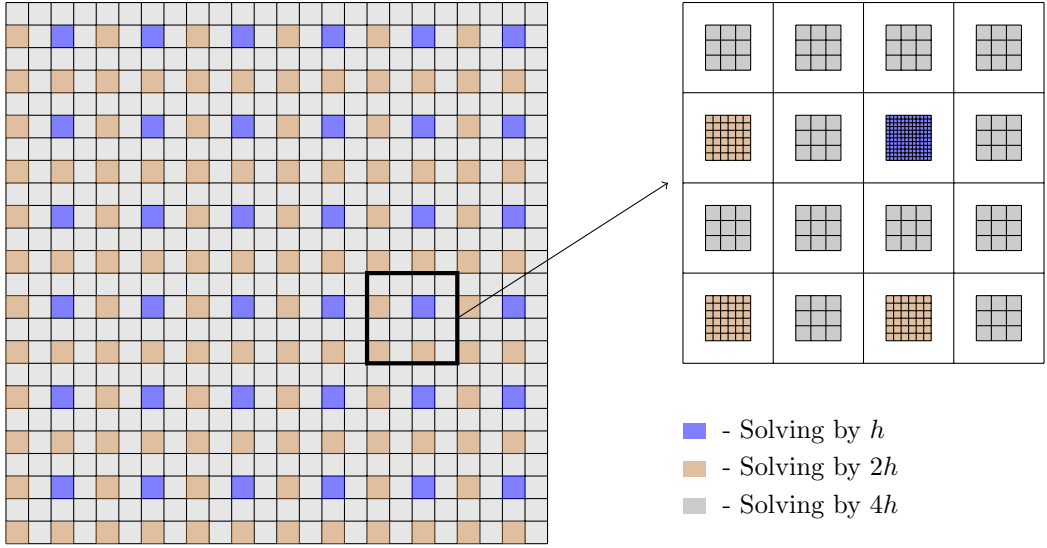
\begin{figure}[h!]
\centering
\begin{tikzpicture}
\draw[step=0.3] (0,0) grid (7.2,7.2);
\foreach \i in {0,...,5} {
\foreach \j in {0,...,5} {
\fill[blue,opacity=0.5] (\i*1.2+0.6,\j*1.2+0.6) rectangle (\i*1.2+0.9,\j*1.2+0.9);
}}
\foreach \i in {0,...,5} {
\foreach \j in {0,...,11} {
\fill[brown,opacity=0.5] (\i*1.2,\j*0.6) rectangle (\i*1.2+0.3,\j*0.6+0.3);
}}
\foreach \i in {0,...,5} {
\foreach \j in {0,...,5} {
\fill[brown,opacity=0.5] (\i*1.2+0.6,\j*1.2) rectangle (\i*1.2+0.9,\j*1.2+0.3);
}}
\foreach \i in {0,...,11} {
\foreach \j in {0,...,23} {
\fill[opacity=0.1] (\i*0.6+0.3,\j*0.3) rectangle (\i*0.6+0.6,\j*0.3+0.3);
}}
\foreach \i in {0,...,11} {
\foreach \j in {0,...,11} {
\fill[opacity=0.1] (\i*0.6,\j*0.6+0.3) rectangle (\i*0.6+0.3,\j*0.6+0.6);
}}
\draw[ultra thick,shift ={(4.8,2.4)}] (0,0) rectangle (1.2,1.2);
\draw[step=1.2,shift ={(8,2.4)}] (0,0) grid (4.8,4.8);
\draw[->] (6,3) -- (7.8,4.8);
\fill[blue,opacity=0.5,shift={(10.7,5.1)}] (0,0) rectangle (0.6,0.6);
\draw[step=0.05,shift={(10.7,5.1)}] (0,0) grid (0.6,0.6);
\fill[brown,opacity=0.5,shift={(8.3,2.7)}] (0,0) rectangle (0.6,0.6);
\draw[step=0.1,shift={(8.3,2.7)}] (0,0) grid (0.6,0.6);
\fill[brown,opacity=0.5,shift={(8.3,5.1)}] (0,0) rectangle (0.6,0.6);
\draw[step=0.1,shift={(8.3,5.1)}] (0,0) grid (0.6,0.6);
\fill[brown,opacity=0.5,shift={(10.7,2.7)}] (0,0) rectangle (0.6,0.6);
\draw[step=0.1,shift={(10.7,2.7)}] (0,0) grid (0.6,0.6);
\foreach \i in {0,1} {
\foreach \j in {0,1} {
\fill[opacity=0.2,shift={(\i*2.4+8.3,\j*2.4+3.9)}] (0,0) rectangle (0.6,0.6);
\draw[step=0.2,shift={(\i*2.4+8.3,\j*2.4+3.9)}] (0,0) grid (0.6,0.6);
}}
\foreach \i in {0,1} {
\foreach \j in {0,...,3} {
\fill[opacity=0.2,shift={(\i*2.4+9.5,\j*1.2+2.7)}] (0,0) rectangle (0.6,0.6);
\draw[step=0.2,shift={(\i*2.4+9.5,\j*1.2+2.7)}] (0,0) grid (0.6,0.6);
}}
\node at (8.4,1.5) [anchor=west] {- Solving by $h$};
\node at (8.4,1) [anchor=west] {- Solving by $2h$};
\node at (8.4,0.5) [anchor=west] {- Solving by $4h$};
\fill[opacity=0.2] (8,0.4) rectangle (8.3,0.6);
\fill[brown,opacity=0.5] (8,0.9) rectangle (8.3,1.1);
\fill[blue,opacity=0.5] (8,1.4) rectangle (8.3,1.6);
\end{tikzpicture}
\caption{Illustration of 3-level hierarchical macropoint and fine-grid structure for cell problems.}
\label{fig:grid_3level_h2h4h}
\end{figure}

In Figure \ref{fig:grid_3level_h2h4h}, the coarse blocks are color-coded to illustrate the three-level hierarchical macrogrids. The color scheme is as follows: blue for the first level $\mathcal{S}_1$, brown for the second level $\mathcal{S}_2$, and light gray for the third level $\mathcal{S}_3$. The figure also shows the FE grids used to resolve the cell problems for constructing the correction terms \eqref{eq:phi_avar_corrector} and \eqref{eq:phi_grad_corrector}.  
As observed in Figure \ref{fig:grid_3level_h2h4h}, the number of macropoints is smallest at the first level and largest at the third, reflecting the strategy of performing the most expensive fine-scale computations on a limited number of macropoints while relying on coarser approximations for the majority.

In the numerical experiments, only macropoints at the first level are employed in the interpolation operator, i.e., $\mathcal{I}_p = \Phi_t$ with $x_t^* \in \mathcal{S}_1$, although more sophisticated combinations can be considered. Figure \ref{fig:grid_3level_h2h4h} highlights a first-level macropoint along with neighboring points that inherit its highly resolved local solutions. We observe that only one macropoint requires solving the cell problem on a grid of size $\mathcal{O}(h^{-2})$, three macropoints on $\mathcal{O}(4h^{-2})$, and twelve on $\mathcal{O}(16h^{-2})$. Compared with fully resolving all sixteen macropoints at the finest level $\mathcal{O}(h^{-2})$, this hierarchical strategy achieves significant computational savings. More general computational costs are discussed in the next subsection, and the performance and cost reductions associated with different macrogrid levels are reported in the numerical examples section.

\subsection{Computational savings}

Now, we compare the computational complexity of the original multicontinuum homogenization method and our hierarchical algorithm.
Again, let $h$ denote the finest FE mesh size for the constraint problems \eqref{eq:phi_abstract}, $H$ denote the initial coarse mesh size, and $\eta$ be the coarsening factor.
In the hierarchical approach, for the $n$-th level, the degrees of freedom of the local problems at one macropoint are $\mathcal{O}\left((\eta^{n-1}h)^{-d}\right)$, while the number of $n$-th macropoints is  $\mathcal{O}\left( (\eta^{L-n}H)^{-d} \right)$.
The computational complexity for one level is $\mathcal{O}\left( (hH)^{-d} \eta^{(1-L)d} \right)$.
Since there are $L$ levels, the total computational cost is $\mathcal{O}\left( L (hH)^{-d} \eta^{(1-L)d}\right)$.
In contrast, for the original multicontinuum homogenization method,
the degrees of freedom for local problems at one macropoint are $\mathcal{O}(h^{-d})$, while the number of macropoints is $\mathcal{O}(H^{-d})$.
Thus, the total computational cost is $\mathcal{O}\left( (hH)^{-d} \right)$.
Note that $\mathcal{O}(L \eta^{(1-L)d} ) < \mathcal{O}(1)$, which highlights the significant reduction in computational cost achieved by our hierarchical approach.

\section{Numerical examples} \label{sec:numerres}

In this section, we present five numerical examples to demonstrate the efficiency of our method, including both 2D and 3D cases.  
The permeability field is defined as $\kappa(x) = k(x) g(x)$, where $g(x)$ represents a slowly varying function over the global domain: $g(x) = 2+\sin(\pi x_1)\sin(\pi x_2)$ for 2D cases, and $g(x) = e^{(x_1+x_2+x_3)/2}$ for the 3D case. The function $k(x)$ describes a high-contrast distribution across two continua, with $k(x) = \eta$ in the low-conductivity region $\Omega_1$ and $k(x) = 1$ in the high-conductivity region $\Omega_2$. We set $\eta = 10^{-4}$ for 2D examples and $\eta = 10^{-2}$ for the 3D case. The specific configurations of $k(x)$ are provided in the corresponding examples. 

The source term is defined as
\[
f(x) = 
\begin{cases}
\eta \, e^{-40 \left((x_1-0.5)^2 + (x_2-0.5)^2\right)}, & x \in \Omega_1, \\
e^{-40 \left((x_1-0.5)^2 + (x_2-0.5)^2\right)}, & x \in \Omega_2.
\end{cases}
\]

The computational domain is $\Omega = [0,1]^d$. For the 2D examples, the fine-grid solution is obtained using a rectangular mesh of size $h=1/960$, while for the 3D example, a hexahedral mesh with $h=1/120$ is employed. Each coarse element is treated as a representative volume element (RVE), i.e., $R_p = K_p$, and the number of oversampling layers is set to $l = \lceil -2 \log(H) \rceil$ unless otherwise specified.
All numerical experiments were carried out on an HPC cluster node equipped with dual AMD EPYC 9654 CPUs (192 physical cores, 384 threads) and 512 GB of memory.

To quantify the performance of the upscaling techniques, we define the average fine-grid solution over each coarse block as
\begin{equation}
U_h^{(i)} = \frac{1}{|K_p \cap \Omega_i|} \int_{K_p \cap \Omega_i} u, 
\label{eq:U_h_defi}
\end{equation}
and the corresponding multiscale solutions as
\begin{equation}
U_{\textnormal{MH}}^{(i)} =
\frac{1}{|K_p|} \int_{K_p} U_{i,\textnormal{MH}},
\quad
U_{\textnormal{HMH}}^{(i)} =
\frac{1}{|K_p|} \int_{K_p} U_{i,\textnormal{HMH}},
\label{eq:U_MH_defi}
\end{equation}
where ``MH'' and ``HMH'' denote the multicontinuum homogenization and hierarchical multicontinuum homogenization methods, respectively.  

To compare two different average solutions for each continuum, we define
\[
e_2^{(i)}(U, V) = \sqrt{\frac{|U-V|^2}{|U|^2}}.
\]
In this work, we consider three types of errors:
\begin{itemize}
\item \textbf{Type 1:} $e_2^{(i)}(U_h, U_{\textnormal{MH}})$, which measures the error of the multicontinuum homogenization solution using the fine-grid solution as the reference. This quantifies the performance of the original multicontinuum homogenization method.
\item \textbf{Type 2:} $e_2^{(i)}(U_h, U_{\textnormal{HMH}})$, which measures the error of the hierarchical multicontinuum homogenization solution relative to the fine-grid reference. This assesses the accuracy of the hierarchical approach.
\item \textbf{Type 3:} $e_2^{(i)}(U_{\textnormal{MH}}, U_{\textnormal{HMH}})$, which quantifies the difference between the hierarchical solution and the standard multicontinuum homogenization solution. This isolates the error introduced by the hierarchical algorithm itself.
\end{itemize}

\subsection{Example 1}
\begin{figure}[h!]
\centering
\includegraphics[width=0.46\linewidth]{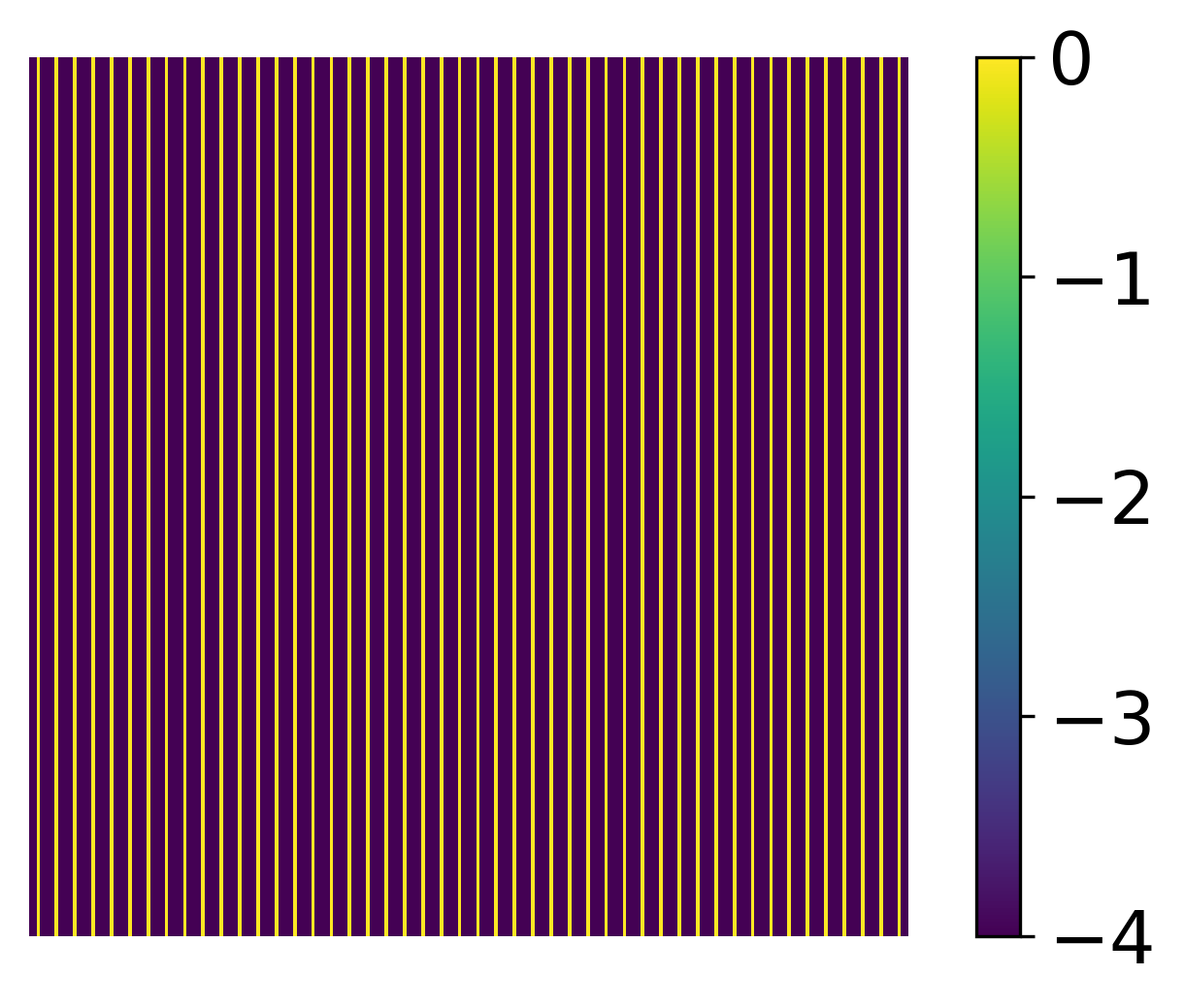}
\includegraphics[width=0.5\linewidth]{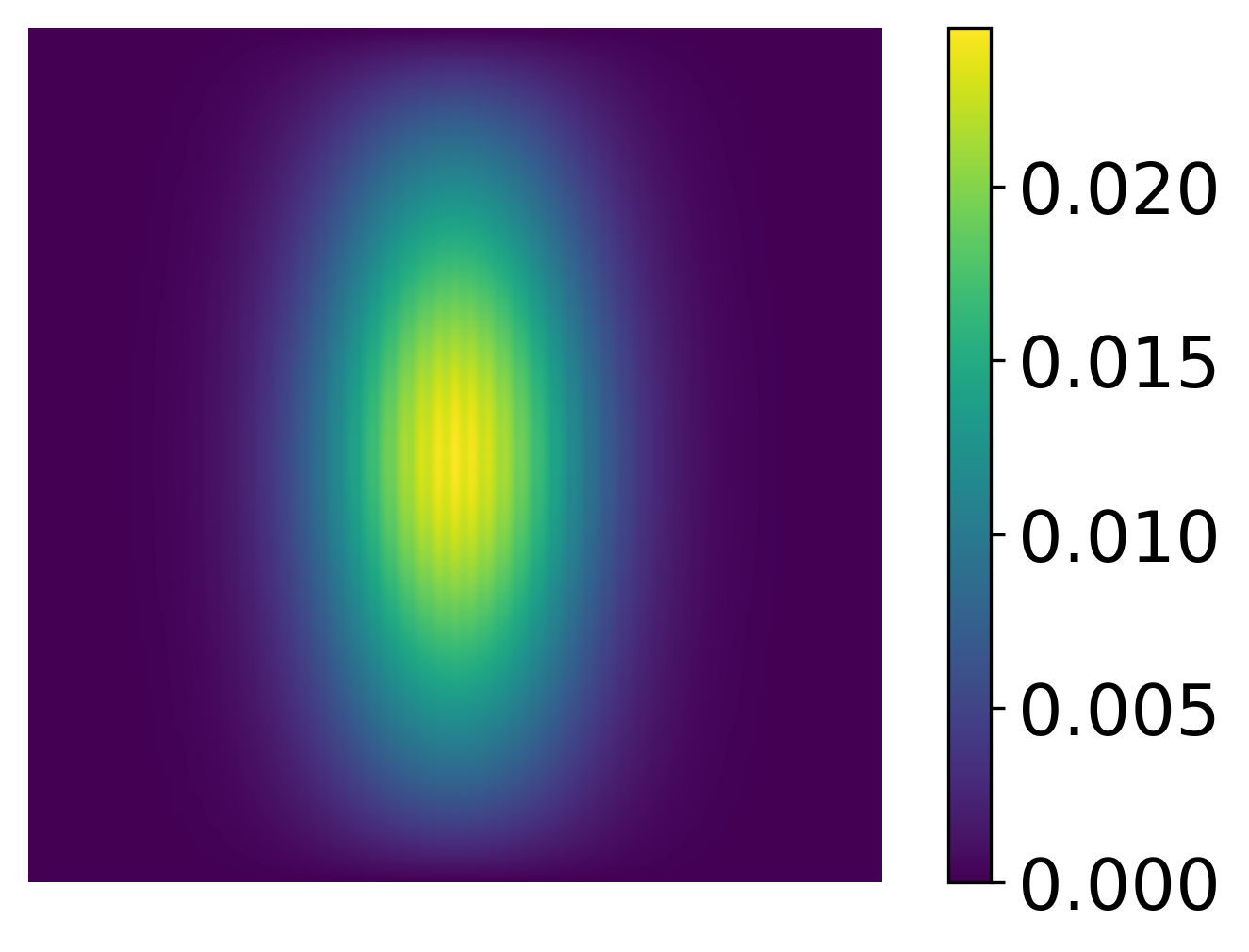}
\caption{High-contrast part $\log(k)$ (left) and reference solution (right) in Example 1.}
\label{fig:k1u}
\end{figure}

\begin{figure}[h!]
\centering
\begin{subfigure}[b]{0.48\textwidth}
\includegraphics[width=\linewidth]{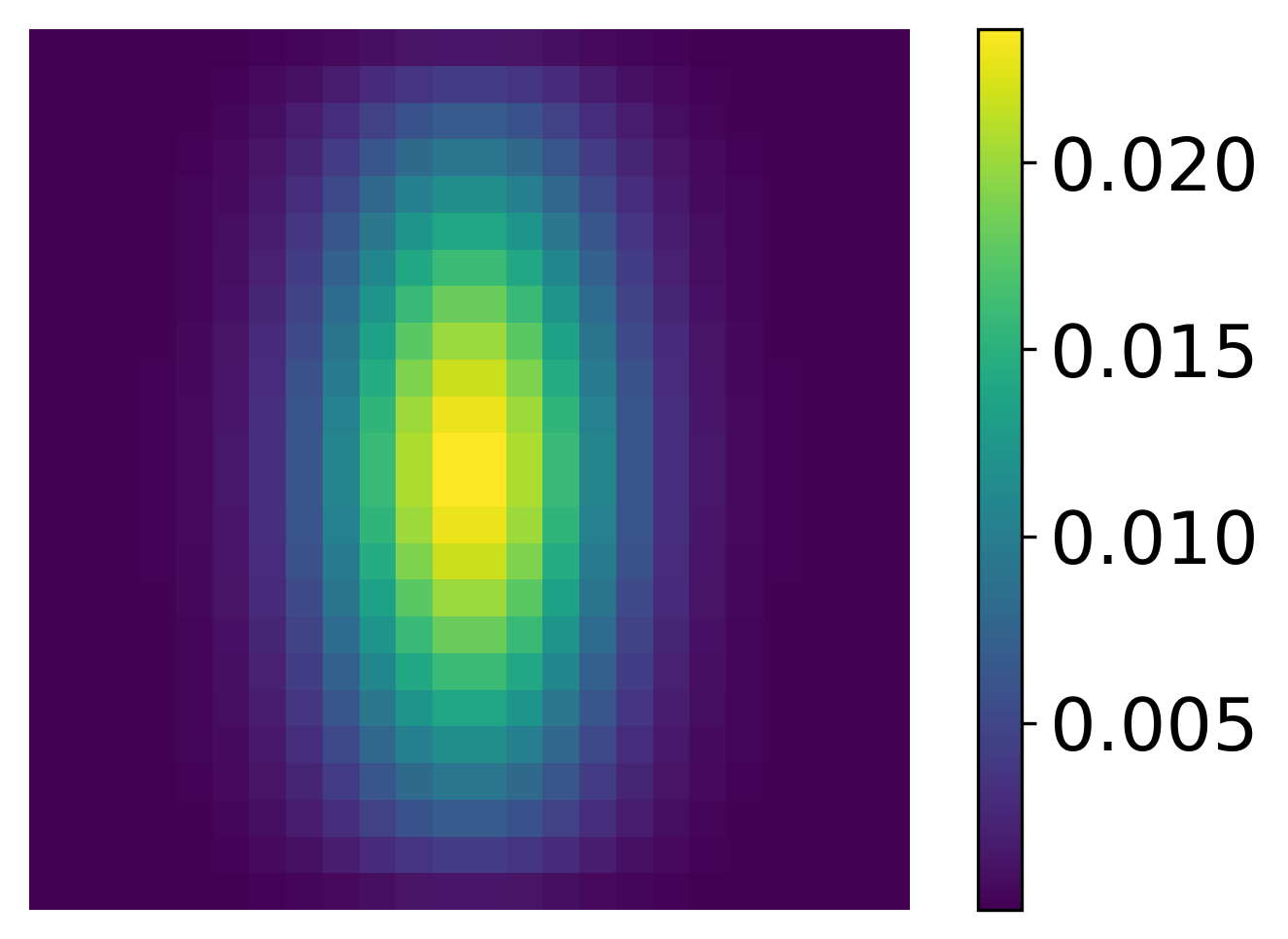}
\caption{$U_h^{(1)}$}
\end{subfigure}
\begin{subfigure}[b]{0.48\textwidth}
\includegraphics[width=\linewidth]{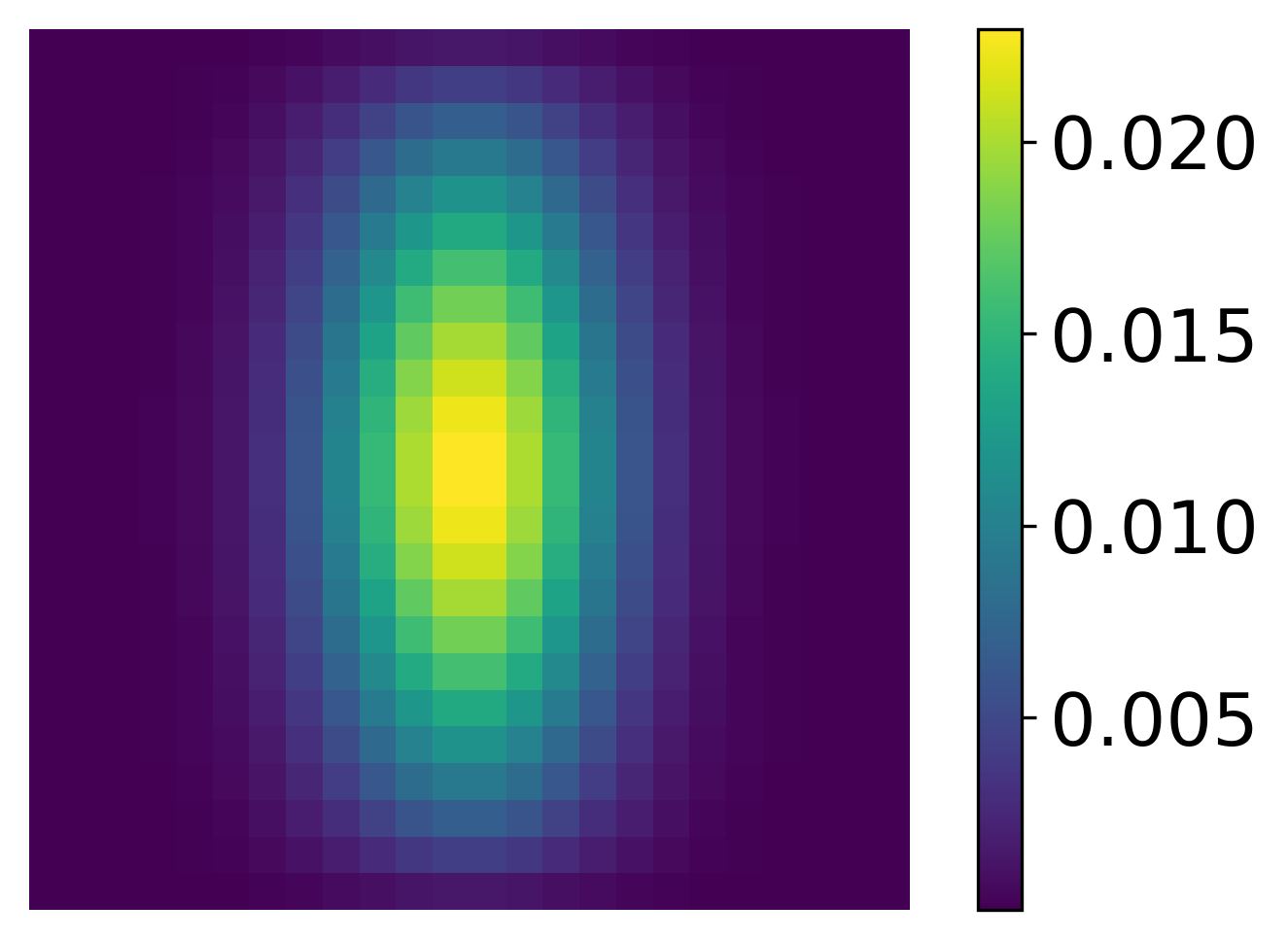}
\caption{$U_h^{(2)}$}
\end{subfigure}
\begin{subfigure}[b]{0.48\textwidth}
\includegraphics[width=\linewidth]{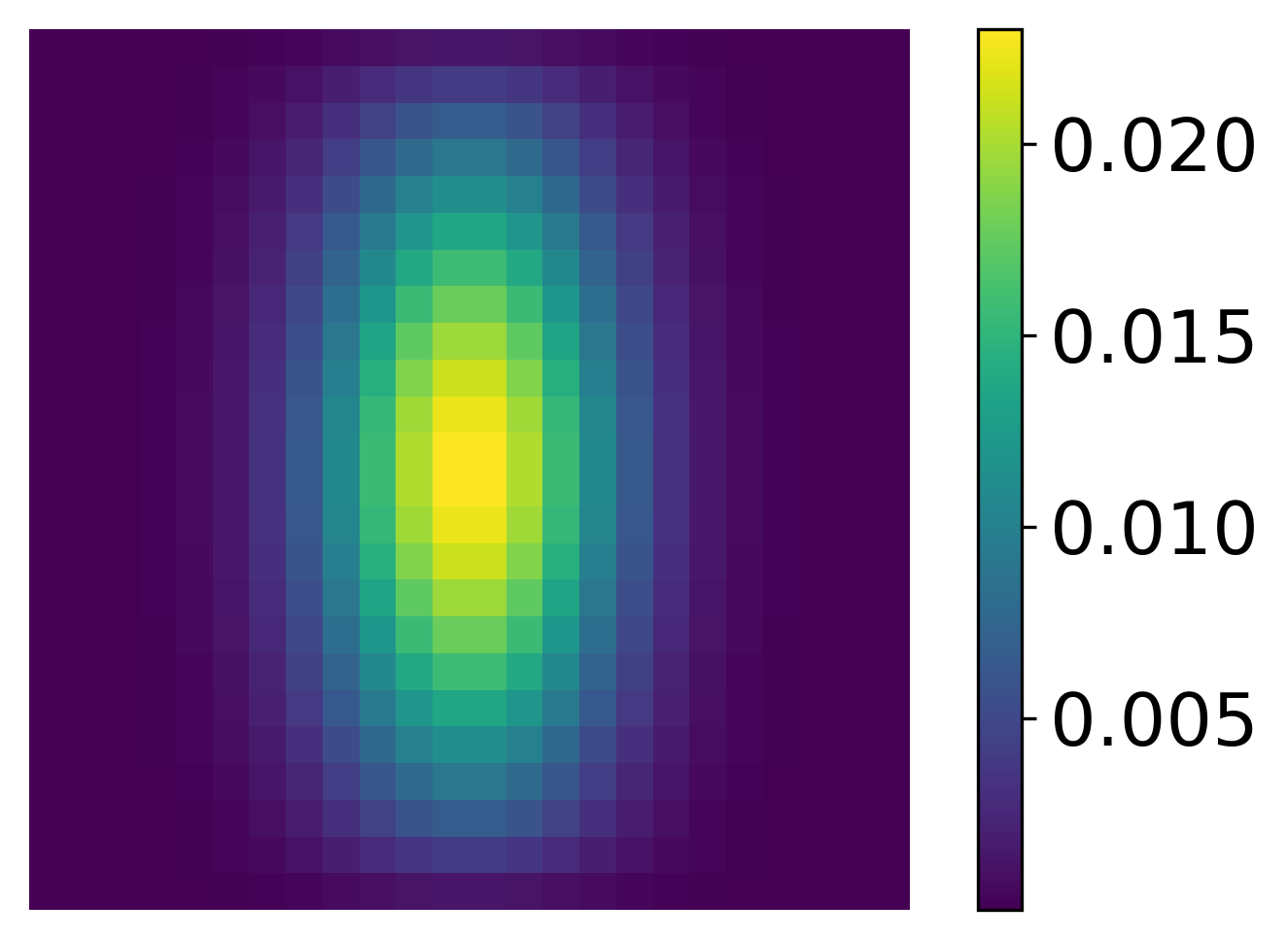}
\caption{$U_{\textnormal{MH}}^{(1)}$}
\end{subfigure}
\begin{subfigure}[b]{0.48\textwidth}
\includegraphics[width=\linewidth]{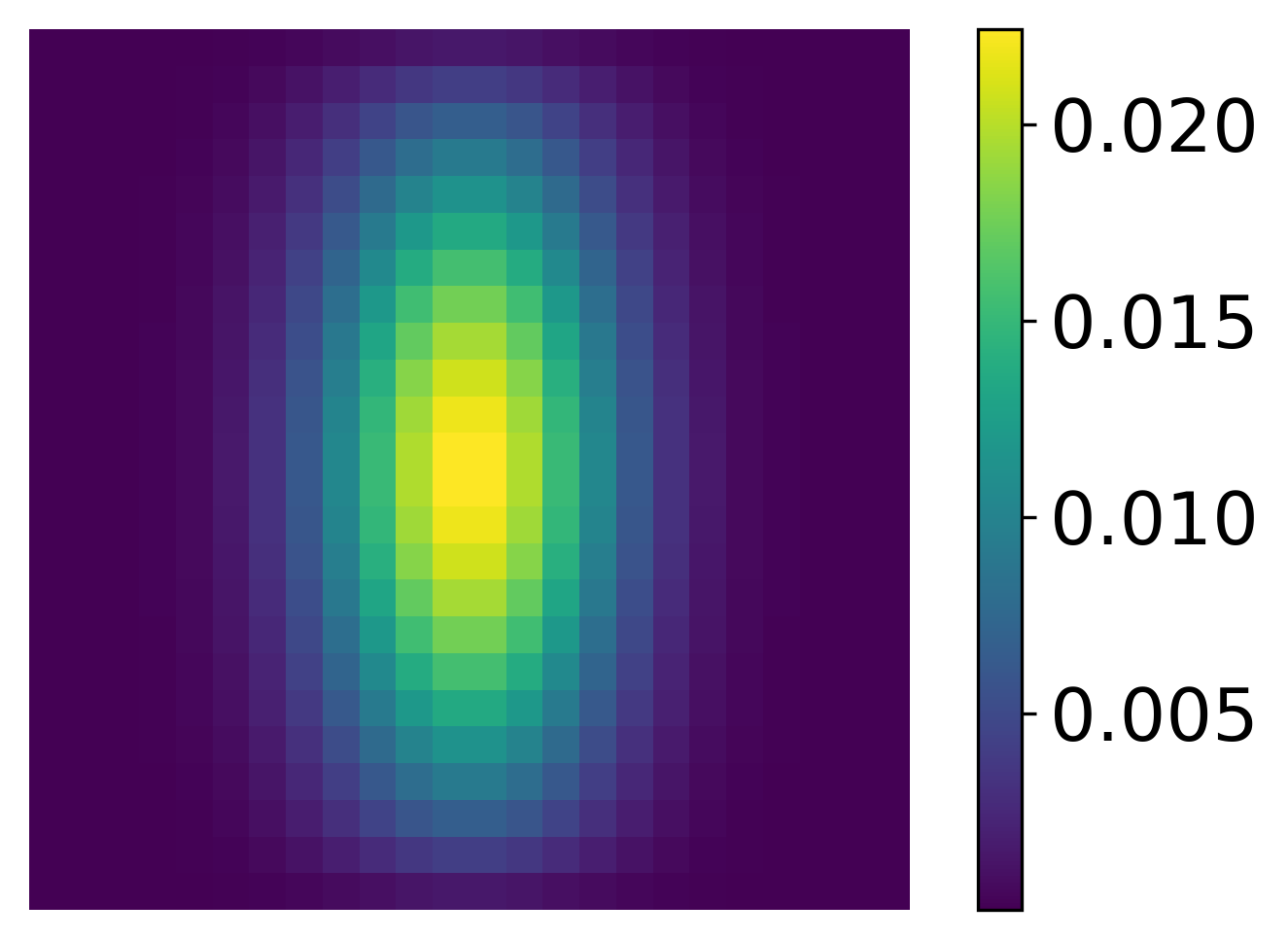}
\caption{$U_{\textnormal{MH}}^{(2)}$}
\end{subfigure}
\begin{subfigure}[b]{0.48\textwidth}
\includegraphics[width=\linewidth]{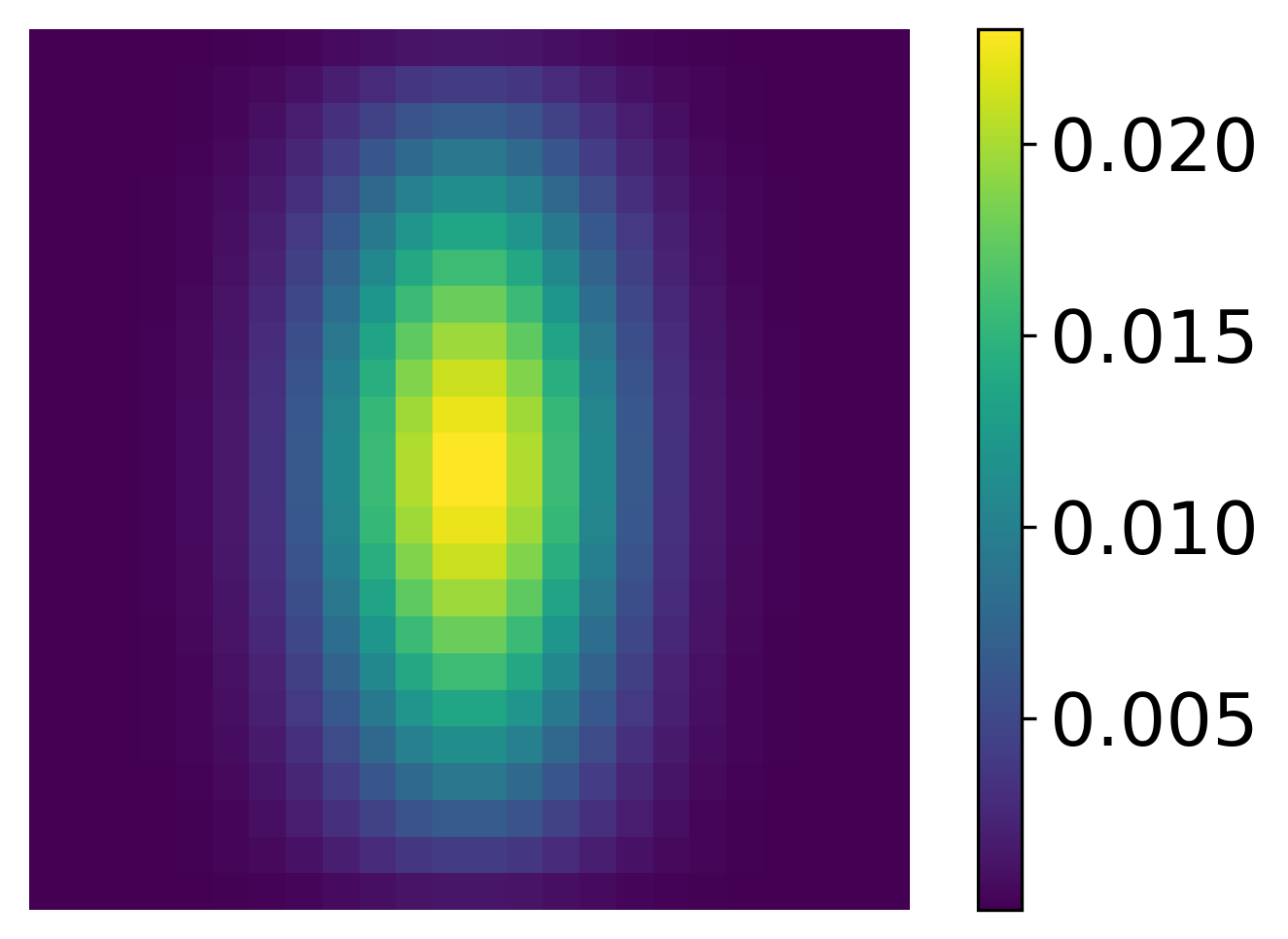}
\caption{$U_{\textnormal{HMH}}^{(1)}$}
\end{subfigure}
\begin{subfigure}[b]{0.48\textwidth}
\includegraphics[width=\linewidth]{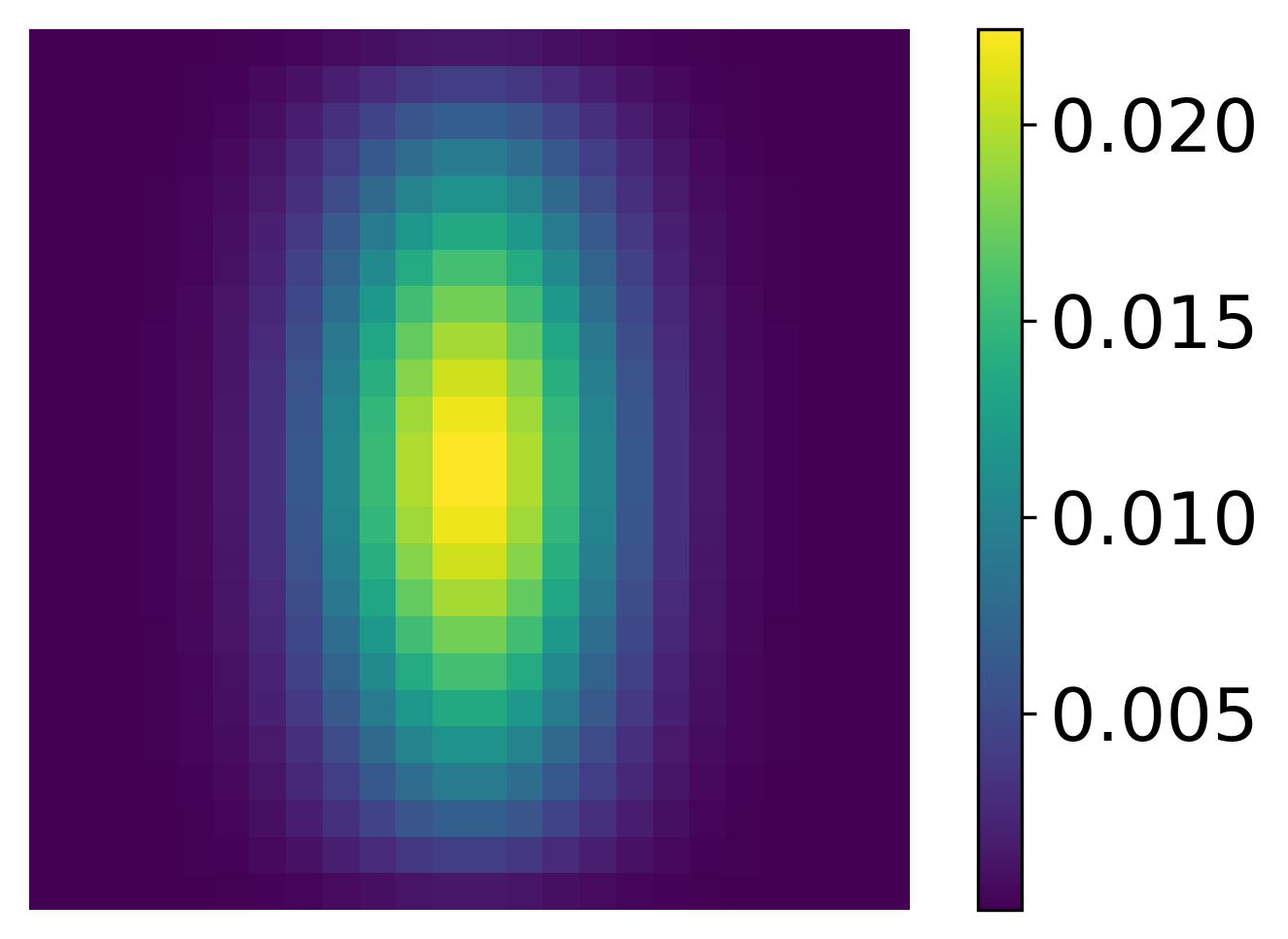}
\caption{$U_{\textnormal{HMH}}^{(2)}$}
\end{subfigure}
\caption{Average solution in Example 1.}
\label{fig:ex1_auhs}
\end{figure}

We begin with a layered medium characterized by a high-contrast function $k$. The logarithm of $k$ together with the reference solution is displayed in Figure \ref{fig:k1u}, while the corresponding three-level hierarchical structure ($L=3$) is illustrated in Figure \ref{fig:grid_3level_h2h4h}. 
CPU timings for MH and HMH at one macropoint and its neighbors are reported in Table \ref{tab:cpu2dL2} for $H=1/12$ with oversampling of $l=5$ layers. 
The hierarchical approach achieves a reduction of nearly half in cost (101.21 seconds for MH vs.~46.9 seconds for HMH). Table \ref{tab:Herror1} reports three error measures. The first two, corresponding to the two upscaling methods, decrease consistently as the coarse mesh is refined. The third, which quantifies the difference between MH and HMH, does not exhibit convergence but stays at a negligible level for all tested mesh sizes.
 The hierarchical method delivers accuracy comparable to MH, with differences remaining negligible; for instance, at $H=1/24$, both methods produce errors on the order of $10^{-2}$, differing by only $10^{-7}$. The average solutions of the reference, MH, and hierarchical approaches are plotted in Figure \ref{fig:ex1_auhs}, confirming both accuracy and robustness.

\begin{table}[h!]
\renewcommand\arraystretch{1.5}
\centering
\begin{tabular}{ccccc}
\hline
& Macropoints & DOF & CPU (s) \\ \hline
MH & 16 & $\mathcal{O}\left( (\frac{1}{h})^2 \right)$ & 101.21  \\ \hline
\multirow{3}*{HMH} & 1 & $\mathcal{O}\left( (\frac{1}{h})^2 \right)$ & \multirow{3}*{46.90}  \\ \cline{2-3}
 & 3 & $\mathcal{O}\left( (\frac{2}{h})^2 \right)$ &  \\ \cline{2-3}
 & 12 & $\mathcal{O}\left( (\frac{4}{h})^2 \right)$ &  \\ \hline
\end{tabular}
\caption{Comparison of DOF and CPU time of cell problems for MH and HMH methods in Example 1.}
\label{tab:cpu2dL3}
\end{table}

\begin{table}[h!]
\renewcommand\arraystretch{1.5}
\centering
\begin{tabular}{cccccccc}
\hline
\multirow{2}*{$H$}
& \multirow{2}*{$l$}
& \multicolumn{2}{c}{Type 1} 
& \multicolumn{2}{c}{Type 2} 
& \multicolumn{2}{c}{Type 3}
\\ \cline{3-8}
& & $e_2^{(1)}$ & $e_2^{(2)}$ & $e_2^{(1)}$ & $e_2^{(2)}$ & $e_2^{(1)}$ & $e_2^{(2)}$ \\ \hline
1/12 & 5
& 8.34e-02 & 8.04e-02 & 8.34e-02 & 8.04e-02 & 1.90e-06 & 2.12e-06
\\ \hline
1/24 & 7
& 1.92e-02 & 1.74e-02 & 1.92e-02 & 1.74e-02 & 1.02e-07 & 4.22e-08
\\ \hline
1/48 & 8
& 2.39e-03 & 2.32e-03 & 2.39e-03 & 2.32e-03 & 1.16e-07 & 1.16e-07
\\ \hline
\end{tabular}
\caption{Error comparison by different coarse mesh in Example 1.}
\label{tab:Herror1}
\end{table}

\subsection{Example 2}
\begin{figure}[h!]
\centering
\includegraphics[width=0.45\linewidth]{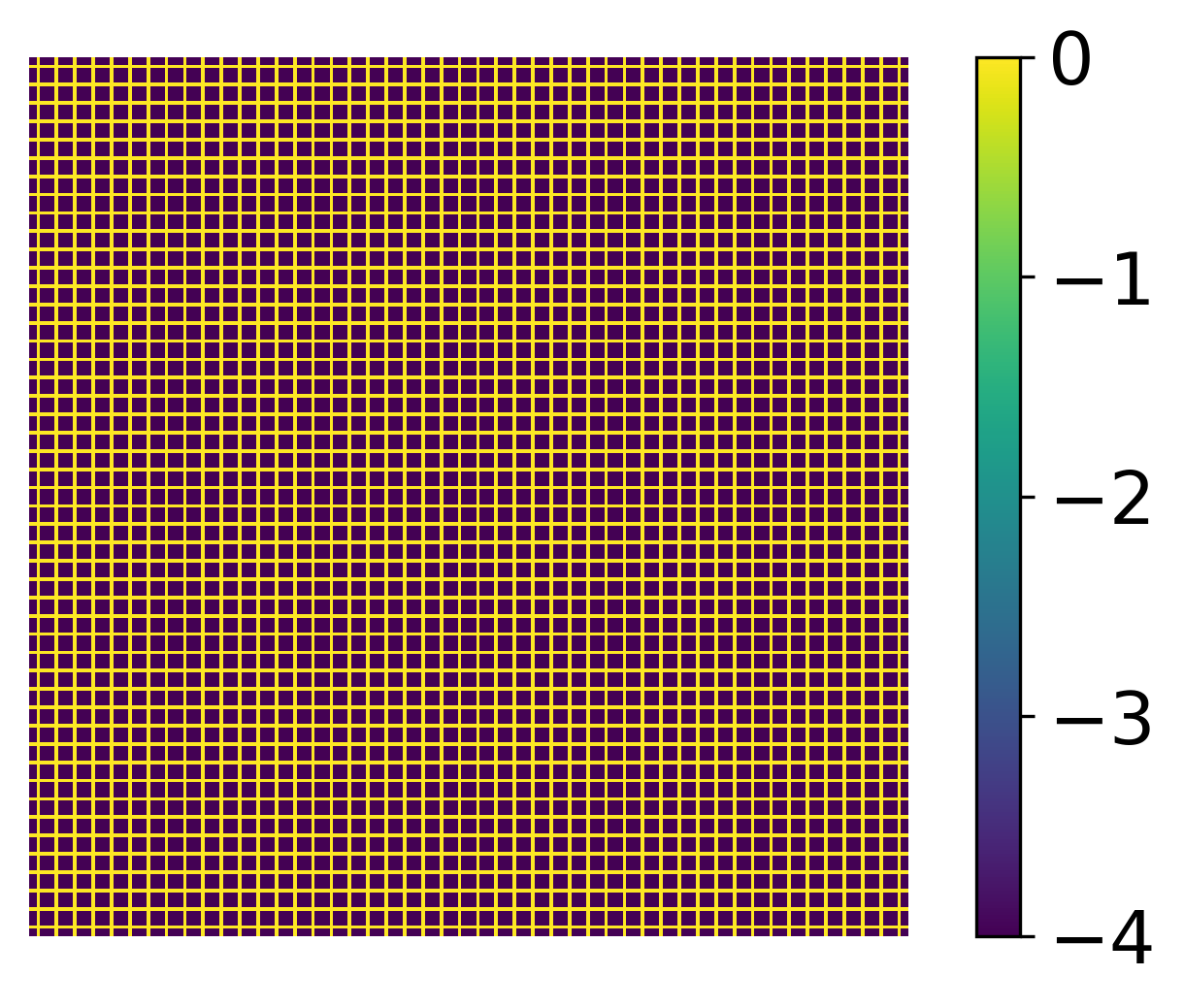}
\includegraphics[width=0.5\linewidth]{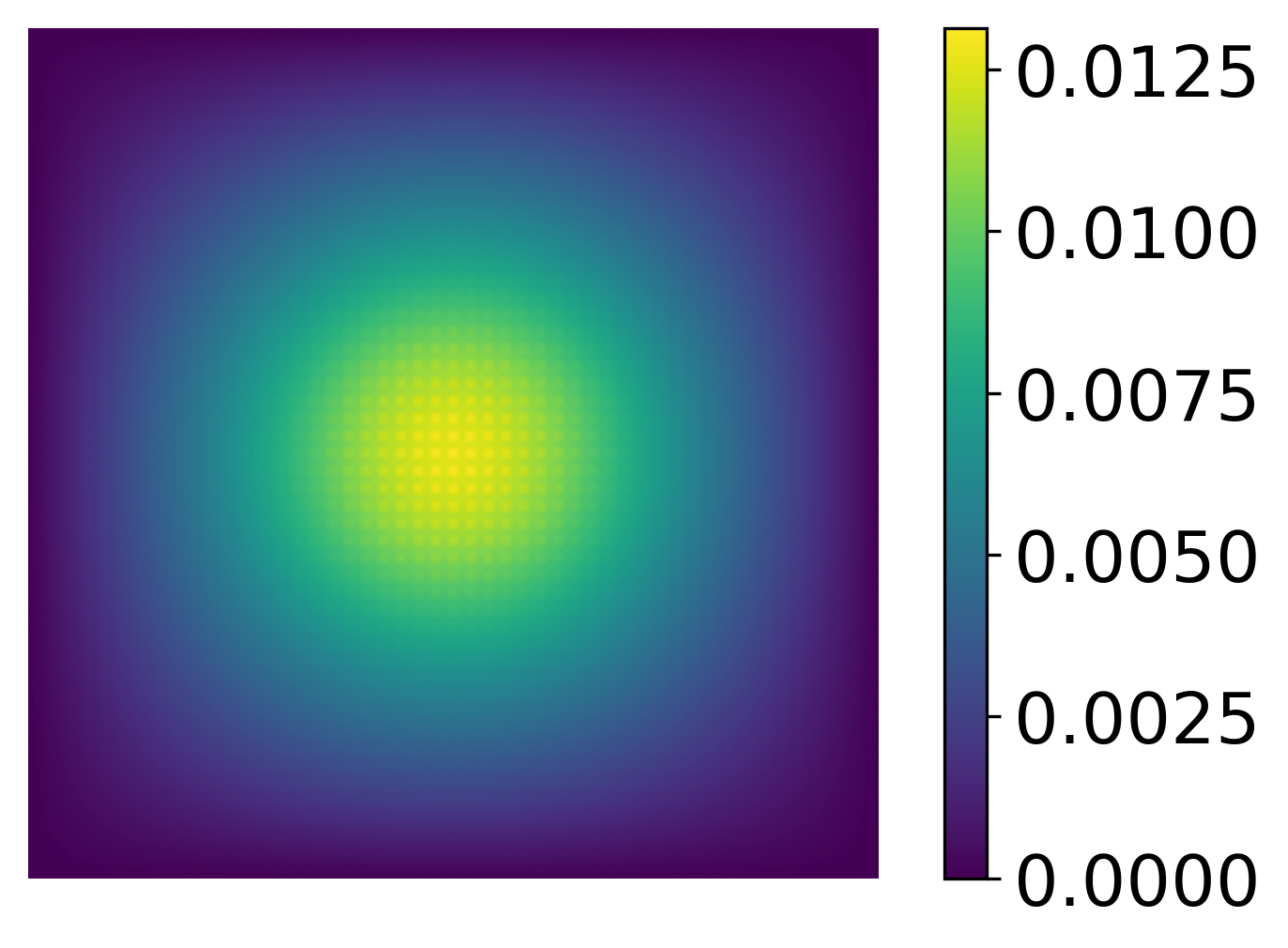}
\caption{High-contrast part $\log(k)$ (left) and reference solution (right) in Example 2.}
\label{fig:k2u}
\end{figure}

\begin{figure}[h!]
\centering
\begin{subfigure}[b]{0.48\textwidth}
\includegraphics[width=\linewidth]{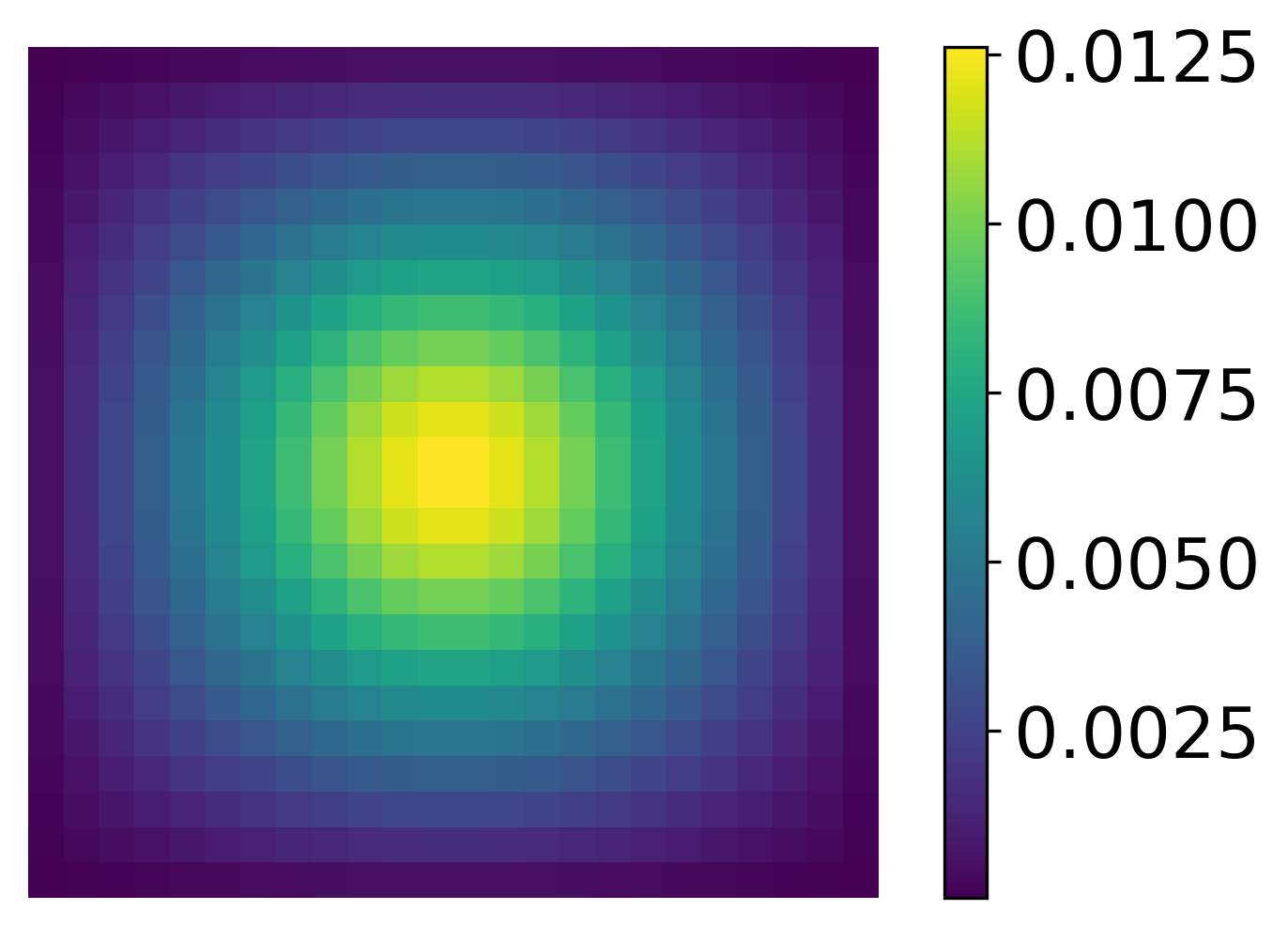}
\caption{$U_h^{(1)}$}
\end{subfigure}
\begin{subfigure}[b]{0.48\textwidth}
\includegraphics[width=\linewidth]{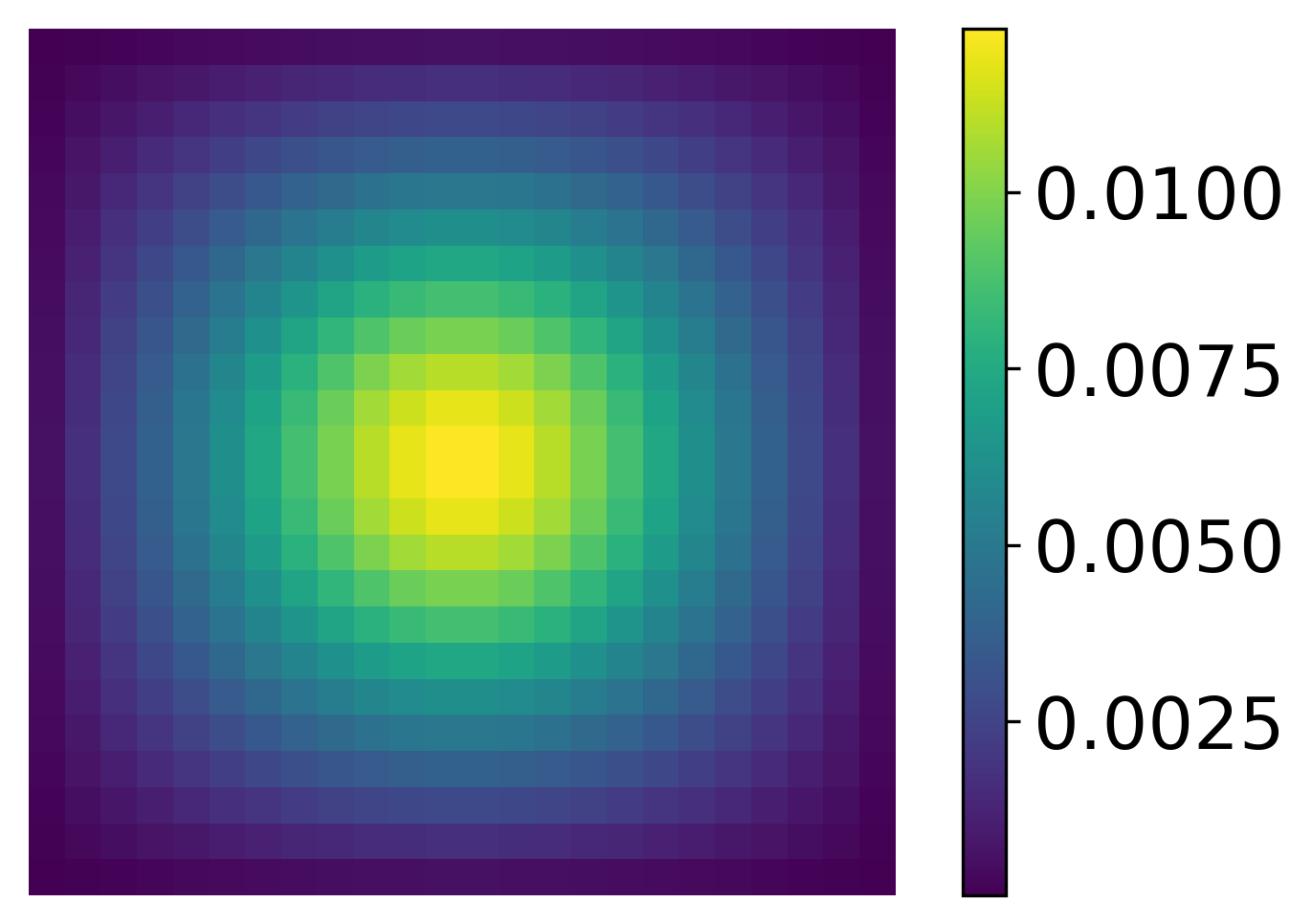}
\caption{$U_h^{(2)}$}
\end{subfigure}
\begin{subfigure}[b]{0.48\textwidth}
\includegraphics[width=\linewidth]{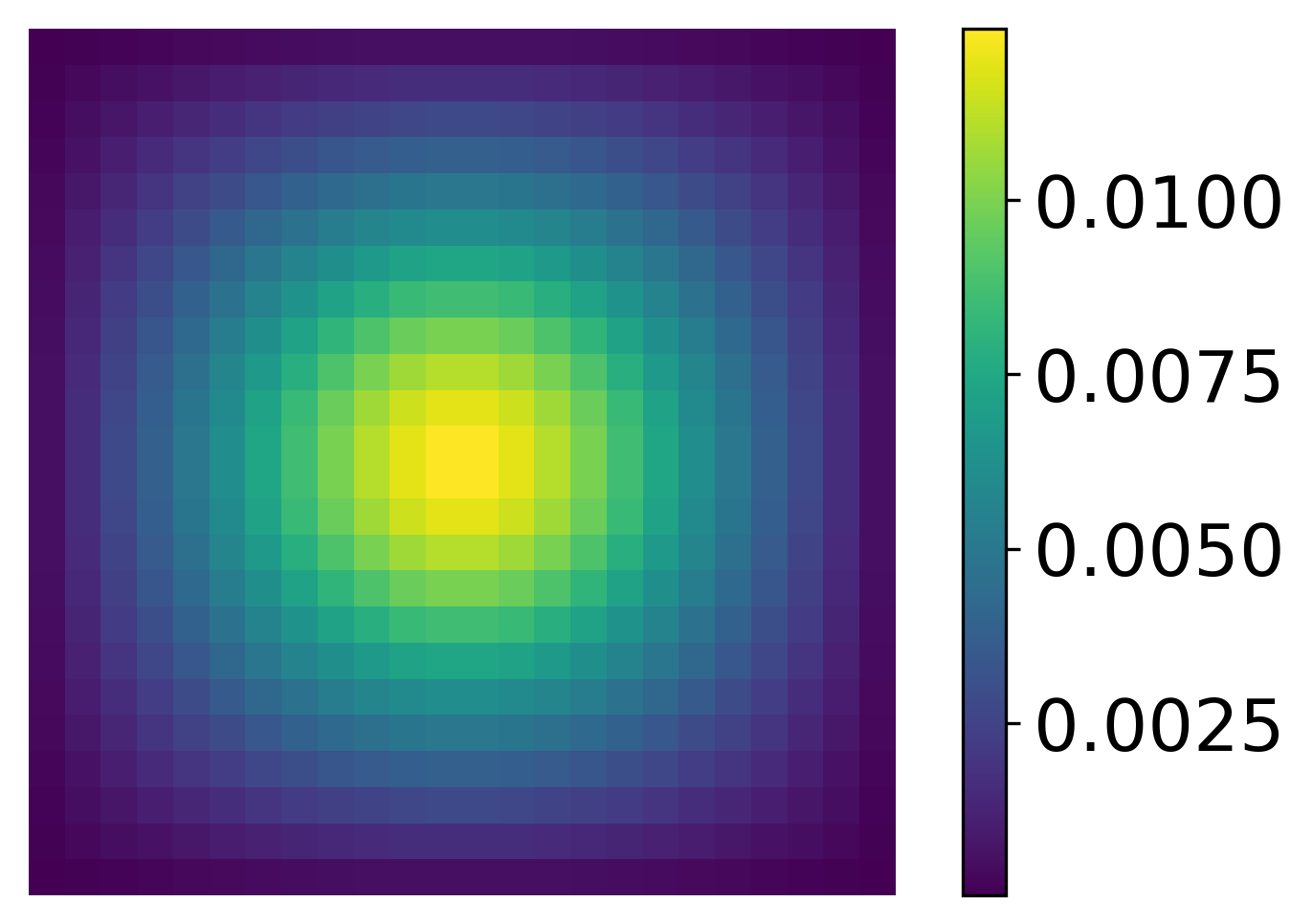}
\caption{$U_{\textnormal{MH}}^{(1)}$}
\end{subfigure}
\begin{subfigure}[b]{0.48\textwidth}
\includegraphics[width=\linewidth]{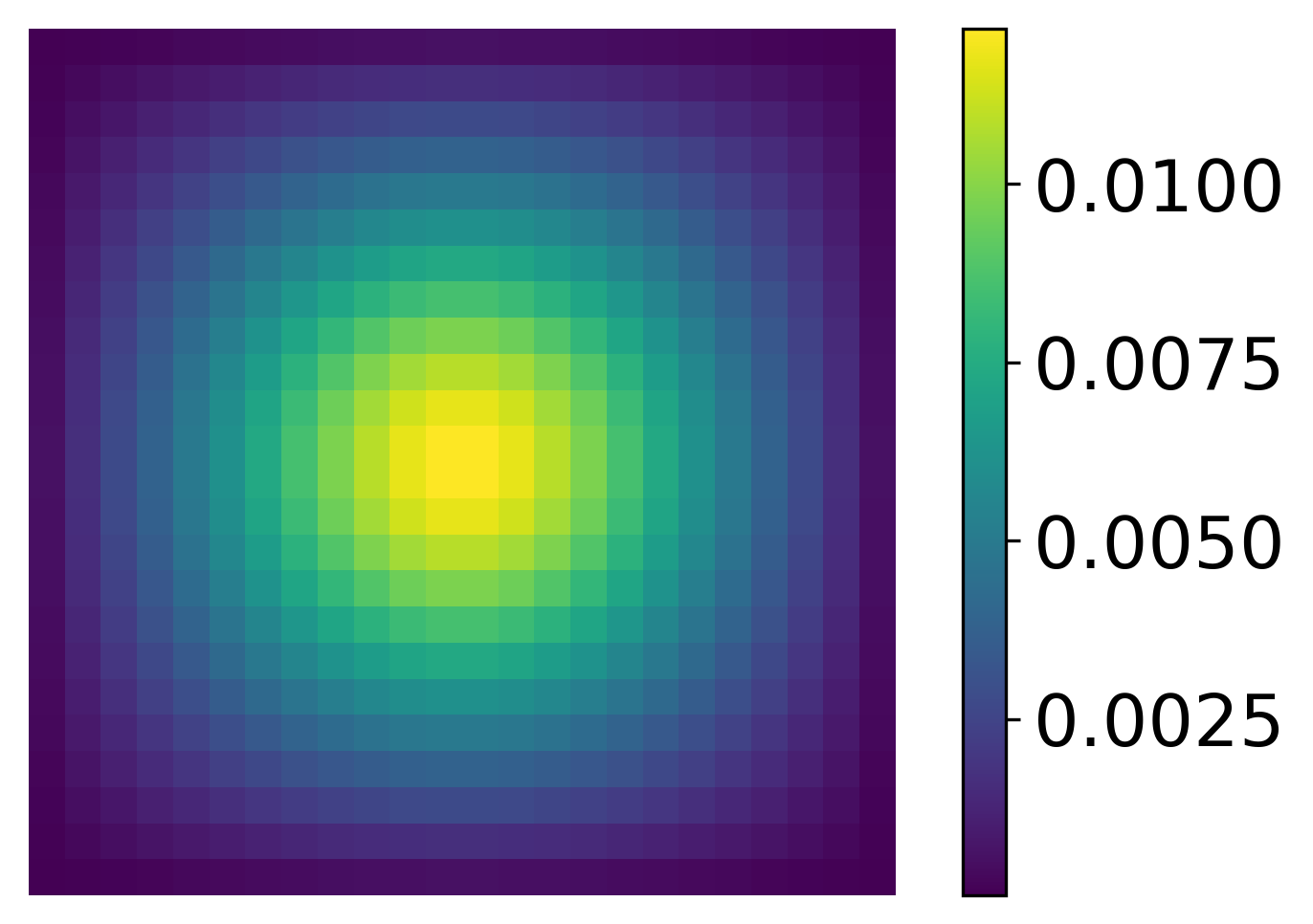}
\caption{$U_{\textnormal{MH}}^{(2)}$}
\end{subfigure}
\begin{subfigure}[b]{0.48\textwidth}
\includegraphics[width=\linewidth]{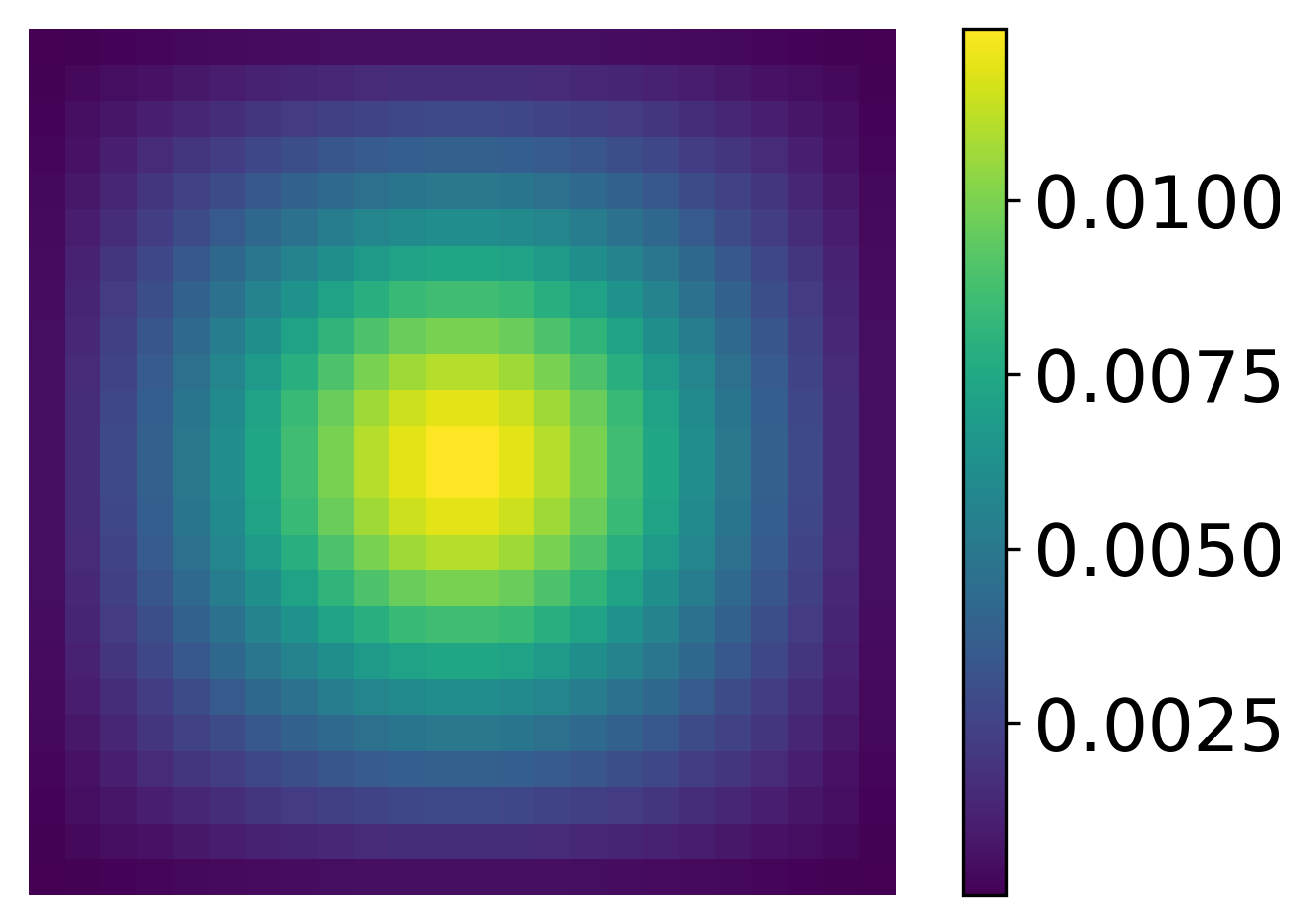}
\caption{$U_{\textnormal{HMH}}^{(1)}$}
\end{subfigure}
\begin{subfigure}[b]{0.48\textwidth}
\includegraphics[width=\linewidth]{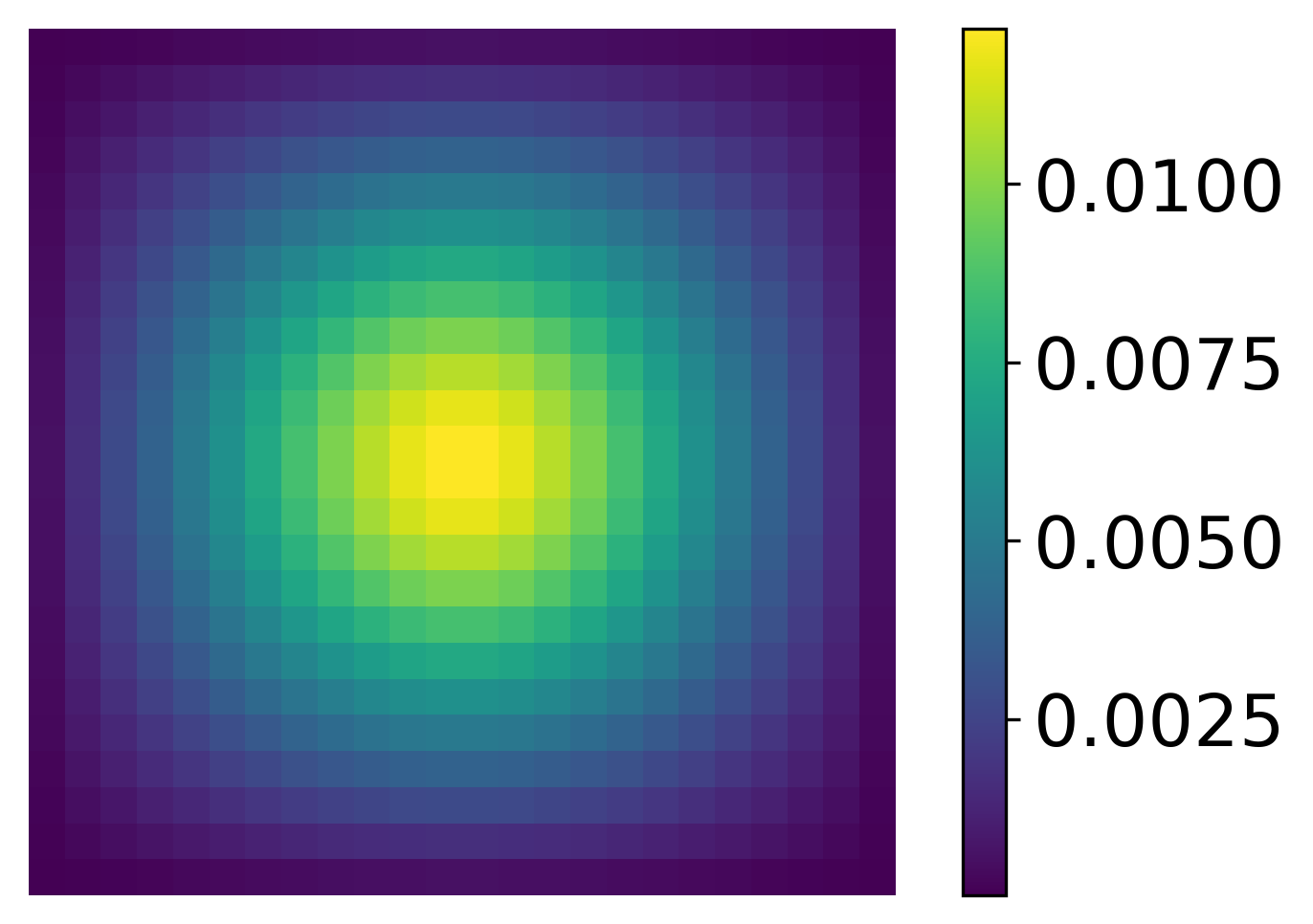}
\caption{$U_{\textnormal{HMH}}^{(2)}$}
\end{subfigure}
\caption{Average solution in Example 2.}
\label{fig:ex2_auhs}
\end{figure}

The second test involves a crossed geometry for the coefficient $k$. As in Example 1, a three-level hierarchical approach ($L=3$) is used. The logarithm of $k$ and the reference solution appear in Figure \ref{fig:k2u}, while Table \ref{tab:Herror2} presents the errors for several coarse mesh sizes. The hierarchical method again performs on par with MH, with errors decreasing systematically under refinement. The discrepancy between the two approaches is insignificant compared to the homogenization error; for instance, the homogenization error can be as large as $10^{-3}$, whereas the difference between MH and HMH is only about $10^{-7}$ for $H=1/48$. Selected average solutions are shown in Figure \ref{fig:ex2_auhs}, illustrating the consistency of the results.

\begin{table}[h!]
\renewcommand\arraystretch{1.5}
\centering
\begin{tabular}{cccccccc}
\hline
\multirow{2}*{$H$}
& \multirow{2}*{$l$}
& \multicolumn{2}{c}{Type 1} 
& \multicolumn{2}{c}{Type 2} 
& \multicolumn{2}{c}{Type 3}
\\ \cline{3-8}
& & $e_2^{(1)}$ & $e_2^{(2)}$ & $e_2^{(1)}$ & $e_2^{(2)}$ & $e_2^{(1)}$ & $e_2^{(2)}$ \\ \hline
1/12 & 5
& 4.53e-02 & 4.53e-02 & 4.53e-02 & 4.53e-02 & 6.15e-06 & 6.35e-06
\\ \hline
1/24 & 7
& 6.37e-03 & 6.23e-03 & 6.37e-03 & 6.23e-03 & 4.10e-07 & 3.77e-07
\\ \hline
1/48 & 8
& 1.47e-03 & 1.41e-03 & 1.47e-03 & 1.41e-03 & 1.59e-07 & 1.55e-07
\\ \hline
\end{tabular}
\caption{Error comparison by different coarse mesh in Example 2.}
\label{tab:Herror2}
\end{table}

\subsection{Example 3}
\begin{figure}[h!]
\centering
\includegraphics[width=0.46\linewidth]{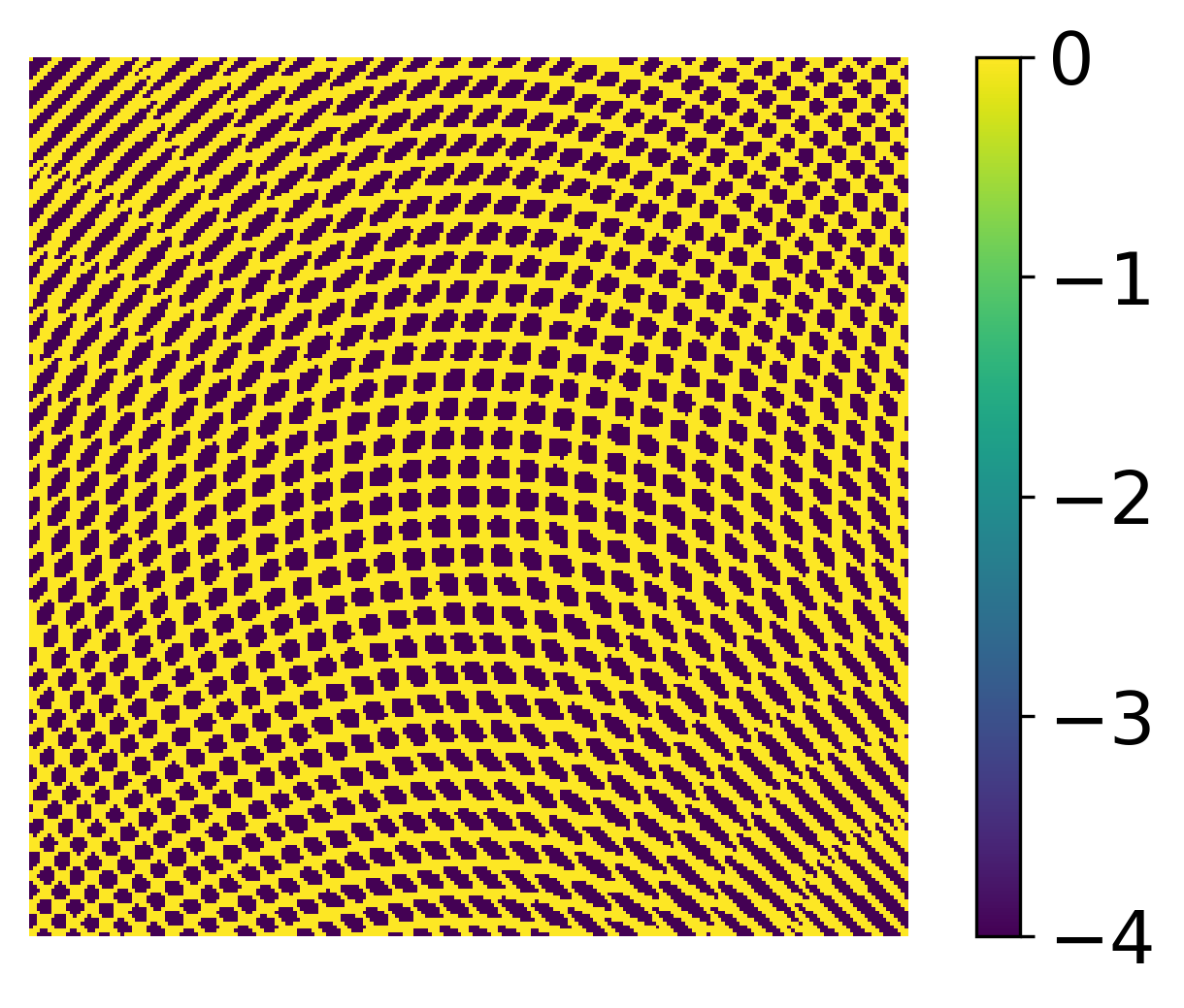}
\includegraphics[width=0.5\linewidth]{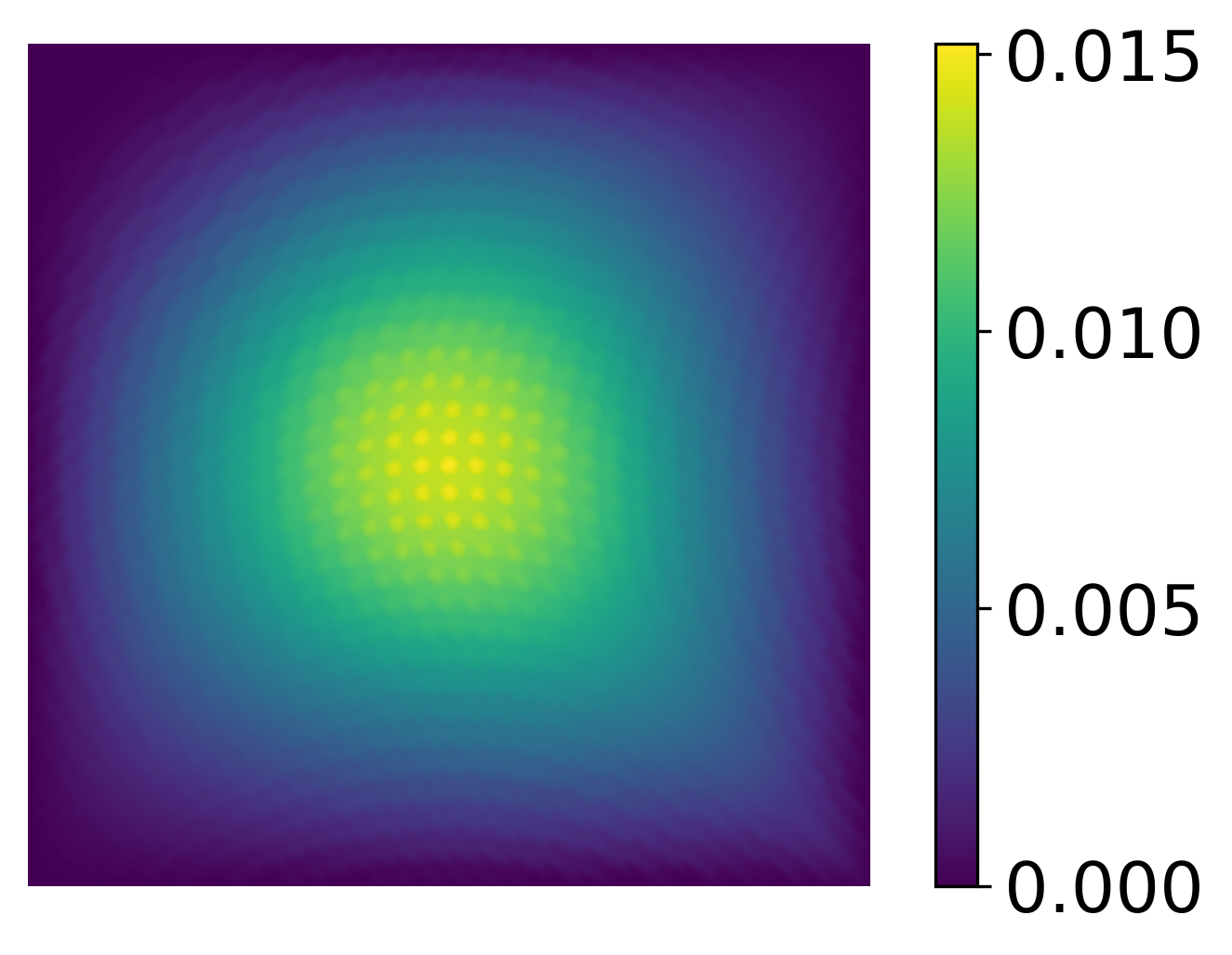}
\caption{High-contrast part $\log(k)$ (left) and reference solution (right) in Example 3.}
\label{fig:k3u}
\end{figure}

\begin{figure}[h!]
\centering
\begin{subfigure}[b]{0.48\textwidth}
\includegraphics[width=\linewidth]{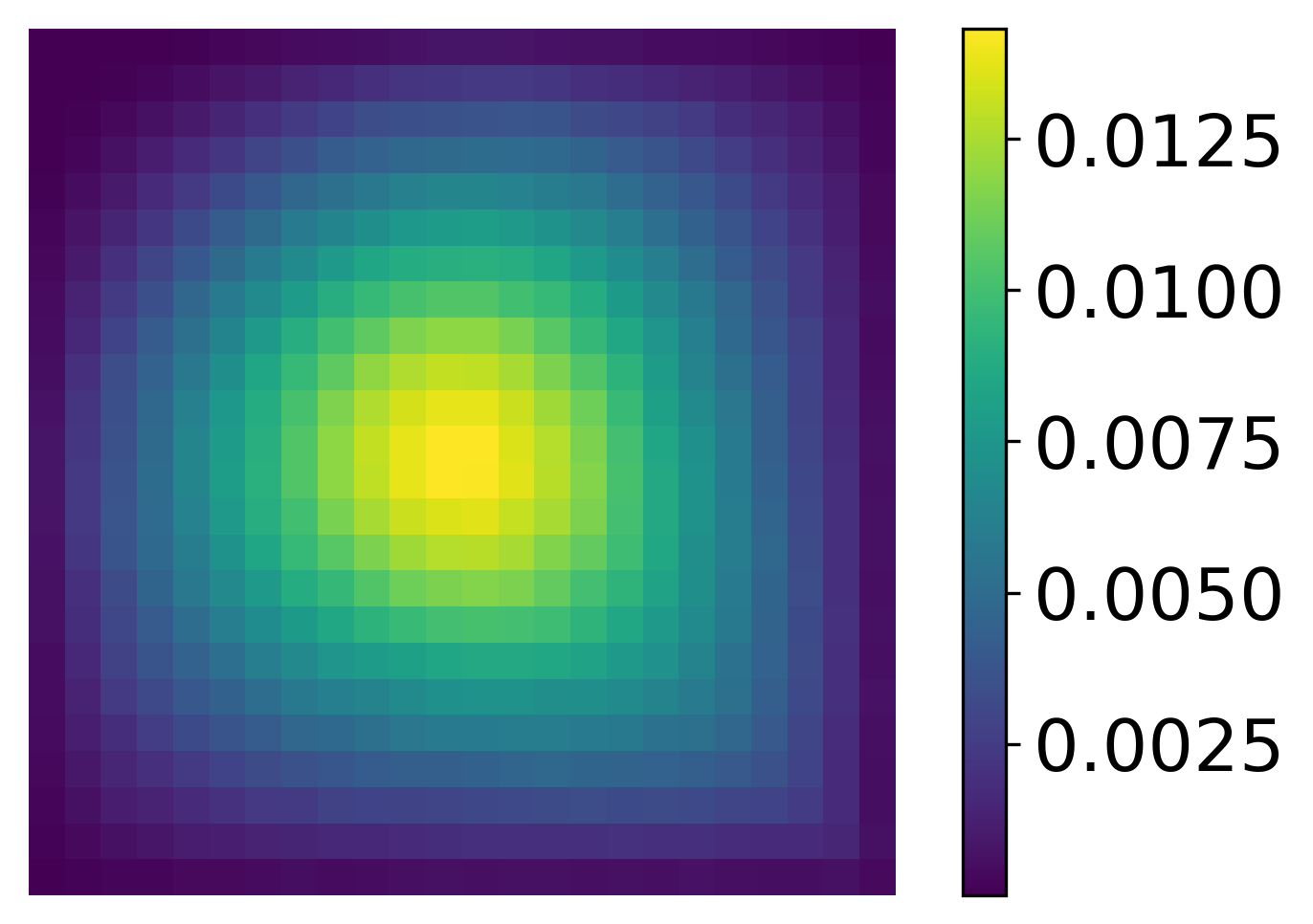}
\caption{$U_h^{(1)}$}
\end{subfigure}
\begin{subfigure}[b]{0.48\textwidth}
\includegraphics[width=\linewidth]{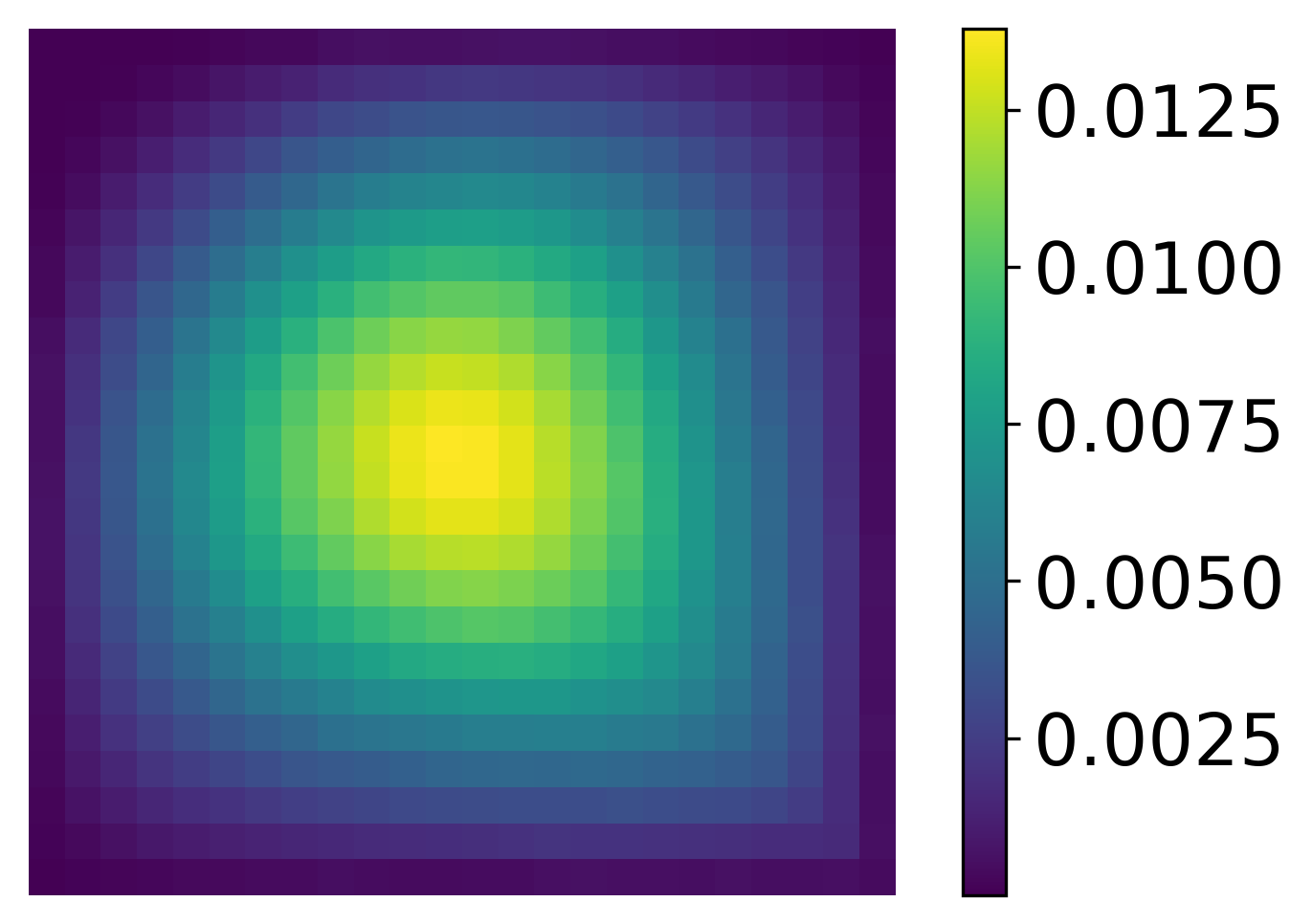}
\caption{$U_h^{(2)}$}
\end{subfigure}
\begin{subfigure}[b]{0.48\textwidth}
\includegraphics[width=\linewidth]{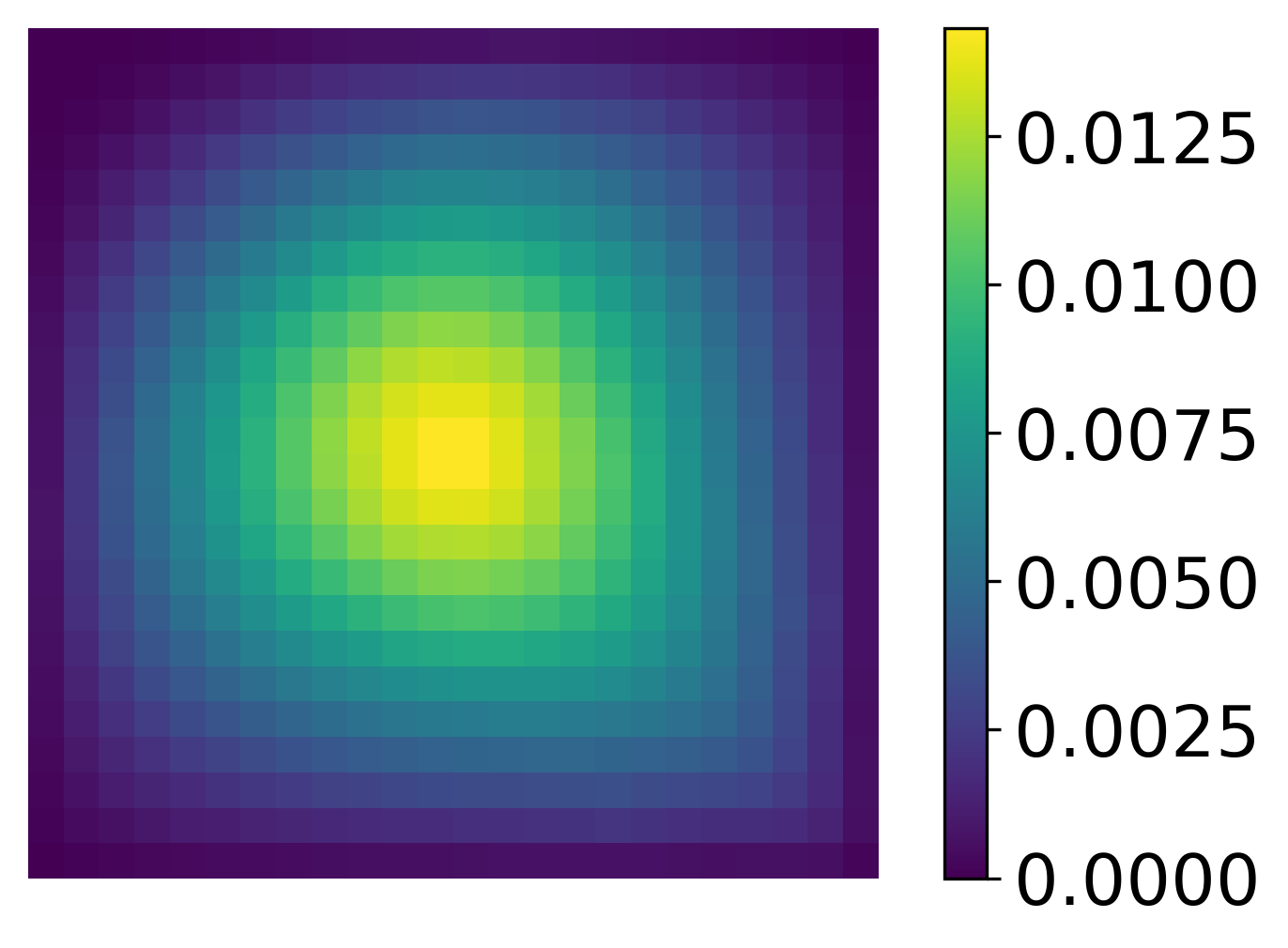}
\caption{$U_{\textnormal{MH}}^{(1)}$}
\end{subfigure}
\begin{subfigure}[b]{0.48\textwidth}
\includegraphics[width=\linewidth]{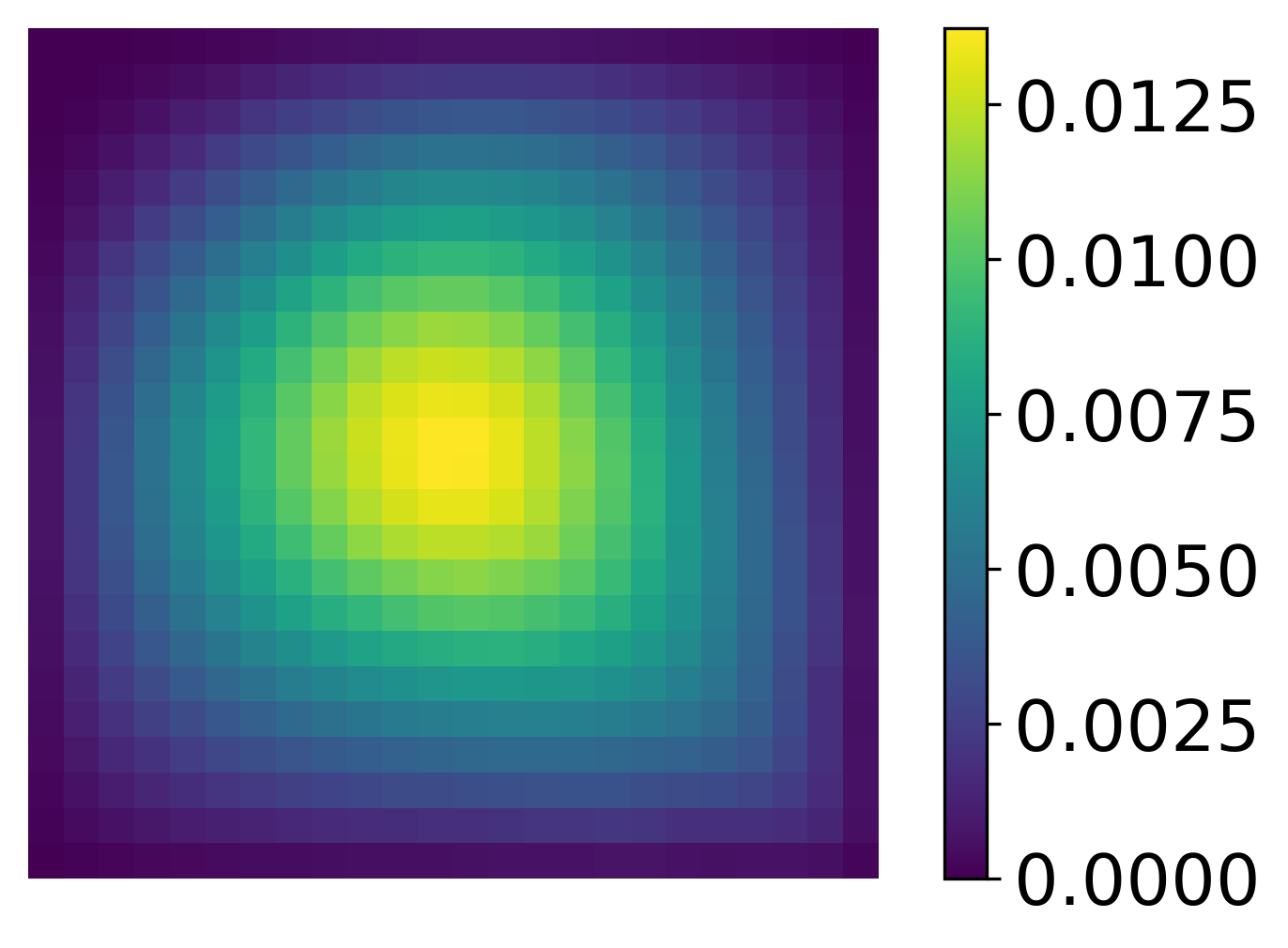}
\caption{$U_{\textnormal{MH}}^{(2)}$}
\end{subfigure}
\begin{subfigure}[b]{0.48\textwidth}
\includegraphics[width=\linewidth]{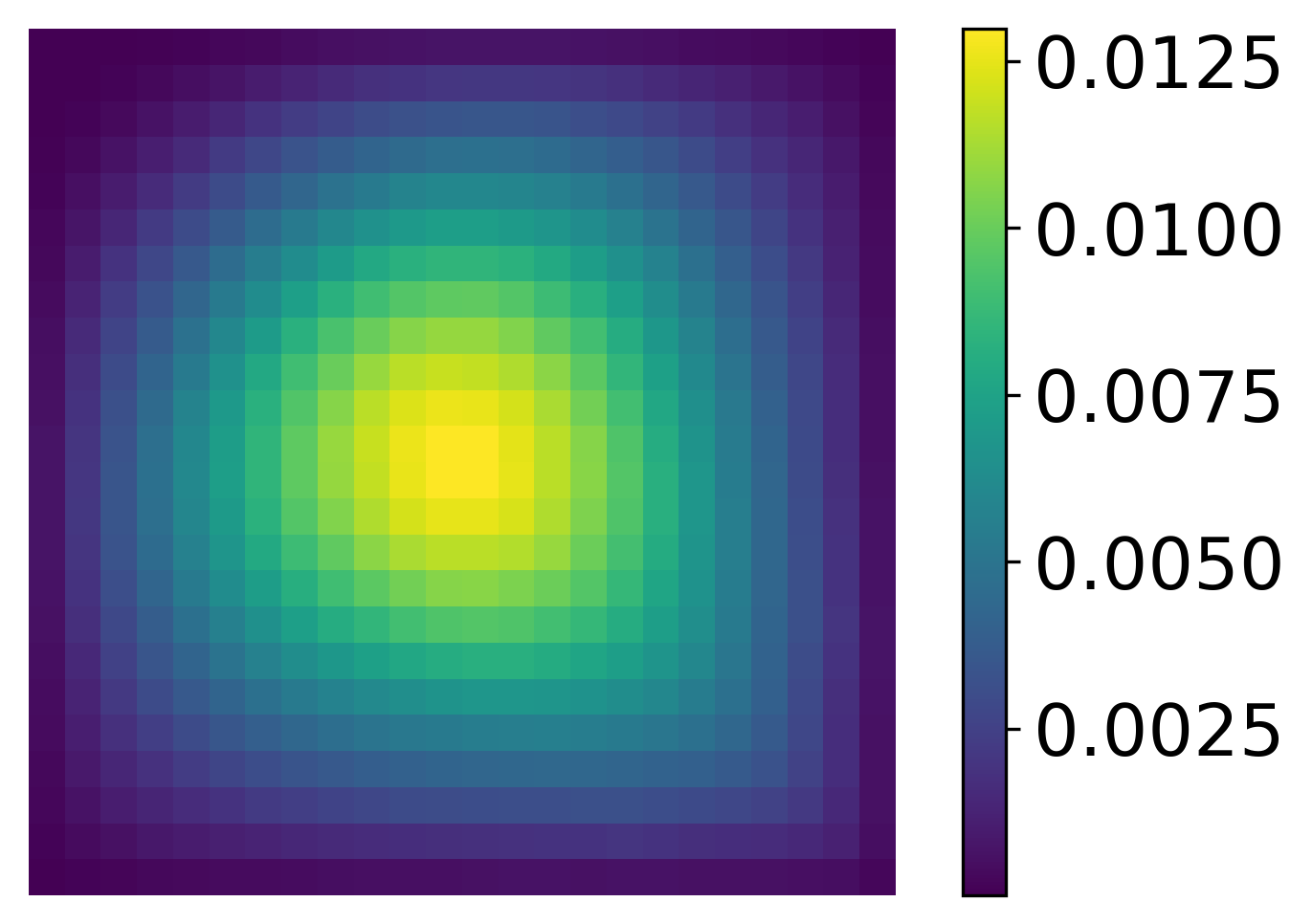}
\caption{$U_{\textnormal{HMH}}^{(1)}$}
\end{subfigure}
\begin{subfigure}[b]{0.48\textwidth}
\includegraphics[width=\linewidth]{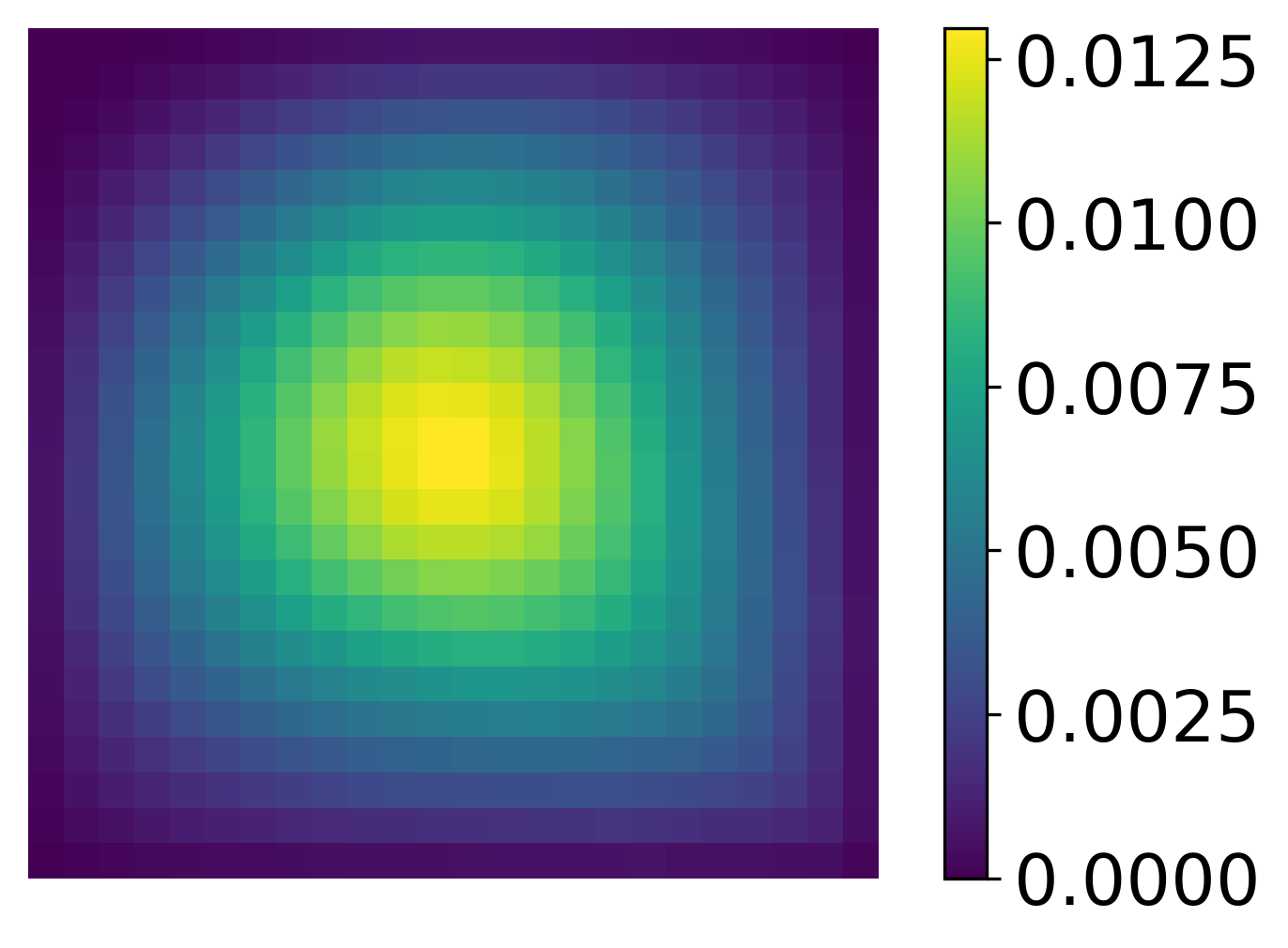}
\caption{$U_{\textnormal{HMH}}^{(2)}$}
\end{subfigure}
\caption{Average solution in Example 3.}
\label{fig:ex3_auhs}
\end{figure}

We next turn to a non-periodic high-contrast function $k$. Figure \ref{fig:k3u} depicts the logarithm of $k$ alongside the fine-grid solution. A two-level hierarchical approach ($L=2$) is adopted in this case. 
For $H=1/12$ with oversampling of $l=5$ layers, CPU timings at one macropoint and its neighbors are given in Table \ref{tab:cpu2dL2}. The hierarchical strategy reduces the computational cost by nearly a factor of two (26.42 seconds for MH vs.~12.01 seconds for HMH). Errors for two coarse mesh sizes are listed in Table \ref{tab:Herror3}. While the homogenization error decreases with mesh refinement, the hierarchical error does not exhibit the same convergence; nevertheless, it remains within a reasonable range, at roughly $10\%$ for the tested meshes. Figure \ref{fig:ex3_auhs} shows the average solutions from the reference, MH, and hierarchical methods, which align well with the numerical observations.

\begin{table}[h!]
\renewcommand\arraystretch{1.5}
\centering
\begin{tabular}{ccccc}
\hline
& Macropoints & DOF & CPU (s) \\ \hline
MH & 4 & $\mathcal{O}\left( (\frac{1}{h})^2 \right)$ & 26.42  \\ \hline
\multirow{2}*{HMH} & 1 & $\mathcal{O}\left( (\frac{1}{h})^2 \right)$ & \multirow{2}*{12.01}  \\ \cline{2-3}
 & 3 & $\mathcal{O}\left( (\frac{2}{h})^2 \right)$ &  \\ \hline
\end{tabular}
\caption{Comparison of DOF and CPU time of cell problems for MH and HMH methods in Example 3.}
\label{tab:cpu2dL2}
\end{table}

\begin{table}[h!]
\renewcommand\arraystretch{1.5}
\centering
\begin{tabular}{cccccccc}
\hline
\multirow{2}*{$H$}
& \multirow{2}*{$l$}
& \multicolumn{2}{c}{Type 1} 
& \multicolumn{2}{c}{Type 2} 
& \multicolumn{2}{c}{Type 3}
\\ \cline{3-8}
& & $e_2^{(1)}$ & $e_2^{(2)}$ & $e_2^{(1)}$ & $e_2^{(2)}$ & $e_2^{(1)}$ & $e_2^{(2)}$ \\ \hline
1/12 & 5
& 4.71e-02 & 4.69e-02 & 1.10e-01 & 9.68e-02 & 6.72e-02 & 5.44e-02
\\ \hline
1/24 & 7
& 1.28e-02 & 1.06e-02 & 7.46e-02 & 6.24e-02 & 7.73e-02 & 6.24e-02
\\ \hline
\end{tabular}
\caption{Error comparison by different coarse mesh in Example 3.}
\label{tab:Herror3}
\end{table}

\subsection{Example 4}
\begin{figure}[h!]
\centering
\includegraphics[width=0.46\linewidth]{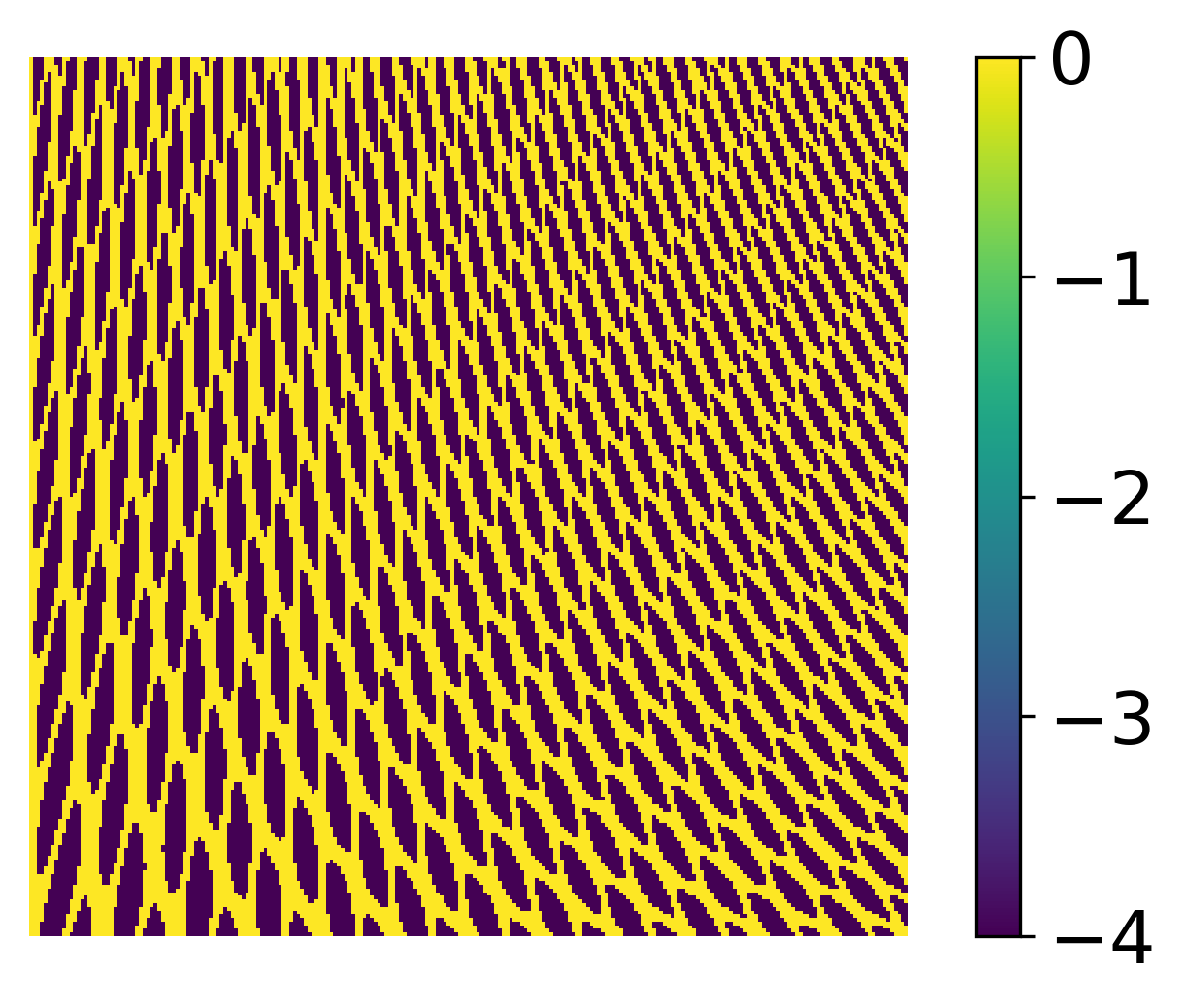}
\includegraphics[width=0.5\linewidth]{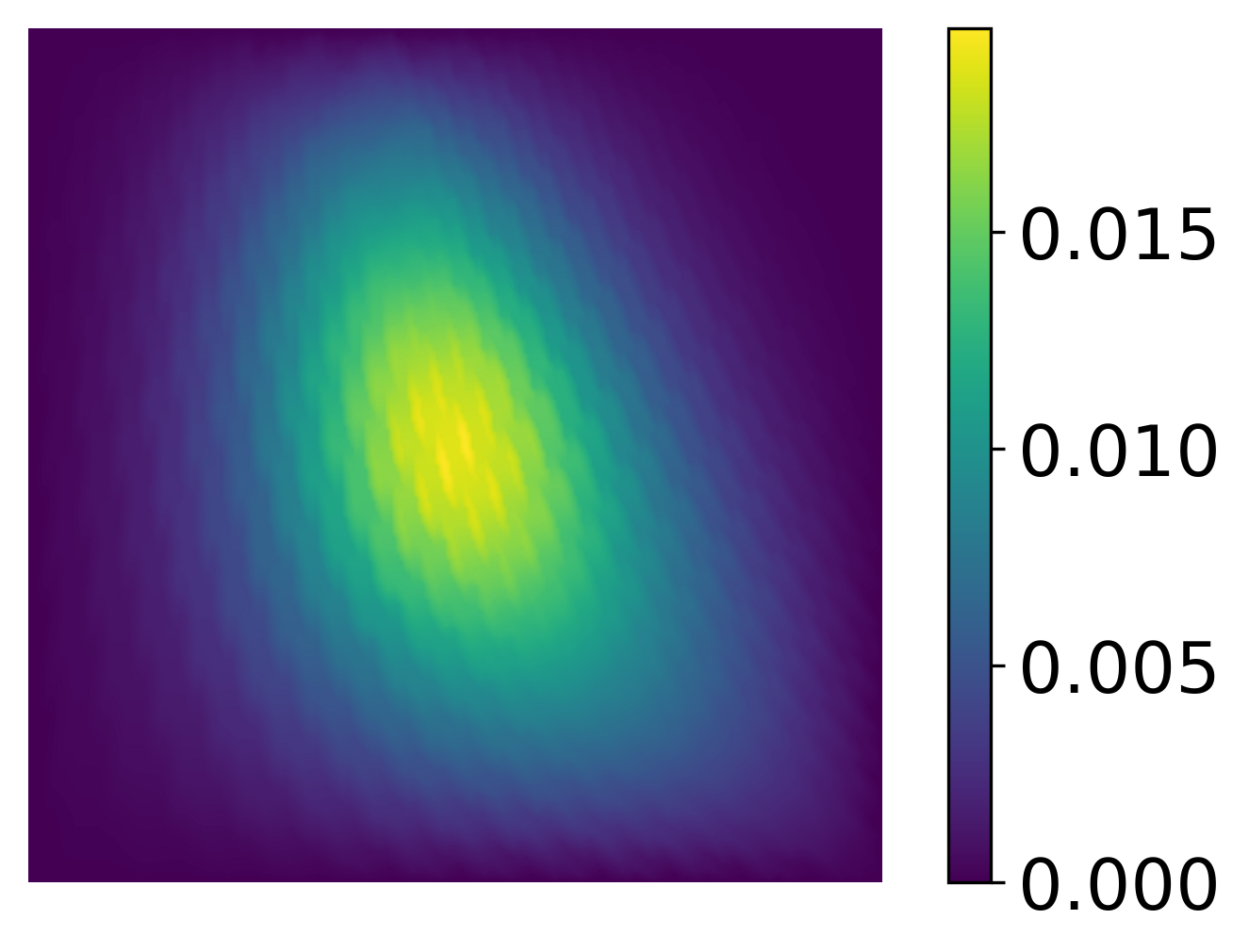}
\caption{High-contrast part $\log(k)$ (left) and reference solution (right) in Example 4.}
\label{fig:k4u}
\end{figure}

\begin{figure}[h!]
\centering
\begin{subfigure}[b]{0.48\textwidth}
\includegraphics[width=\linewidth]{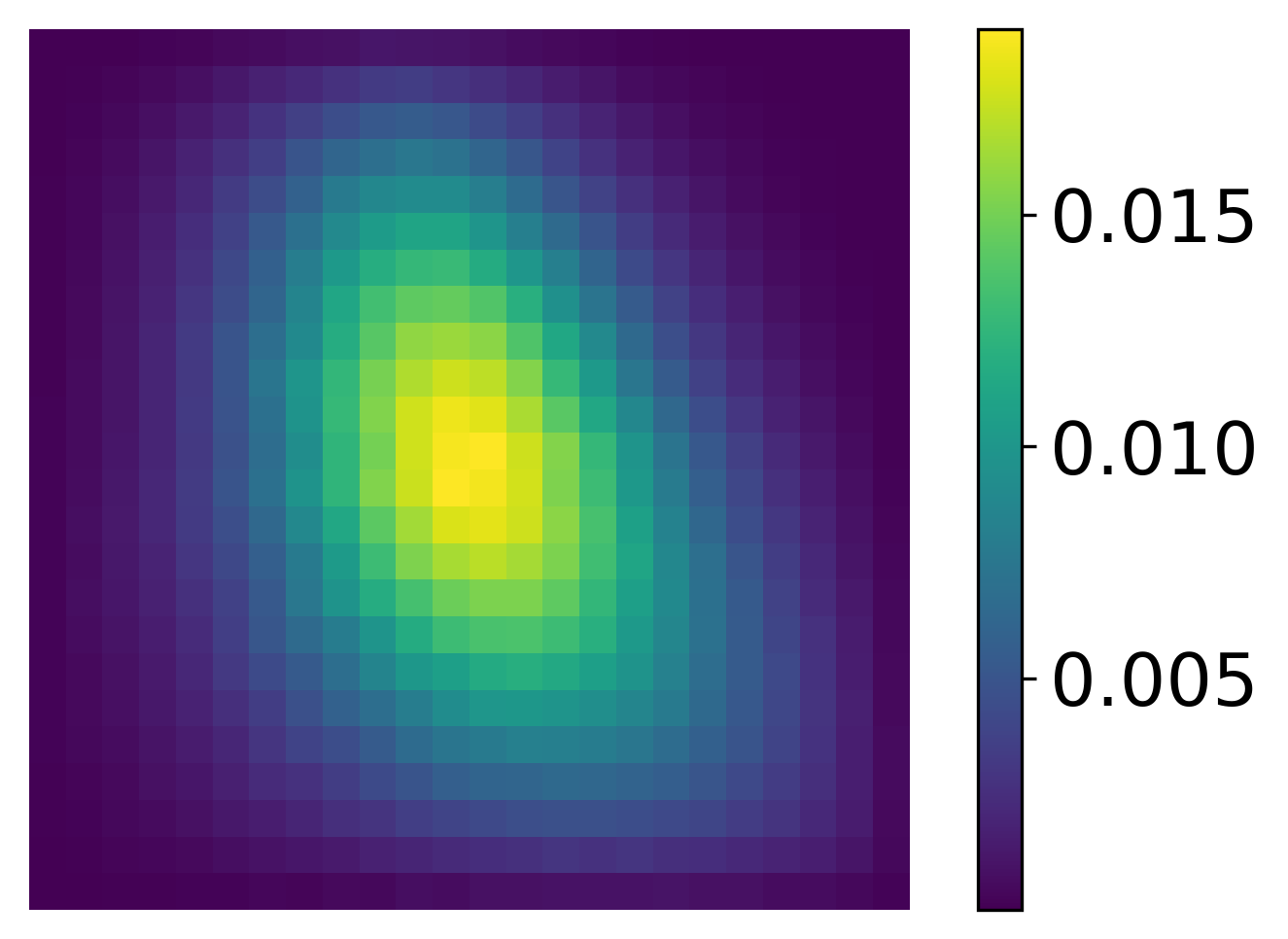}
\caption{$U_h^{(1)}$}
\end{subfigure}
\begin{subfigure}[b]{0.48\textwidth}
\includegraphics[width=\linewidth]{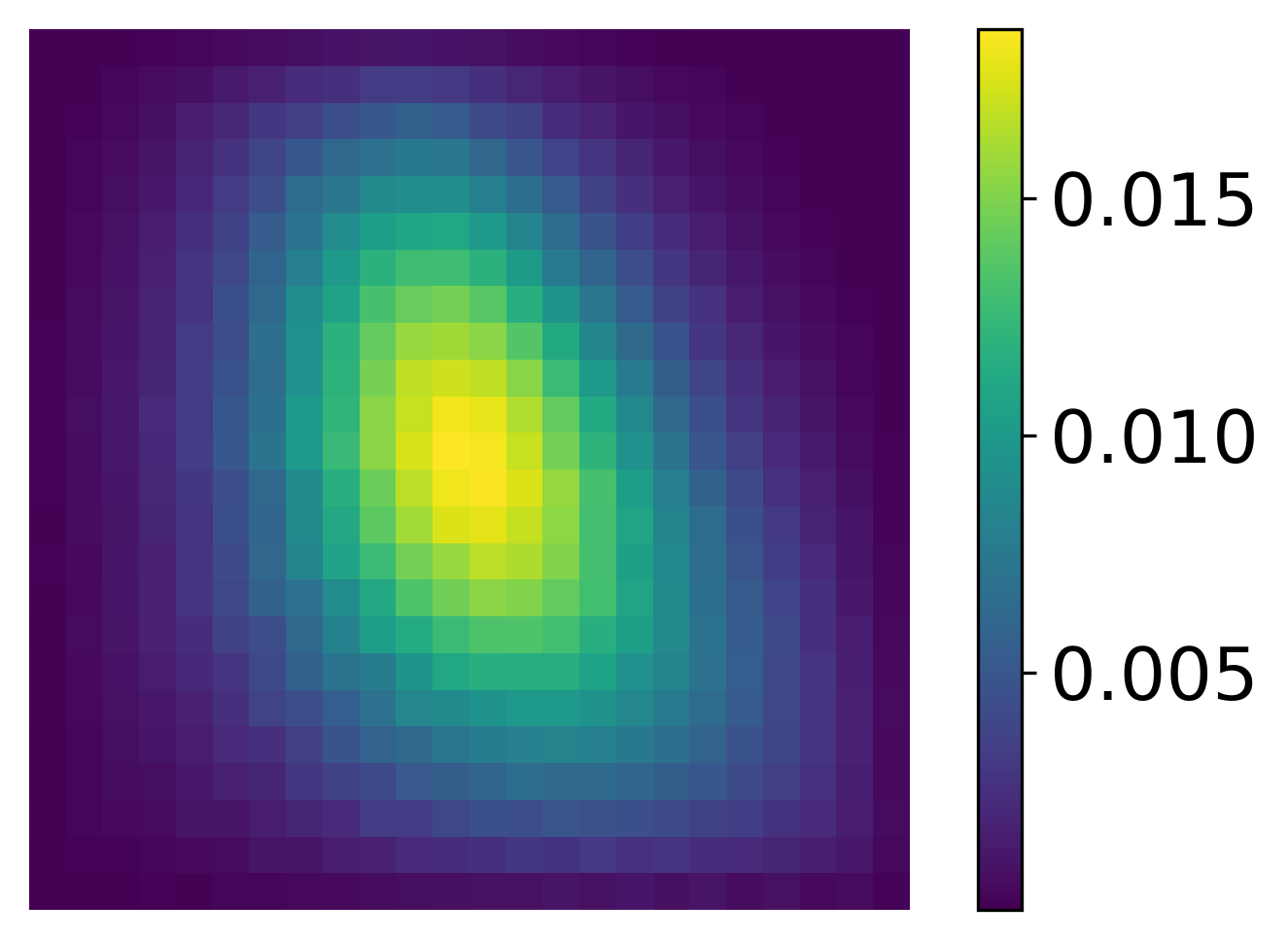}
\caption{$U_h^{(2)}$}
\end{subfigure}
\begin{subfigure}[b]{0.48\textwidth}
\includegraphics[width=\linewidth]{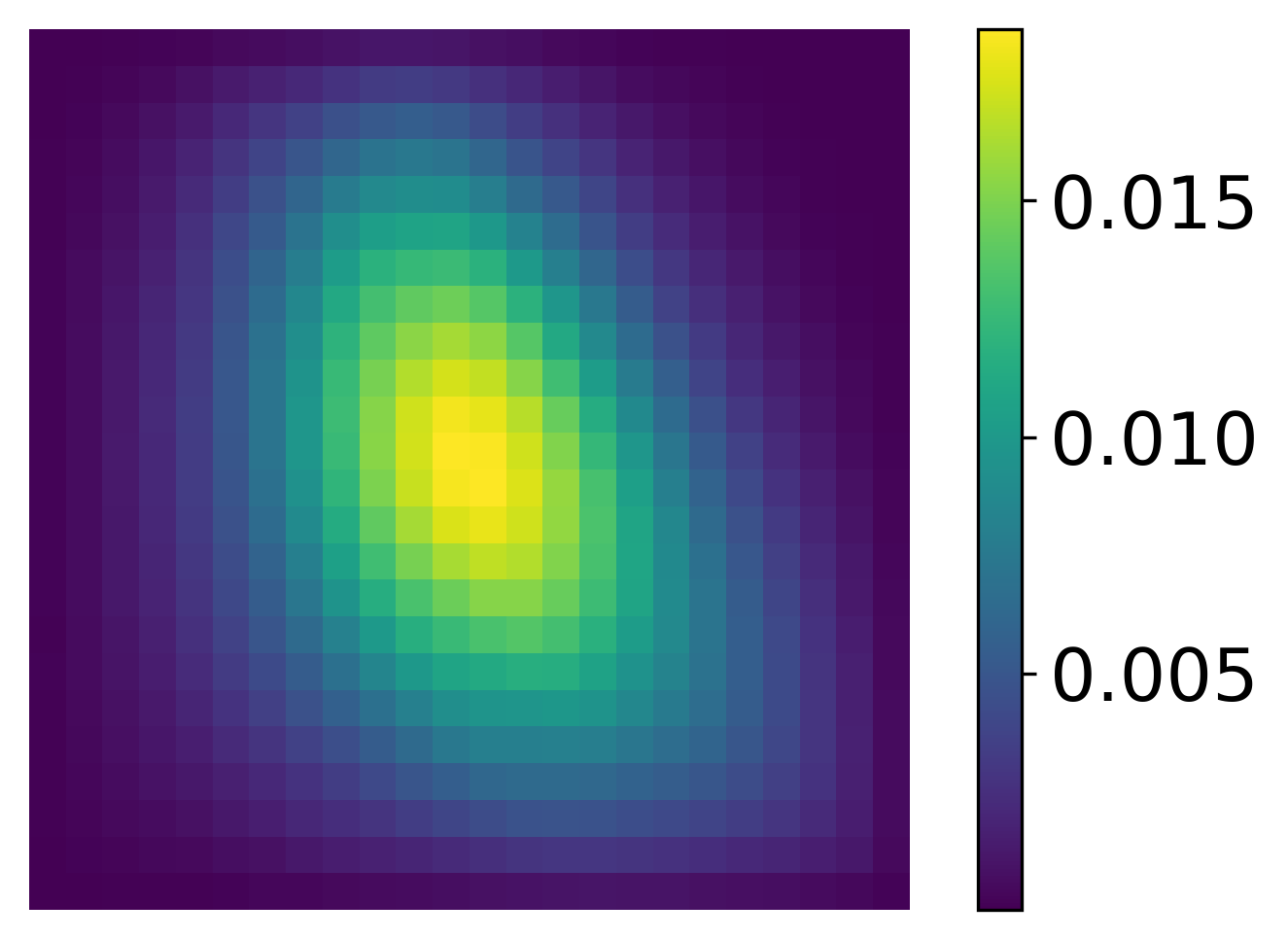}
\caption{$U_{\textnormal{MH}}^{(1)}$}
\end{subfigure}
\begin{subfigure}[b]{0.48\textwidth}
\includegraphics[width=\linewidth]{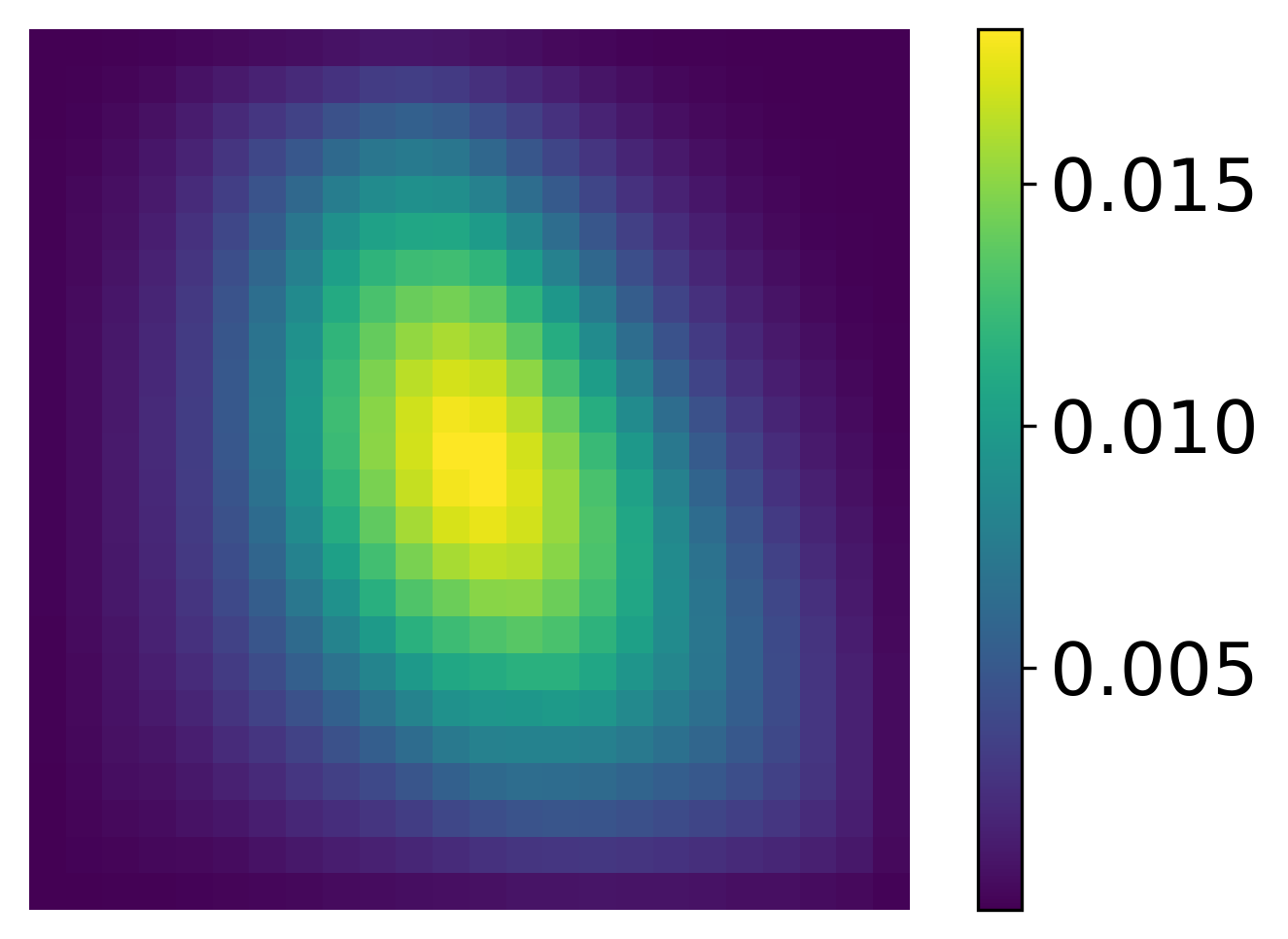}
\caption{$U_{\textnormal{MH}}^{(2)}$}
\end{subfigure}
\begin{subfigure}[b]{0.48\textwidth}
\includegraphics[width=\linewidth]{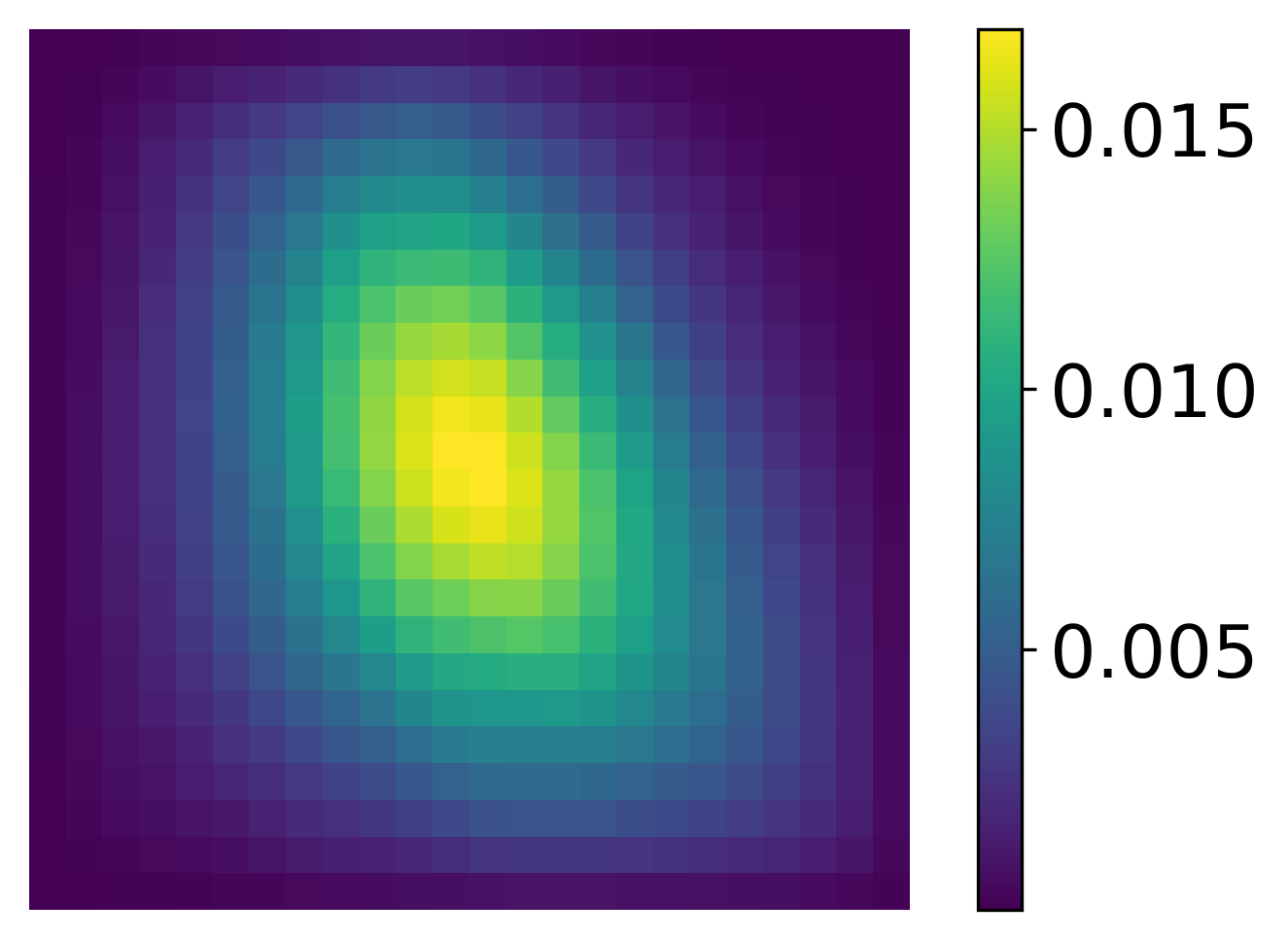}
\caption{$U_{\textnormal{HMH}}^{(1)}$}
\end{subfigure}
\begin{subfigure}[b]{0.48\textwidth}
\includegraphics[width=\linewidth]{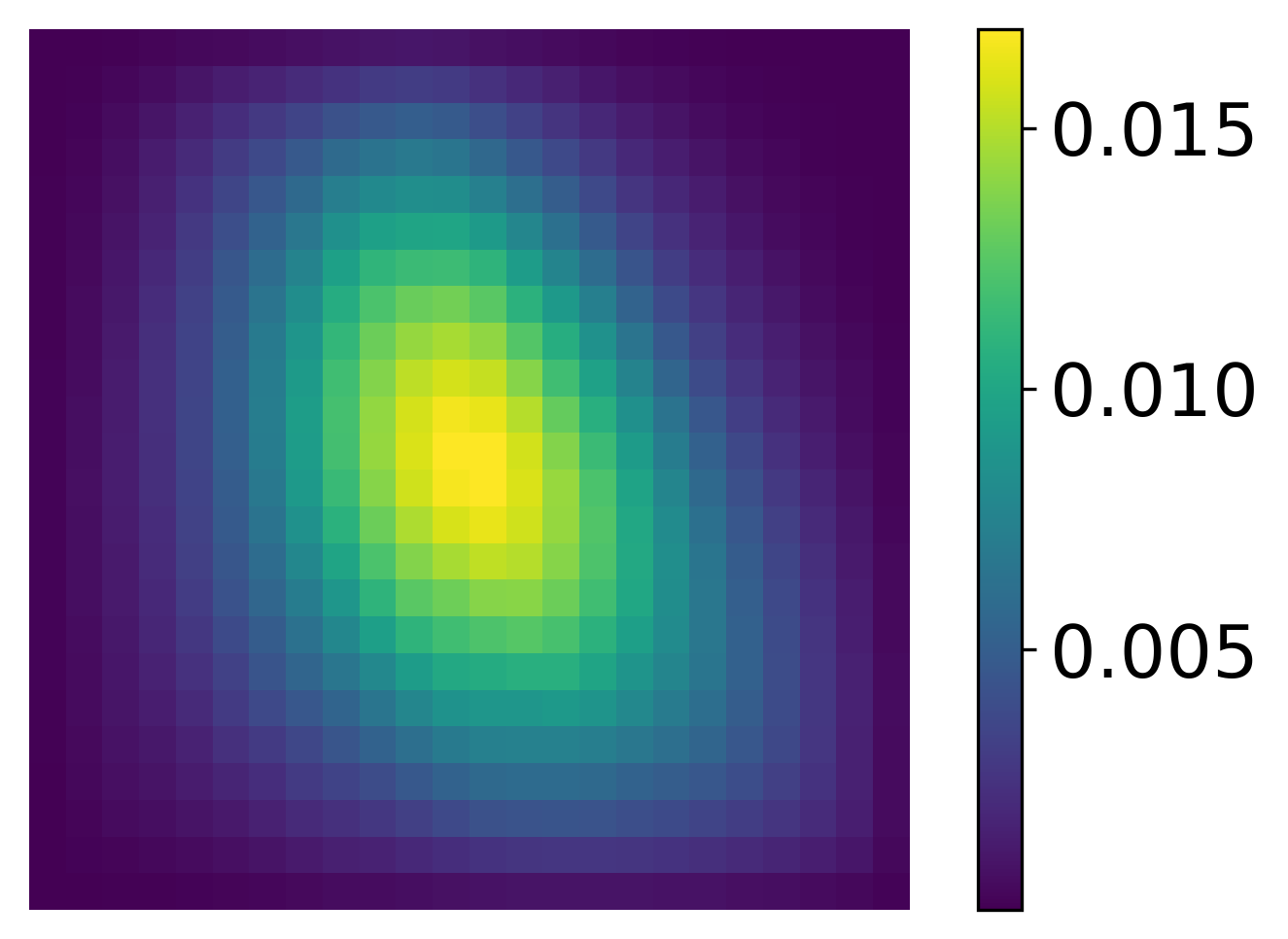}
\caption{$U_{\textnormal{HMH}}^{(2)}$}
\end{subfigure}
\caption{Average solution in Example 4.}
\label{fig:ex4_auhs}
\end{figure}

The fourth example considers another non-periodic configuration of $k$. A two-level hierarchical approach ($L=2$) is employed. Figure \ref{fig:k4u} illustrates the distribution of the two continua together with the reference solution. The errors for different coarse mesh sizes are summarized in Table \ref{tab:Herror4}. As the coarse mesh is refined, the homogenization error decreases, whereas the hierarchical error remains within an acceptable range (approximately $10\%$). Figure \ref{fig:ex4_auhs} displays the average solutions of the fine-grid, MH, and hierarchical approaches, showing consistent agreement among the three.

\begin{figure}[h!]
\centering
\includegraphics[width=0.46\linewidth]{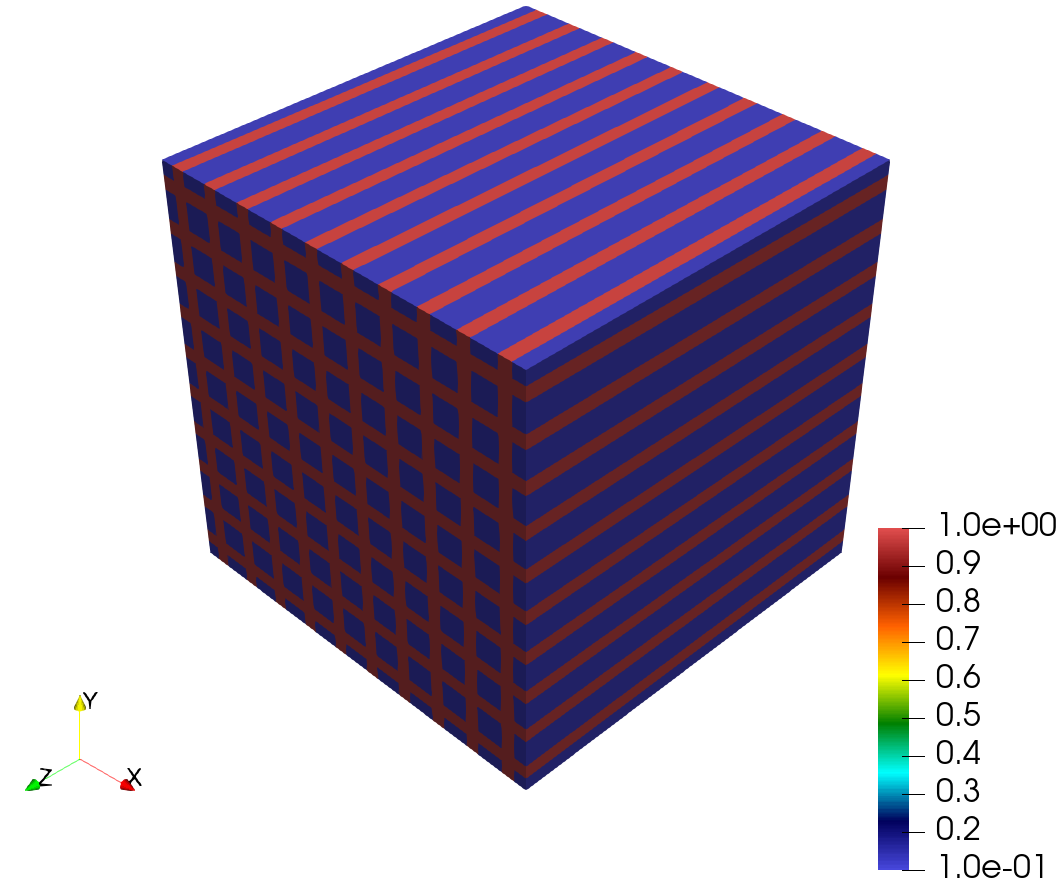}
\includegraphics[width=0.5\linewidth]{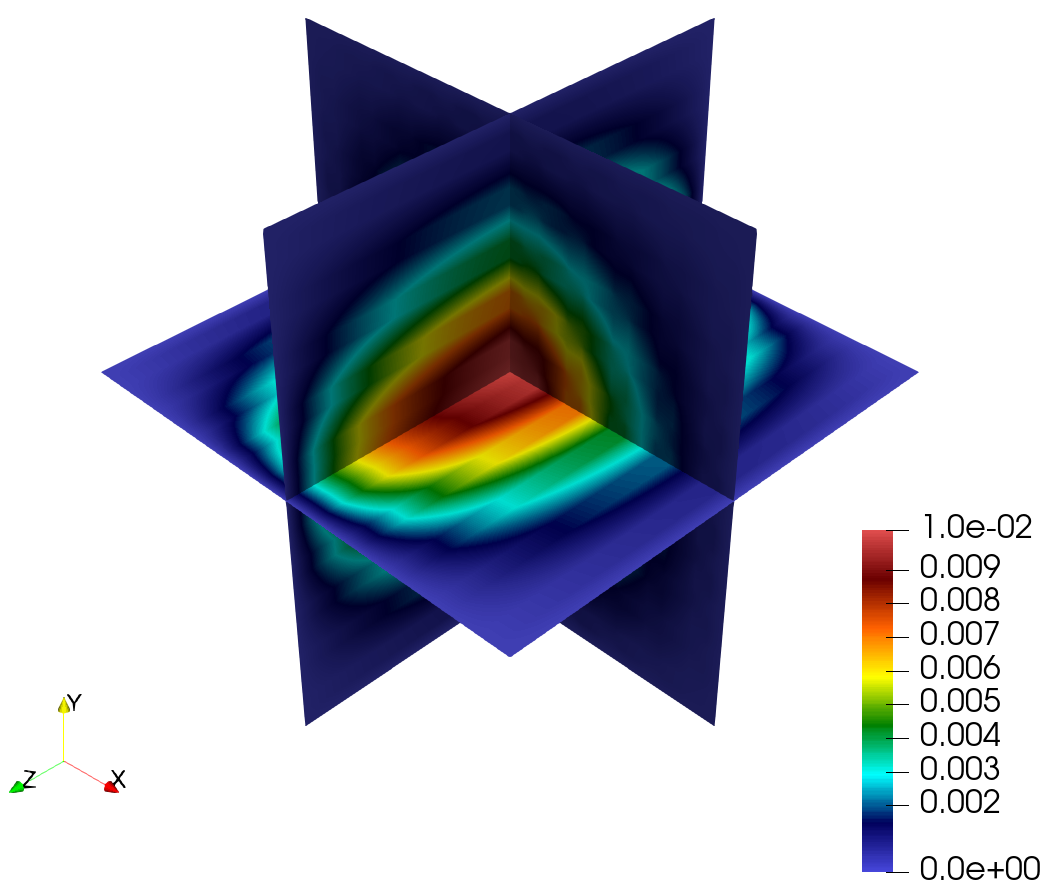}
\caption{High-contrast part $k$ (left) and reference solution (right) in Example 5.}
\label{fig:k5u}
\end{figure}

\begin{table}[h!]
\renewcommand\arraystretch{1.5}
\centering
\begin{tabular}{cccccccc}
\hline
\multirow{2}*{$H$}
& \multirow{2}*{$l$}
& \multicolumn{2}{c}{Type 1} 
& \multicolumn{2}{c}{Type 2} 
& \multicolumn{2}{c}{Type 3}
\\ \cline{3-8}
& & $e_2^{(1)}$ & $e_2^{(2)}$ & $e_2^{(1)}$ & $e_2^{(2)}$ & $e_2^{(1)}$ & $e_2^{(2)}$ \\ \hline
1/12 & 5
& 6.26e-02 & 6.41e-02 & 1.20e-01 & 1.09e-01 & 6.23e-02 & 4.89e-02
\\ \hline
1/24 & 7
& 1.92e-02 & 2.26e-02 & 8.49e-02 & 7.46e-02 & 7.66e-02 & 6.52e-02
\\ \hline
\end{tabular}
\caption{Error comparison by different coarse mesh in Example 4.}
\label{tab:Herror4}
\end{table}

\subsection{Example 5}
Finally, we test a three-dimensional case. The high-contrast function $k$ and the fine-grid solution are presented in Figure \ref{fig:k5u}. CPU timings for MH and HMH are listed in Table \ref{tab:cpu3dL2}, demonstrating that the hierarchical strategy reduces the cost by more than half (460.22 seconds for MH vs.~211.18 seconds for HMH). In this example we report results only for $H=1/10$, with the errors summarized in Table \ref{tab:Herror5}. The accuracy is satisfactory, with an error of about $7\%$, and the hierarchical approach achieves nearly machine precision when compared with MH. These findings confirm that the proposed method provides a reliable and efficient approximation also in three dimensions.

\begin{table}[h!]
\renewcommand\arraystretch{1.5}
\centering
\begin{tabular}{ccccc}
\hline
& Macropoints & DOF & CPU (s) \\ \hline
MH & 8 & $\mathcal{O}\left( (\frac{1}{h})^3 \right)$ & 460.22  \\ \hline
\multirow{2}*{HMH} & 1 & $\mathcal{O}\left( (\frac{1}{h})^3 \right)$ & \multirow{2}*{211.18}  \\ \cline{2-3}
 & 7 & $\mathcal{O}\left( (\frac{2}{h})^3 \right)$ &  \\ \hline
\end{tabular}
\caption{Comparison of DOF and CPU time of cell problems for MH and HMH methods in Example 5.}
\label{tab:cpu3dL2}
\end{table}

\begin{table}[h!]
\renewcommand\arraystretch{1.5}
\centering
\begin{tabular}{cccccccc}
\hline
\multirow{2}*{$H$}
& \multirow{2}*{$l$}
& \multicolumn{2}{c}{Type 1} 
& \multicolumn{2}{c}{Type 2} 
& \multicolumn{2}{c}{Type 3}
\\ \cline{3-8}
& & $e_2^{(1)}$ & $e_2^{(2)}$ & $e_2^{(1)}$ & $e_2^{(2)}$ & $e_2^{(1)}$ & $e_2^{(2)}$ \\ \hline
1/10 & 4
& 7.21e-02 & 7.82e-02 & 7.21e-02 & 7.82e-02 & 1.17e-16 & 8.43e-16
\\ \hline
\end{tabular}
\caption{Error comparison by different coarse mesh in Example 5.}
\label{tab:Herror5}
\end{table}

\begin{figure}[h!]
\centering
\begin{subfigure}[b]{0.48\textwidth}
\includegraphics[width=\linewidth]{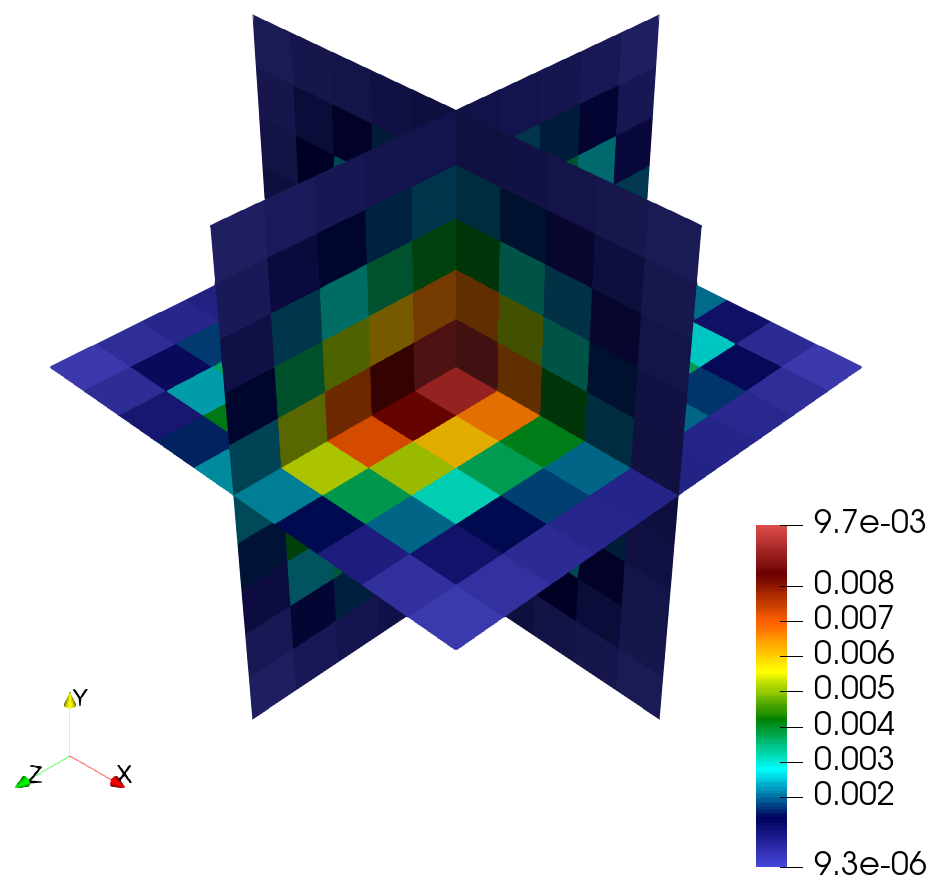}
\caption{$U_h^{(1)}$}
\end{subfigure}
\begin{subfigure}[b]{0.48\textwidth}
\includegraphics[width=\linewidth]{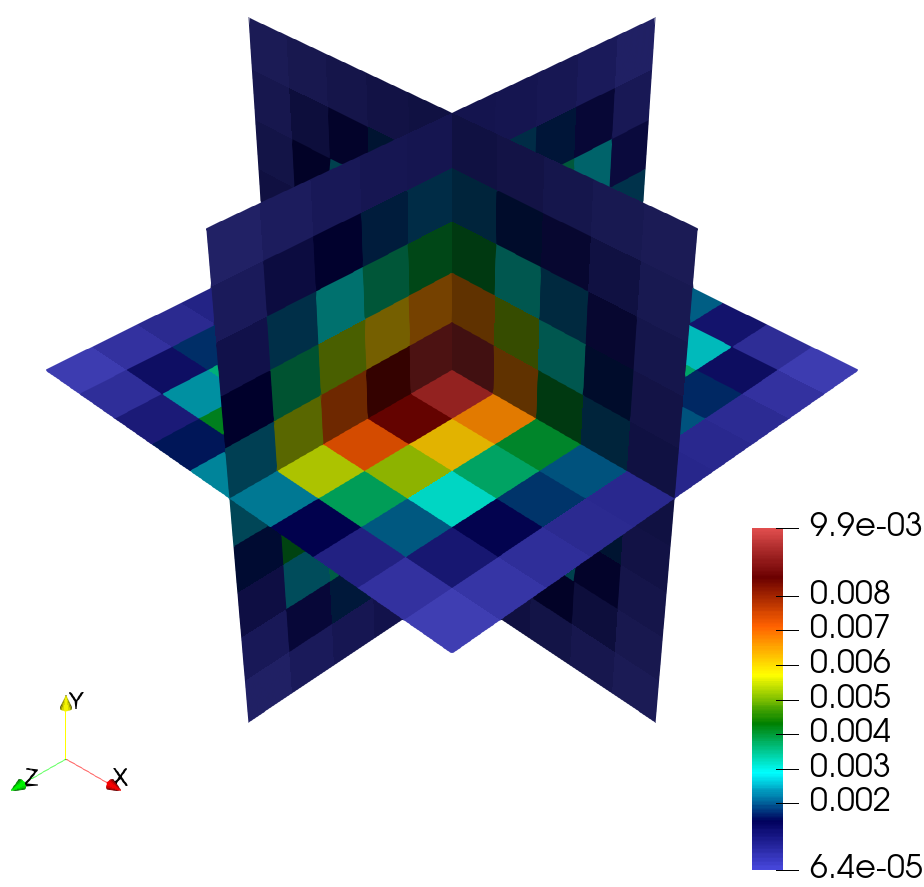}
\caption{$U_h^{(2)}$}
\end{subfigure}
\begin{subfigure}[b]{0.48\textwidth}
\includegraphics[width=\linewidth]{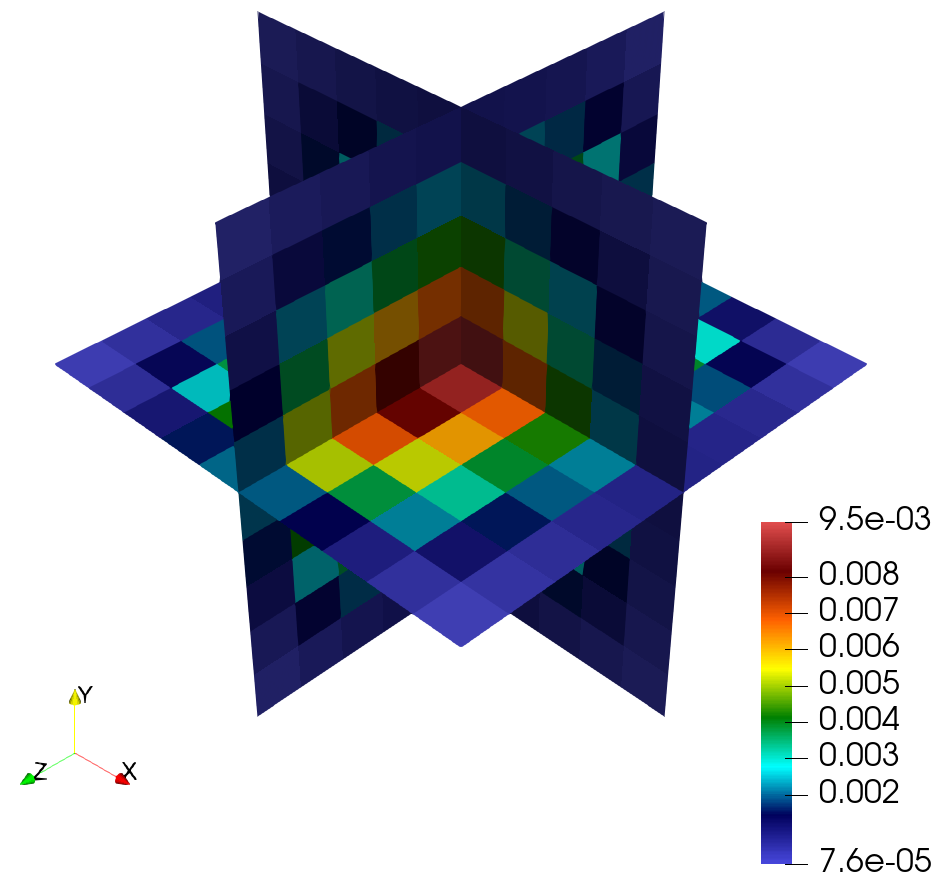}
\caption{$U_{\textnormal{MH}}^{(1)}$}
\end{subfigure}
\begin{subfigure}[b]{0.48\textwidth}
\includegraphics[width=\linewidth]{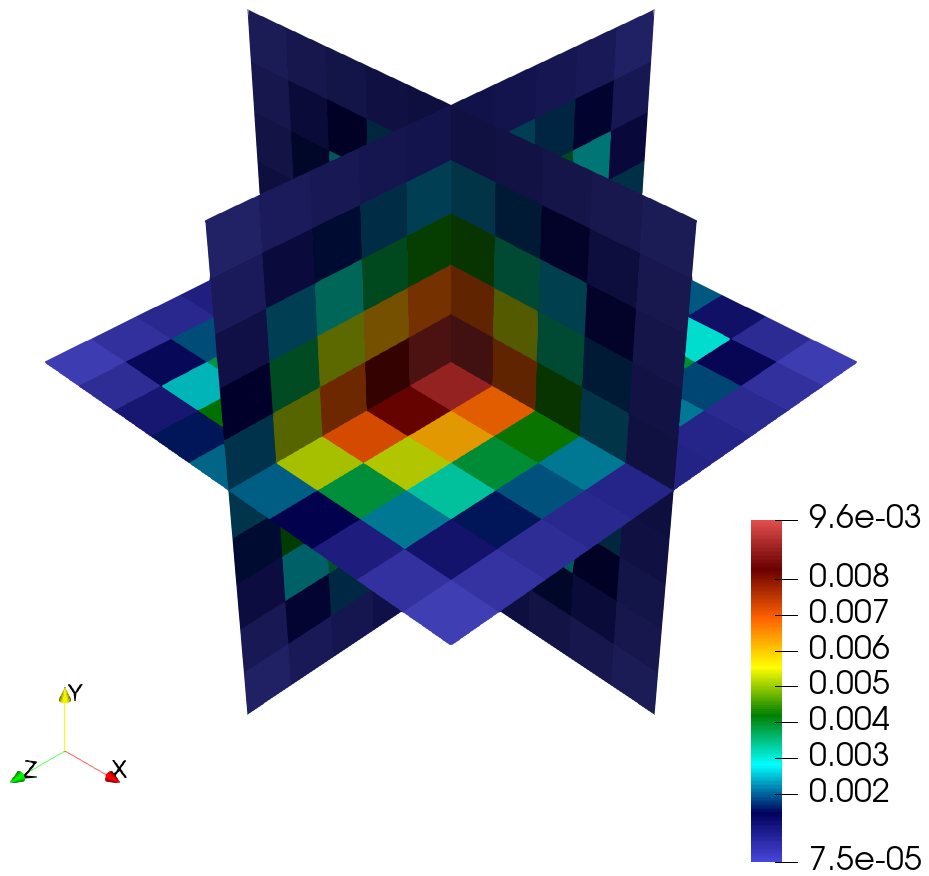}
\caption{$U_{\textnormal{MH}}^{(2)}$}
\end{subfigure}
\begin{subfigure}[b]{0.48\textwidth}
\includegraphics[width=\linewidth]{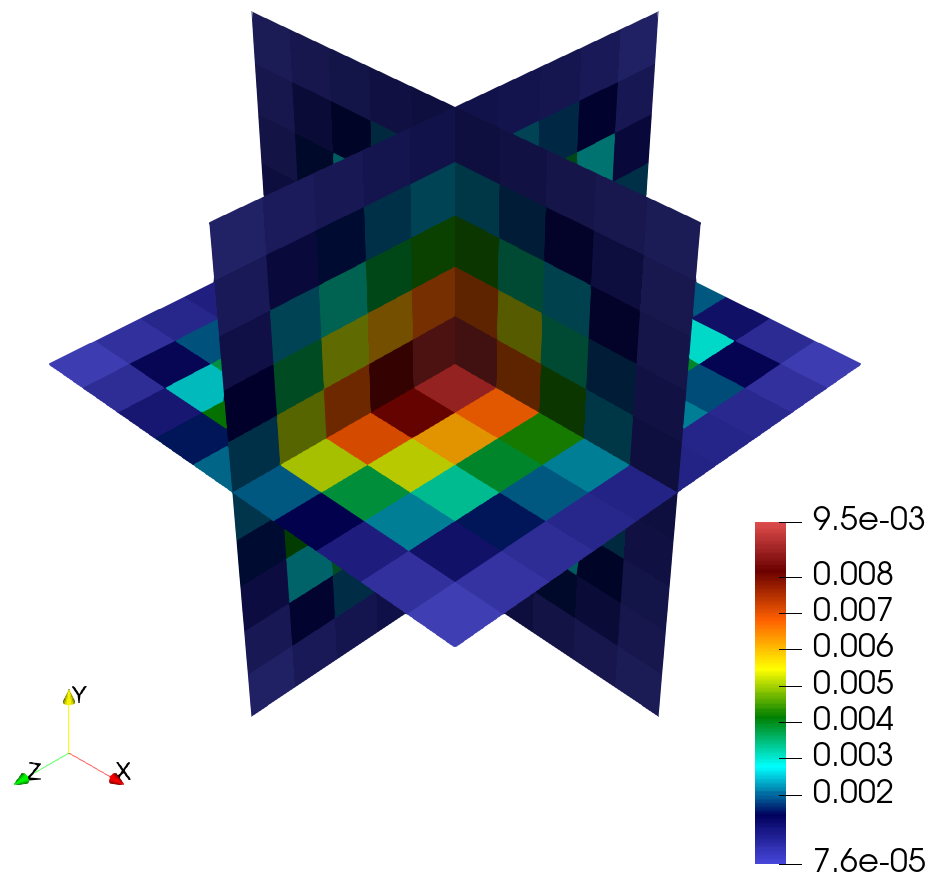}
\caption{$U_{\textnormal{HMH}}^{(1)}$}
\end{subfigure}
\begin{subfigure}[b]{0.48\textwidth}
\includegraphics[width=\linewidth]{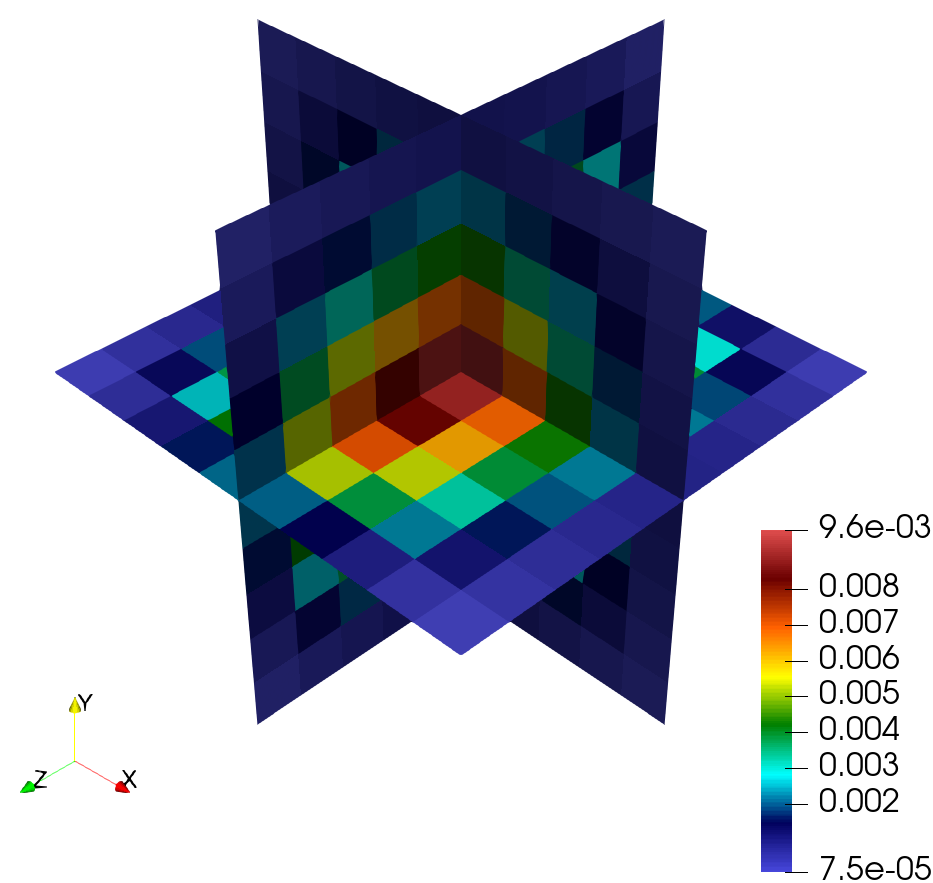}
\caption{$U_{\textnormal{HMH}}^{(2)}$}
\end{subfigure}
\caption{Average solution in Example 5.}
\label{fig:ex5_auhs}
\end{figure}

\section{Conclusions} \label{sec:conclusions}

In this paper, we propose a hierarchical algorithm for the multicontinuum homogenization method to reduce computational cost by resolving the local solutions with a very fine mesh size at each macropoint.
We outline the framework of the hierarchical multicontinuum homogenization method.
Compared to the original multicontinuum homogenization method, we assume that the local solutions can be expressed as a combination of solutions at already-constructed macropoints and a correction function.
Our approach consists of three main steps:
First, we construct a hierarchical macropoint structure with dense properties.
Second, we define a nested finite element (FE) space for the hierarchical macropoints.
Lastly, the correction terms are computed within the corresponding FE spaces.
We further analyze the computational complexity of the original multicontinuum homogenization method and our hierarchical approach.
Numerical tests show that for permeability fields without significant orders of magnitude variations, our method is fully sufficient.
Even when the permeability field exhibits large variations, our approach still achieves acceptable accuracy.

\subsection*{Acknowledgement}
Wei Xie acknowledges support from the China Scholarship Council (CSC, Project ID: 202308430231) for funding during the visit to Nanyang Technological University (NTU). Sincere thanks to NTU for providing an outstanding working and learning environment that greatly advanced the research and professional development.
Viet Ha Hoang is supported by the Tier 2 grant  T2EP20123-0047 awarded by the Singapore Ministry of Education. 
Yin Yang is supported by the grant Scientific Research In-novation Capability Support Project for Young Facult (No. SRICSPYF-BS2025132), National Natural Science Foundation of China Project (No. 12571469), the Project of Scientific Research Fund of the Hunan Provincial Science and Technology Department (No. 2024JJ1008), the 111 Project (No. D23017), and Program for Science and Technology Innovative Research Team in Higher Educational Institutions of Hunan Province of China.

\bibliographystyle{abbrv}
\bibliography{references}

@book{hornung2012homogenization,
  title={Homogenization and porous media},
  author={Hornung, Ulrich},
  volume={6},
  year={2012},
  publisher={Springer Science \& Business Media}
}

@book{cioranescu1999introduction,
  title={An introduction to homogenization},
  author={Cioranescu, Doina and Donato, Patrizia},
  year={1999},
  publisher={Oxford university press}
}

@article{hou1997multiscale,
  title={A multiscale finite element method for elliptic problems in composite materials and porous media},
  author={Hou, Thomas Y and Wu, XiaoHui},
  journal={Journal of Computational Physics},
  volume={134},
  number={1},
  pages={169--189},
  year={1997},
  publisher={Elsevier}
}

@article{hou1999convergence,
  title={Convergence of a multiscale finite element method for elliptic problems with rapidly oscillating coefficients},
  author={Hou, Thomas and Wu, XiaoHui and Cai, Zhiqiang},
  journal={Mathematics of Computation},
  volume={68},
  number={227},
  pages={913--943},
  year={1999}
}

@article{ye2020convergence,
  title={Convergence rate of multiscale finite element method for various boundary problems},
  author={Ye, Changqing and Dong, Hao and Cui, Junzhi},
  journal={Journal of Computational and Applied Mathematics},
  volume={374},
  pages={112754},
  year={2020},
  publisher={Elsevier}
}

@article{jiang2017reduced,
  title={Reduced multiscale computation on adapted grid for the convection-diffusion Robin problem},
  author={Jiang, Shan and Sun, Meiling and Yang, Yin},
  journal={Journal of Applied Analysis and Computation},
  volume={7},
  number={4},
  pages={1488--1502},
  year={2017}
}

@incollection{huang2001partition,
  title={A partition-of-unity finite element method for elliptic problems with highly oscillating coefficients},
  author={Huang, Yunqing and Xu, Jinchao},
  editor={Shi, Zhong-Ci and Mu, Mo and Xue, Weimin and Zou, Jun},
  booktitle={Advances in Scientific Computing},
  publisher={Science Press},
  address={Beijing and New York},
  pages={61--70},
  year={2001}
}

@article{efendiev2013generalized,
  title={Generalized multiscale finite element methods (GMsFEM)},
  author={Efendiev, Yalchin and Galvis, Juan and Hou, Thomas Y},
  journal={Journal of Computational Physics},
  volume={251},
  pages={116--135},
  year={2013},
  publisher={Elsevier}
}

@article{calo2016randomized,
  title={Randomized oversampling for generalized multiscale finite element methods},
  author={Calo, Victor M and Efendiev, Yalchin and Galvis, Juan and Li, Guanglian},
  journal={Multiscale Modeling \& Simulation},
  volume={14},
  number={1},
  pages={482--501},
  year={2016},
  publisher={SIAM}
}

@article{alikhanov2025multiscale,
  title={Multiscale model reduction for the time fractional thermoporoelasticity problem in fractured and heterogeneous media},
  author={Alikhanov, Anatoly and Bai, Huiran and Huang, Jian and Tyrylgin, Aleksei and Yang, Yin},
  journal={Journal of Computational and Applied Mathematics},
  volume={455},
  pages={116157},
  year={2025},
  publisher={Elsevier}
}

@article{xie2025time,
  title={On time integrators for Generalized Multiscale Finite Element Methods applied to advection--diffusion in high-contrast multiscale media},
  author={Xie, Wei and Galvis, Juan and Yang, Yin and Huang, Yunqing},
  journal={Journal of Computational and Applied Mathematics},
  volume={460},
  pages={116363},
  year={2025},
  publisher={Elsevier}
}

@article{xie2025multiscale,
  title={Multiscale modeling of wave propagation with exponential integration and GMsFEM},
  author={Xie, Wei and Galvis, Juan and Yang, Yin and Huang, Yunqing},
  journal={Communications in Nonlinear Science and Numerical Simulation},
  volume={147},
  pages={108825},
  year={2025},
  publisher={Elsevier}
}

@article{chung2014adaptive,
  title={An adaptive GMsFEM for high-contrast flow problems},
  author={Chung, Eric T and Efendiev, Yalchin and Li, Guanglian},
  journal={Journal of Computational Physics},
  volume={273},
  pages={54--76},
  year={2014},
  publisher={Elsevier}
}

@article{ye2024robust,
  title={A robust two-level overlapping preconditioner for Darcy flow in high-contrast media},
  author={Ye, Changqing and Fu, Shubin and Chung, Eric T and Huang, Jizu},
  journal={SIAM Journal on Scientific Computing},
  volume={46},
  number={5},
  pages={A3151--A3176},
  year={2024},
  publisher={SIAM}
}

@article{fu2024efficient,
  title={An Efficient Multiscale Preconditioner for Large-Scale Highly Heterogeneous Flow},
  author={Fu, Shubin and Chung, Eric and Zhao, Lina},
  journal={SIAM Journal on Scientific Computing},
  volume={46},
  number={2},
  pages={S352--S377},
  year={2024},
  publisher={SIAM}
}

@article{chung2018constraint,
  title={Constraint energy minimizing generalized multiscale finite element method},
  author={Chung, Eric T and Efendiev, Yalchin and Leung, Wing Tat},
  journal={Computer Methods in Applied Mechanics and Engineering},
  volume={339},
  pages={298--319},
  year={2018},
  publisher={Elsevier}
}

@article{ye2023constraint,
  title={Constraint energy minimizing generalized multiscale finite element method for inhomogeneous boundary value problems with high contrast coefficients},
  author={Ye, Changqing and Chung, Eric T},
  journal={Multiscale Modeling \& Simulation},
  volume={21},
  number={1},
  pages={194--217},
  year={2023},
  publisher={SIAM}
}

@article{wang2024multiscale,
  title={A multiscale method for inhomogeneous elastic problems with high contrast coefficients},
  author={Wang, Zhongqian and Ye, Changqing and Chung, Eric T},
  journal={Journal of Computational and Applied Mathematics},
  volume={436},
  pages={115397},
  year={2024},
  publisher={Elsevier}
}

@article{xie2024cem,
  title={{CEM-GMsFEM for Poisson equations in heterogeneous perforated domains}},
  author={Xie, Wei and Yang, Yin and Chung, Eric and Huang, Yunqing},
  journal={Multiscale Modeling \& Simulation},
  volume={22},
  number={4},
  pages={1683-1708},
  year={2024},
  publisher={SIAM}
}

@article{chung2018non,
  title={Non-local multi-continua upscaling for flows in heterogeneous fractured media},
  author={Chung, Eric T and Efendiev, Yalchin and Leung, Wing Tat and Vasilyeva, Maria and Wang, Yating},
  journal={Journal of Computational Physics},
  volume={372},
  pages={22--34},
  year={2018},
  publisher={Elsevier}
}

@article{chung2021nonlocal,
  title={Nonlocal multicontinua with representative volume elements. Bridging separable and non-separable scales},
  author={Chung, Eric T and Efendiev, Yalchin and Leung, Wing T and Vasilyeva, Maria},
  journal={Computer Methods in Applied Mechanics and Engineering},
  volume={377},
  pages={113687},
  year={2021},
  publisher={Elsevier}
}

@article{abdulle2012heterogeneous,
  title={The heterogeneous multiscale method},
  author={Abdulle, Assyr and Weinan, E and Engquist, Bj{\"o}rn and Vanden-Eijnden, Eric},
  journal={Acta Numerica},
  volume={21},
  pages={1--87},
  year={2012},
  publisher={Cambridge University Press}
}

@article{henning2009heterogeneous,
  title={The heterogeneous multiscale finite element method for elliptic homogenization problems in perforated domains},
  author={Henning, Patrick and Ohlberger, Mario},
  journal={Numerische Mathematik},
  volume={113},
  pages={601--629},
  year={2009},
  publisher={Springer}
}

@article{maalqvist2014localization,
  title={Localization of elliptic multiscale problems},
  author={M{\aa}lqvist, Axel and Peterseim, Daniel},
  journal={Mathematics of Computation},
  volume={83},
  number={290},
  pages={2583--2603},
  year={2014}
}

@article{henning2014localized,
  title={Localized orthogonal decomposition techniques for boundary value problems},
  author={Henning, Patrick and M{\aa}lqvist, Axel},
  journal={SIAM Journal on Scientific Computing},
  volume={36},
  number={4},
  pages={A1609--A1634},
  year={2014},
  publisher={SIAM}
}

@article{chung2024multicontinuum,
  title={Multicontinuum homogenization. General theory and applications},
  author={Chung, E and Efendiev, Yalchin and Galvis, Juan and Leung, Wing Tat},
  journal={Journal of Computational Physics},
  volume={510},
  pages={112980},
  year={2024},
  publisher={Elsevier}
}

@article{efendiev2023multicontinuum,
  title={Multicontinuum homogenization and its relation to nonlocal multicontinuum theories},
  author={Efendiev, Yalchin and Leung, Wing Tat},
  journal={Journal of Computational Physics},
  volume={474},
  pages={111761},
  year={2023},
  publisher={Elsevier}
}

@article{xie2025multicontinuum,
  title={Multicontinuum homogenization in perforated domains},
  author={Xie, Wei and Efendiev, Yalchin and Huang, Yunqing and Leung, Wing Tat and Yang, Yin},
  journal={Journal of Computational Physics},
  pages={113845},
  year={2025},
  publisher={Elsevier}
}

@article{efendiev2025multicontinuum,
  title={Multicontinuum splitting scheme for multiscale flow problems},
  author={Efendiev, Yalchin and Leung, Wing Tat and Shan, Buzheng and Wang, Min},
  journal={Journal of Computational Physics},
  pages={114265},
  year={2025},
  publisher={Elsevier}
}

@article{bai2026multicontinuum,
  title={Multicontinuum modeling of time-fractional diffusion-wave equation in heterogeneous media},
  author={Bai, Huiran and Ammosov, Dmitry and Yang, Yin and Xie, Wei and Al Kobaisi, Mohammed},
  journal={Journal of Computational and Applied Mathematics},
  volume={473},
  pages={116846},
  year={2026},
  publisher={Elsevier}
}

@article{park2020hierarchical,
  title={Hierarchical multiscale finite element method for multi-continuum media},
  author={Park, Jun Sur Richard and Hoang, Viet Ha},
  journal={Journal of Computational and Applied Mathematics},
  volume={369},
  pages={112588},
  year={2020},
  publisher={Elsevier}
}

@article{brown2013efficient,
  title={An efficient hierarchical multiscale finite element method for Stokes equations in slowly varying media},
  author={Brown, Donald L and Efendiev, Yalchin and Hoang, Viet Ha},
  journal={Multiscale Modeling \& Simulation},
  volume={11},
  number={1},
  pages={30--58},
  year={2013},
  publisher={SIAM}
}

@article{brown2017hierarchical,
  title={A hierarchical finite element Monte Carlo method for stochastic two-scale elliptic equations},
  author={Brown, Donald L and Hoang, Viet Ha},
  journal={Journal of Computational and Applied Mathematics},
  volume={323},
  pages={16--35},
  year={2017},
  publisher={Elsevier}
}

@book{briggs2000multigrid,
  title={A multigrid tutorial},
  author={Briggs, William L and Henson, Van Emden and McCormick, Steve F},
  year={2000},
  publisher={SIAM}
}

@article{xu2017algebraic,
  title={Algebraic multigrid methods},
  author={Xu, Jinchao and Zikatanov, Ludmil},
  journal={Acta Numerica},
  volume={26},
  pages={591--721},
  year={2017},
  publisher={Cambridge University Press}
}

@article{cheng2022scalable,
  title={Scalable semismooth Newton methods with multilevel domain decomposition for subsurface flow and reactive transport in porous media},
  author={Cheng, Tianpei and Yang, Haijian and Yang, Chao and Sun, Shuyu},
  journal={Journal of Computational Physics},
  volume={467},
  pages={111440},
  year={2022},
  publisher={Elsevier}
}

@article{efendiev2015multilevelmms,
  title={Multilevel Monte Carlo approaches for numerical homogenization},
  author={Efendiev, Yalchin and Kronsbein, Cornelia and Legoll, Fr{\'e}d{\'e}ric},
  journal={Multiscale Modeling \& Simulation},
  volume={13},
  number={4},
  pages={1107--1135},
  year={2015},
  publisher={SIAM}
}

@article{efendiev2015multilevelcicp,
  title={Multilevel Markov chain Monte Carlo method for high-contrast single-phase flow problems},
  author={Efendiev, Yalchin and Jin, Bangti and Michael, Presho and Tan, Xiaosi},
  journal={Communications in Computational Physics},
  volume={17},
  number={1},
  pages={259--286},
  year={2015},
  publisher={Cambridge University Press}
}

@article{li2018multilevel,
  title={Multilevel space-time additive Schwarz methods for parabolic equations},
  author={Li, Shishun and Shao, Xinping and Cai, Xiao-Chuan},
  journal={SIAM Journal on Scientific Computing},
  volume={40},
  number={5},
  pages={A3012--A3037},
  year={2018},
  publisher={SIAM}
}

@article{ye2025highly,
  title={A highly parallelized multiscale preconditioner for Darcy flow in high-contrast media},
  author={Ye, Changqing and Fu, Shubin and Chung, Eric T and Huang, Jizu},
  journal={Journal of Computational Physics},
  volume={522},
  pages={113603},
  year={2025},
  publisher={Elsevier}
}
\end{document}